\documentclass[a4paper,11pt,reqno]{amsart}
\usepackage[centering,includeheadfoot,margin=2.5 cm]{geometry}
\usepackage{graphics,enumitem,epsfig,textcomp}
\usepackage{amsfonts,amsmath,amssymb,euscript,color,mathrsfs}
\usepackage[utf8]{inputenc}	

\usepackage[caption=false,subrefformat=parens,labelformat=parens]{subfig}

\usepackage{url}

\usepackage[symbol]{footmisc}

\numberwithin{equation}{section}
\newtheorem{theorem}{Theorem}%[section]

\newtheorem{lemma}{Lemma}%[section]
\newtheorem{remark}{Remark}%[section]

\def\d{\,\mathrm{d}}
\def\eps{\varepsilon}
\def\N{\mathbb{N}}
\def\R{\mathbb{R}}
\def\C{\hbox{\rlap{\kern.24em\raise.1ex\hbox
			{\vrule height1.3ex width.9pt}}C}}
\def\P{\hbox{\rlap{I}\kern.16em P}}
\def\Q{\hbox{\rlap{\kern.24em\raise.1ex\hbox
			{\vrule height1.3ex width.9pt}}Q}}
\def\M{\hbox{\rlap{I}\kern.16em\rlap{I}M}}
\def\Z{\hbox{\rlap{Z}\kern.20em Z}}
\def\({\begin{eqnarray}}
\def\){\end{eqnarray}}
\def\[{\begin{eqnarray*}}
\def\]{\end{eqnarray*}}
\def\part#1#2{\frac{\partial #1}{\partial #2}}

\def\grad{\nabla}
\def\Norm#1{\left\| #1 \right\|}
\def\pmb#1{\setbox0=\hbox{$#1$}
	\kern-.025em\copy0\kern-\wd0
	\kern-.05em\copy0\kern-\wd0
	\kern-.025em\raise.0433em\box0 }
\def\bar{\overline}

\def\tot#1#2{\frac{\d #1}{\d #2}} 
\def\laplace{\Delta}
\def\d{\,\mathrm{d}}
\def\N{\mathbb{N}}
\def\R{\mathbb{R}}

\def\epsilon{\varepsilon}

\def\E{\mathcal{E}}
\def\P{\mathbb{P}}
\def\Q{\mathbb{Q}}

\newcommand*\di{\mathop{}\!\mathrm{d}}
\newcommand\bl{\left(}
\newcommand\br{\right)}

\newcommand\Vset{V}
\newcommand\Eset{E}

\newcommand\Neigh{\mathcal{N}}

\def\review#1{#1}

\begin{document}
	%%%%%%%%%%%%%%%%

	%%%%%%%%%%%%%%%%%%%%%%%%%%%%%%%%%
	%% Title
	%%%%%%%%%%%%%%%%%%%%%%%%%%%%%%%%%
	%\title[]{Discrete modelling of biological transport networks}
	
	\centerline{{\Large\textbf{Auxin transport model for leaf venation}}}
	\vskip 7mm

	%%%%%%%%%%%%%%%%%%%%%%%%%%%%%%%%%
	%% Authors
	%%%%%%%%%%%%%%%%%%%%%%%%%%%%%%%%%
	\centerline{
		{\large Jan Haskovec}\footnote{Mathematical and Computer Sciences and Engineering Division,
			King Abdullah University of Science and Technology,
			Thuwal 23955-6900, Kingdom of Saudi Arabia; 
			{\it jan.haskovec@kaust.edu.sa}}\qquad
		{\large Henrik J\"onsson}\footnote{Sainsbury Laboratory, University of Cambridge,
			Bateman Street, Cambridge CB2 1LR, UK; Department of Applied Mathematics and Theoretical Physics (DAMTP), University of Cambridge, Wilberforce Road, Cambridge CB3 0WA, UK;
			{\it henrik.jonsson@slcu.cam.ac.uk}}\qquad
		{\large Lisa Maria Kreusser}\footnote{Department of Applied Mathematics and Theoretical Physics (DAMTP), University of Cambridge, Wilberforce Road, Cambridge CB3 0WA, UK;
			{\it L.M.Kreusser@damtp.cam.ac.uk}}\qquad
		{\large Peter Markowich}\footnote{Mathematical and Computer Sciences and Engineering Division,
			King Abdullah University of Science and Technology,
			Thuwal 23955-6900, Kingdom of Saudi Arabia;
			Faculty of Mathematics, University of Vienna, Oskar-Morgenstern-Platz 1, 1090 Vienna, Austria;
			{\it peter.markowich@kaust.edu.sa; peter.markowich@univie.ac.at}}
	}
	\vskip 10mm

	\noindent{\bf Abstract.} 
The plant hormone auxin controls many aspects of the development of plants.
One striking dynamical feature is the self-organisation of leaf
venation patterns which is driven by high levels of auxin within vein cells.
The auxin transport is mediated by specialised
membrane-localised proteins. Many venation models have been based on polarly localised efflux-mediator
proteins of the PIN family.
%Since PINs have an a priori defined direction, PIN-based models consider polar/directional transport activity where the orientation of the flow is given.
%However, pattern formation may also arise without PINs. }
%Highly important for correct
%venation patterns are polarly localised efflux-mediator
%proteins of the PIN family.
%, expressed in the cells during
%vein initiation and contributing to a generation of directed
%auxin fluxes within the veins. Still, veins form also when
%PIN-mediated transport is blocked chemically or in mutants,
%and the mechanism driving vein formation is not fully
%understood.
Here, we investigate a modeling framework for auxin transport with a
positive feedback between auxin fluxes and transport
capacities that are not necessarily polar, i.e.\ directional across a cell wall.
Our approach is derived from a discrete graph-based model
for biological transportation networks,
where cells are represented by graph nodes
and intercellular membranes by edges.
The edges are not a-priori oriented and the direction of auxin flow
is determined by its concentration gradient along the edge.
We prove global existence of solutions to the model and the validity of Murray's law for its steady states.
Moreover, we demonstrate with numerical simulations that the model is able connect an auxin
source-sink pair with a mid-vein and that it can also produce
branching vein patterns.
A significant innovative aspect of our approach is that it allows the
passage to a formal macroscopic limit which can be extended to include network growth.
We perform mathematical analysis of the macroscopic formulation, showing the global existence
of weak solutions for an appropriate parameter range.
	
	\vskip 7mm
	
	%\allowdisplaybreaks
	
	\subjclass[]{}
	
	%\thanks{\textbf{Acknowledgments.} This work is supported by Engineering and Physical Sciences Research Council(EP/K008404/1).}
	
	%\maketitle \centerline{\date}

	%\tableofcontents

%%%%%%%%%%%%%%%%%%%%%%%%%%%%%%%%%
%% Introduction
%%%%%%%%%%%%%%%%%%%%%%%%%%%%%%%%%
\section{Introduction}

The hormone auxin plays a central role in \review{many} developmental
processes in plants
\cite{Heisler2010,Scarpella2006,Sawchuk2013,Sawchuk2013:2}. During the development of
a leaf, a connected network of veins is formed in a highly
predictable order, generating a well defined pattern in the final leaf
\cite{hickey1973}. 
\review{High levels of auxin are present}  in the  forming
vein cells compared to the neighboring 
tissues. It has been shown that the membrane localized
PIN-FORMED (PIN) family of auxin transport mediators 
is essential for the correct patterning of the vein network
\cite{Sachs1969, Scarpella2006}. \review{The
	patterns} could result from a canalisation mechanism where the auxin
flux feeds back itself to a polarised transport connecting sources and
sinks of auxin \cite{Sachs1981,Mitchison1980,Mitchison1981}. This idea has been revisited recently and has led to
models with polarised PIN transporters
\cite{RollandLagan2005,Feugier2005,Feller2015}. No flux-sensing
mechanism has been identified 
but models have been used to suggest
alternatives \cite{Kramer2009, Cieslak2015}. While newer models have
solved the \review{issue of unrealistically low levels of auxin within veins in flux-based models}
% having low auxin
%within the veins
\cite{Feugier2005}, \review{it is still an open question} how
looped veins can form \cite{RollandLagan2005,Dimitrov2006} \review{and if}
specified auxin production can provide \review{an answer}.

PIN proteins are involved in several patterning processes in plants. Alternative models, not based on auxin flux, have been proposed, for
\review{instance for producing Turing-like dynamics} in the context of
phyllotaxis \cite{Jonsson2006, Smith2006, Bhatia2016}, and for single
cell polarity resulting in planar polarity \cite{Abley2016}.

\review{Since the discovery of PINs,  many
	venation	models have been based}  on polarised transport via PINs, \review{while recent data suggests} that polar auxin transport mediated by PINs
is not crucial for forming veins \cite{Sawchuk2013,Sawchuk2013:2}. Although \review{
	characteristic vein patterns and leaf shapes can be obtained with these PIN-based models, veins can also form in chemical
	perturbations}
when PIN-mediated auxin transport is blocked, or when multiple membrane-localised PIN proteins are
mutated. This raises the question if alternative mechanisms 
work in parallel or together with the PIN-based polar transporters
during the initiation of veins. This motivates to consider a more general modelling approach where alternative feedbacks between auxin, auxin fluxes and auxin transport can be included.
%This motivates to consider  a more general modelling approach where PINs are interpreted as heterogeneous auxin transporters.

The ultimate goal for modelling vein networks is to accurately predict vein network geometries seen in different plants. Our novel dynamical description could complement the PIN-based models which have focused on more basic dynamic patterns of veins, such as connecting sources and sinks, and breaking the symmetry of graded diffusion into veins. \review{Examples of these PIN-based models include}  the traditional PIN-based flux models that have been studied since approximately 40 years, see \cite{Sachs1981,Mitchison1980,Mitchison1981}.
%We show some similarities in crude vein formation, but  more elaborate investigations are essential, such as combining the current description with the PIN-based mechanisms
%and testing with more complex configurations of auxin sources and sinks.
The impact of auxin concentration \review{on} the pattern formation \review{has been studied} in \cite{McAdam351}. \review{It would be very interesting to investigate the emergence of patterns in the setting where PINs are removed}. As noted above, the traditional PIN-based flux models are yet to provide a full description of the diverse patterns seen in plants.

Given the strong directional distribution of PINs  and the ability of veins to form \review{without PINs}, it is  \review{important} to introduce and analyse alternative mechanisms. Whether these \review{mechanisms are} identical/redundant to PIN mechanisms in terms of \review{their} dynamical behaviour \review{or whether other mechanisms need to be considered}
%, or alternative mechanisms
is still unknown. Hence, \review{it would be interesting to show} that  polar/directional transport activity and  \review{directional flux measurements}  are not required, and \review{that} vein-like patterns can also result from mere measurement of  magnitudes. \review{This may also} inspire scientists to \review{reconsider} their current data or \review{design new}  experiments. 

In this paper we study a modeling framework for leaf venation which does
not assume polarity of auxin transport mediators across cell walls.
%In particular, we assume that the transport capacity across cell walls does not have to be directional.
The model is introduced in Section \ref{sec:desc}, and
is based on a positive feedback loop between auxin fluxes and transport
capacities that are not necessarily polar.
Our approach is derived from a recent discrete graph-based model
for biological transportation networks introduced by Hu and Cai \cite{hu}.
We represent cells by graph nodes
and intercellular membranes (connections) by edges.
The edges are not a-priori oriented and the direction of auxin flow
is determined by its concentration gradient along the edge.
The transport capacity of each edge is represented by the
local concentration of the auxin mediator.
Our approach can be understood as a modeling framework, which can be equipped or extended
with various biologically relevant features that will produce experimentally testable hypotheses.
\review{We admit that in its present setting it does not capture all relevant biological features,}
however, its main advantage is a rather simple form that facilitates rigorous mathematical analysis.
In particular, the first aim of this paper is the proof of global existence and nonnegativity of solutions
of the discrete model (Section \ref{sec:AnalysisD}).
Moreover, in Section \ref{sec:Murray} we show that the stationary solutions \review{satisfy} a generalized
Murray's law.
The second aim of the paper is to gain \review{a better} understanding of the pattern formation capacity
of the model by means of numerical simulations (Section \ref{sec:num}).
In particular, we show that it is capable of generating  patterns
connecting an auxin source-sink pair with a mid-vein and that it can  produce
branching vein patterns.
The main novelty of our modelling approach is that it facilitates a (formal) passage to a continuum limit,
which is the subject of Section \ref{sec:continuum}.
The resulting system of partial differential equations captures network growth and is expected to exhibit a rich patterning capacity
(see \cite{bionetworksbookchapter} for results of numerical simulations of a related continuum model).
Here we prove the existence of weak solutions of the transient problem and of its steady states.

%%%%%%%%%%%%%%%%%%%%%%%%%%%%%%%%%%%%%
\section{Description of the model}\label{sec:desc}

Hu and Cai considered a discrete model describing the formation of \review{generic} biological transport networks in \cite{hu}.
Existence of transient solutions, their qualitative properties and the formal continuum limit of the Hu and Cai model was studied in \cite{HaKrMa}.
Here we adapt the model to the cellular context to describe auxin transport in plant leafs via transporter proteins,
where the orientation of the flow is determined by auxin concentration gradient.
Our approach shares many similarities with the one introduced by Mitchison in \cite{Mitchison1980}
where the transport capacity is updated as a function of the flux (gradient) between cells. However, while Mitchison suggested an asymmetric update of the transport capacities across a cell wall, our model  assumes a symmetric transport capacity across a cell wall.
In this section we shall first introduce the Hu and Cai model, then shortly discuss the Mitchison model,
and finally describe the adaptation to the cellular context.
%\review{In particular, the orientation of the flow is not a priori given in our modelling approach, while most current plant models consider polar/directional transport via PINs.}
% an asymmetry in transport in one direction versus the opposite is usually considered. 

\subsection{Model of Hu and Cai \cite{hu}}
The discrete model introduced by Hu and Cai \cite{hu} and reformulated in \cite{bionetworksbookchapter}
is posed on a given, fixed undirected connected graph $G=(\Vset,\Eset)$, consisting of a finite set of vertices $\Vset$ of size $N=|\Vset|$ and a finite set of edges $\Eset$.
Any pair of vertices is connected by at most one edge and no vertex is connected to itself.
We denote the edge  between vertices $i\in\Vset$ and $j\in\Vset$ by $(i,j)\in \Eset$.
Since the graph is undirected, $(i,j)$ and $(j,i)$ refer to the same edge.
For each edge $(i,j)\in\Eset$ of the graph $G$ we consider its length and its conductivity, denoted by  $L_{ij}=L_{ji}>0$ and $C_{ij}=C_{ji}\geq 0$, respectively. 
%We denote the conductivity vector by $C=(C_{ij})_{(i,j)\in\Eset}$. 
The edge lengths $L_{ij}>0$ are given as a datum and fixed for all $(i,j)\in\Eset$.
With each vertex $i\in\Vset$, the fluid pressure $P_i\in\R$ is associated.
The pressure drop between vertices $i\in\Vset$ and $j\in\Vset$ connected by an edge $(i,j)\in\Eset$ is given by 
\begin{align}\label{eq:pressuredrop}
(\Delta P)_{ij}:=P_j-P_i.
\end{align}
Note that the pressure drop is antisymmetric, i.e., by definition,  $(\Delta P)_{ij}=-(\Delta P)_{ji}$.
The oriented flux (flow rate) from vertex $i\in\Vset$ to $j\in\Vset$ is denoted by $Q_{ij}$; again, we have $Q_{ij}=-Q_{ji}$.
Since  the Reynolds number of the flow is typically small for biological networks and the flow is predominantly laminar,
the flow rate between vertices $i\in\Vset$ and $j\in\Vset$ along edge $(i,j)\in\Eset$ is proportional to the conductance $C_{ij}$ and the pressure drop $(\Delta P)_{ij}=P_j-P_i$, 
\begin{align}\label{eq:flowrate}
Q_{ij} := C_{ij}\frac{P_j-P_i}{L_{ij}}\qquad\text{for all~}(i,j)\in \Eset.
\end{align}
The local mass conservation in each vertex is expressed in terms of the Kirchhoff law
\begin{align}\label{eq:kirchhoff}
-\sum_{j\in \Neigh(i)} C_{ij}\frac{ P_j-P_i}{L_{ij}}=S_i\qquad \text{for all~}i\in \Vset.
\end{align}
Here $\Neigh(i)$ denotes the set of vertices connected to $i\in\Vset$ through an edge,
and $S=(S_i)_{i\in\Vset}$ is the prescribed strength of the flow source ($S_i>0$) or sink ($S_i<0$) at vertex $i$.
Clearly, a necessary condition for the solvability of \eqref{eq:kirchhoff} is the global mass conservation
\begin{align}\label{eq:conservationmass}
\sum_{i\in\Vset}S_i=0,
\end{align}
which we assume in the following.
Given the vector of conductivities $C=(C_{ij})_{(i,j)\in\Eset}$, the Kirchhoff law \eqref{eq:kirchhoff} is a linear system of equations for the vector of pressures $P=(P_i)_{i\in\Vset}$.
With the global mass conservation \eqref{eq:conservationmass}, the linear system \eqref{eq:kirchhoff} is solvable if and only if the graph with edge weights $C=(C_{ij})_{(i,j)\in\Eset}$ is connected \cite{bionetworksbookchapter}, where only edges with positive conductivities $C_{ij}>0$ are taken into account
(i.e., edges with zero conductivities are discarded).
Note that the solution is unique up to an additive constant.

The conductivities $C_{ij}$ are subject to an energy optimization and adaptation process. %We assume that initially all edges in $\Eset$ have strictly positive conductivities.
Hu and Cai \cite{hu} propose an energy cost functional consisting of a pumping power term and a metabolic cost term.
According to  Joule's law, the power (kinetic energy) needed to pump material through an edge $(i,j)\in\Eset$ is proportional to the pressure drop $(\Delta P)_{ij}=P_j-P_i$
and the flow rate $Q_{ij}$ along the edge, i.e.,
%\begin{align*}
$(\Delta P)_{ij}Q_{ij}=\frac{Q_{ij}^2}{C_{ij}}L_{ij}$. 
%=C_{ij}\bl\frac{P_j-P_i}{L_{ij}}\br^2 L_{ij}.
%\end{align*}
The metabolic cost of maintaining the edge is assumed to be proportional to its length $L_{ij}$ and a power of its conductivity $C_{ij}^{\gamma}$, where the exponent $\gamma>0$ depends on the network.
%For instance, in blood vessels the metabolic cost is proportional to the cross-section area of the vessel \cite{motivationblood}. Modeling the blood flow by Hagen-Poiseuille's law,
%the conductivity is proportional to the square of the cross-section area, implying $\gamma=1/2$ for blood vessel systems.
For models of leaf venation the material cost is proportional to the number of small tubes, which is
proportional to $C_{ij}$, and the metabolic cost is due to the effective loss of the photosynthetic power
at the area of the venation cells, which is proportional to $C_{ij}^{1/2}$. Consequently, the effective value
of $\gamma$ typically used in models of leaf venation lies between $1/2$ and $1$; see \cite{hu}.
The energy cost functional is thus given by
\begin{align}\label{eq:energydisc}
{\E}[C] := \sum_{(i,j)\in\Eset}\bl \frac{Q_{ij}[C]^2}{C_{ij}}+\frac{\nu}{\gamma} C_{ij}^{\gamma}\br L_{ij},
\end{align}
where $Q_{ij}[C]$ is given by \eqref{eq:flowrate} with pressures calculated from the Kirchhoff's law \eqref{eq:kirchhoff},
and $\nu>0$ is the so-called metabolic coefficient.
Note that  every edge of the graph $G$ is counted exactly once in the above sum.
Hu and Cai \cite{hu}	propose an energy optimization and adaptation process for the conductivities $C_{ij}$
based on the gradient flow of the energy \eqref{eq:energydisc},
\begin{align}\label{eq:GF-orig}
\tot{C_{ij}}{t} = \sigma\left( \frac{Q_{ij}[C]^2}{C_{ij}^{\gamma+1}} - \tau^2 \right) C_{ij} L_{ij}
\end{align}
with parameters $\sigma,\tau>0$, constrained by the Kirchhoff law \eqref{eq:kirchhoff},
see \cite{HaKrMa} for details.
%	Note that for $\alpha=2-\gamma$ the ODE system  \eqref{eq:GF-orig} can be rewritten as the gradient flow
%	\begin{align*}%\label{eq:GF-alpha}
%	\tot{C_{ij}}{t} = \left( \frac{Q_{ij}[C]^2}{C_{ij}} - \nu C_{ij}^{\gamma} \right) C_{ij}^{\alpha-1} L_{ij}
%	\end{align*}

\subsection{Mitchison model \cite{Mitchison1981}}

%As described in the introduction, auxins are a class of plant hormones (or plant growth regulators) that play a cardinal role
%in coordination of many growth and behavioral processes in the plant's life cycle
%and are essential for plant body development including for developing its own transport network.
The model proposed by Mitchison \cite{Mitchison1980} describes auxin dynamics within an array of cells with indices $i\in\Vset$.
For two cells $i,j\in\Vset$ with signal concentrations $s_i,s_j$, respectively, the diffusion constant at the interface between the cells is denoted by $D_{ij}=D_{ji}\geq 0$
and can be specified independently for each cell-cell interface. The oriented flux from vertex $i\in \Vset$ to $j\in \Vset$  is given by Fick's law \cite{crank}, 
\begin{align}
\phi_{ij}=D_{ij} \frac{s_i-s_j}{L_{ij}},
\end{align} 
where $L_{ij}=L_{ji}>0$ denotes the (average) length of cells $i$ and  $j$. In particular, we have %the antisymmetry
$\phi_{ij}=-\phi_{ji}$.
The dependence of the diffusion constant $D_{ij}$ on the flux $\phi_{ij}$ is of the form
\begin{align*}
\frac{\di D_{ij}}{\di t}= f(|\phi_{ij}|,D_{ij})
\end{align*}
for a suitable function $f$ such that $|\phi_{ij}|/D_{ij}$  decreases as $|\phi_{ij}|$ increases. For instance, $f$ can be chosen such that $D_{ij}\approx \phi_{ij}^2$ for $\phi_{ij}>0$ and $D_{ij} = 0$ for $\phi_{ij}\leq 0$, resulting in a strictly polar transport capacity across a cell wall.
%at least in a neighborhood of $f^{-1}(0)$. 
Assuming that cell $i\in \Vset$ receives fluxes $\phi_{ji}$ for $j\in \Neigh(i)$, the evolution of the signal $s_i$ is of the form 
\begin{align}\label{eq:signaleq}
\frac{\di s_i}{\di t} =\sigma_i+\frac{1}{v}\sum_{j\in \Neigh(i)} A_{ij}\phi_{ji}.
\end{align}
As before, $\Neigh(i)$ denotes the index set of neighboring cells of cell $i\in\Vset$. The parameter $\sigma_i$ is the source activity for signal production in  cell $i\in\Vset$.
All cells have volume $v>0$ and $A_{ij}=A_{ji}>0$ is the area of the interface between cell $i$ and its neighbor $j\in \Neigh(i)$. Note that the term $\sum_{j\in \Neigh(i)} A_{ij}\phi_{ji}$ can be regarded as the difference between influx and outflux since $\phi_{ij}=-\phi_{ji}$ for $j\in\mathcal{N}(i)$. For the conservation of the signal we require that the source activity $\sigma_i$ for signal production and degradation is chosen such that
%	balances over all cells $i\in\Vset$, i.e.
%	\begin{align*}
%	\sum_{i\in\Vset} \sigma_i=0,
%	\end{align*}	
%	leading to
%\begin{align*}
$	\frac{\di }{\di t}\sum_{i\in\Vset} s_i=0.$
%\end{align*}

It is worth noting that while it was
well established that auxin was important for generating the
vascular or vein patterns (see, e.g., \cite{Sachs1981}), auxin
`transporters' were not identified at the time when the
model was introduced. It received great attention only
later, when auxin transport mediator proteins with similar
polar localisation as predicted by the model were
identified \cite{Scarpella2006}. In particular, PIN proteins are
integral membrane proteins that transport the anionic form
of auxin across membranes. Most of the PIN proteins localize
at the plasma membrane where they serve as secondary active
transporters involved in the efflux of auxin. They show
asymmetrical localizations on the membrane and are therefore
responsible for polar auxin transport. Still, while PIN loss of function mutants generate phenotypes in
venation patterns, they do not completely abolish the
formation of veins \cite{Sawchuk2013}, and as such
alternative mechanisms can contribute to the dynamics of
vein formation. While individual mutants do not show strong
phenotypes, this is also implied by the existence of other auxin
transport proteins, such as AUX1/LAX influx mediators          
\cite{Kramer2004, Peret2012, Sawchuk2013}, regulating
intracellular and intercellular transport. 
In the following discussion we will often use PIN as a descriptor of the auxin transporter protein for simplicity, but it should be seen as a more general description of auxin transport mediated by polar and/or nonpolar membrane proteins, where polar relates  to the difference of transport capacity (PIN localisation) on the two sides of a wall.

\subsection{Adapted Hu-Cai model in cellular context}
Given the known auxin flows generated from sources to sinks in a plant tissue, the sometimes clear expression but unclear polarisation of PIN auxin transporter proteins in these veins, and the ability to generate veins without any PIN transport, it is of interest to investigate alternative mechanisms for the vein dynamics in an auxin context. Such an alternative can be provided
by a proper adaptation of the Hu and Cai model for transport networks \cite{hu}. The mechanism where pressure differences feed back on conductance between elements has similarity with the auxin transport case, as described in the flux-based models \cite{Mitchison1980,Mitchison1981}.
Here auxin sources and concentration differences (pressure in the Hu-Cai model) generate diffusive fluxes between cells (spatial elements) that positively feed back on transport rates between the cells (conductance).
To adapt the Hu-Cai model to a cellular context of plant venation dynamics we consider $n=|V|$ cells with indices $i\in\Vset$ and replace the pressure $P_i$ at vertex $i\in\Vset$ in the Hu-Cai model with the auxin concentration $a_i\geq 0$.

%	to the pressure drop, the auxin drop between cells $i\in\Vset$ and $j\in\Vset$ is defined as
%	\begin{align*}
%		(\Delta a)_{ij}:=a_j-a_i.
%	\end{align*}
%	Like the pressure drop $(\Delta P)_{ij}$ in the original model the auxin drop is antisymmetric, i.e., $(\Delta a)_{ij}=-(\Delta a)_{ji}$.
%We say two cells $i,j\in\Vset$ are neighboring cells and connected to each other if $(i,j)\in\Eset$.
The conductance $C_{ij}$ of edge $(i,j)\in\Eset$ in the Hu and Cai model is replaced by the transport activity $X_{ij}=X_{ji}\geq 0$ in the membrane connecting cells $i\in\Vset$ and $j\in\Vset$ which is the main difference from PIN-based flux models (and experiments) with PINs $\mathcal P_{ij}$ where $\mathcal P_{ij} \neq \mathcal P_{ji}$. Due to this modelling approach auxin transporters are not directional, i.e.\ polar,
and as we shall see, measuring the magnitudes $X_{ij}$ is sufficient for  producing vein-like dynamics.
However, cells, in general, do not transport auxin equally well in all directions (i.e. $X_{ij}$ is typically not equal to $X_{ik}$ for two cell neighbours $i$ and $k$). Based on the definition of $X_{ij}$, 
we define the auxin flow rate $\mathcal{Q}_{ij}=-\mathcal{Q}_{ji}\in\R$ from cell $i\in\Vset$ to cell $j\in\Vset$ by 
%\begin{align*}
$\mathcal{Q}_{ij}=X_{ij}\frac{a_j-a_i}{L_{ij}},$
%\end{align*}
where $L_{ij}=L_{ji}>0$ denotes the (average) length of cells $i$ and  $j$.
Based on the frameworks of Mitchison \eqref{eq:signaleq} and Hu and Cai \eqref{eq:GF-orig} we describe the auxin transport in the cellular context by
the ODE system
\begin{align}
\tot{a_i}{t} %&= S_i - I_ia_i +\delta\sum_{j\in \Neigh(i)} \mathcal{Q}_{ij} \nonumber \\
&= S_i - I_ia_i+\delta\sum_{j\in \Neigh(i)} X_{ij}\frac{a_j-a_i}{L_{ij}}\qquad \text{for all~}i\in \Vset, \label{eq:auxineqproduction}
\end{align}
where $\Neigh(i)$ denotes the index set of neighboring cells of cell $i\in\Vset$ and the parameter $\delta>0$ denotes the (scaled) diffusion rate.
To account for the auxin production and destruction in the cells, we introduced the source terms $S_i\geq 0$ and decay rates $I_i\geq 0$ for $i\in\Vset$.
For simplicity, we assume $S_i$ and $I_i$ to be independent of time.
For the transport activity $X_{ij}$ in the membrane we consider %the evolution
\begin{align}\label{eq:pineq}
\tot{X_{ij}}{t} = \sigma\left( \frac{|\mathcal{Q}_{ij}|^{\kappa}}{X_{ij}^{\gamma+1}} - \tau \right) X_{ij} L_{ij},
\end{align}
where $\gamma>0$ is a control parameter and $\sigma$, $\kappa$, $\tau$ are nonnegative parameters denoting, respectively,
the conductance update rate, the flux feedback and the conductance degradation rate.
In particular, the flux feedback $\kappa$ is an important parameter of the model and is also a relevant parameter in the Mitchison model \cite{Mitchison1980,Mitchison1981}.
The system \eqref{eq:auxineqproduction}--\eqref{eq:pineq} is equipped with the initial datum
\begin{align}
X_{ij}(0) &= X_{ij}^0 =  X_{ji}^0 \geq 0\qquad \text{for all~}i, j \in \Vset,   \label{pin-IC} \\
a_{i}(0) &= a_{i}^0 > 0\qquad \text{for all~}i \in \Vset.   \label{aux-IC}
\end{align}
Clearly,  \eqref{eq:pineq} satisfies the symmetry requirement $X_{ij}=X_{ji}$.
The conductance equation \eqref{eq:GF-orig} and the transport activity equation \eqref{eq:pineq} are of similar form. However, the term $\mathcal{Q}_{ij}^2$
in the conductance equation \eqref{eq:GF-orig} is replaced by the more general term $|\mathcal{Q}_{ij}|^{\kappa}$ in the transport activity equation \eqref{eq:pineq}
so that  \eqref{eq:pineq} reduces to  \eqref{eq:GF-orig} for $\kappa=2$.
Besides, the linear algebraic system \eqref{eq:kirchhoff} is relaxed by the introduction of the time derivative of the auxin concentration in \eqref{aux-rescaled},
leading to a system of linear ordinary differential equations.
%Note the resemblance (up to rescaling) of \eqref{eq:auxineq}--\eqref{eq:pineq} to the model \eqref{eq:GF-orig}--\eqref{eq:kirchhoff}.
%Moreover, due to the symmetry $\mathcal{P}_{ij}=\mathcal{P}_{ji}$ we obtain the conservation of total auxin concentration,
%\begin{align*}
%\tot{}{t}\sum_{i\in\Vset} a_i=0.
%\end{align*}
%\begin{remark}
While the system \eqref{eq:GF-orig}, \eqref{eq:kirchhoff} is a constrained gradient flow
for the energy \eqref{eq:energydisc}, the system \eqref{eq:pineq}, \eqref{eq:auxineqproduction}
does not have a gradient flow structure in full generality.
%\end{remark}

%We usually denote the total amount of transporter in a cell by  $\mathcal{P}_{tot}$. Cells produce or degrade  transporter proteins and these proteins are then cycled back and forth to the membranes such that $\mathcal{P}_{tot} = \sum_{j \in \mathcal{N}(i)} \mathcal{P}_{ij}$. 

%%%%%%%%%%%%%%%%%%%%%%%%%%%%%%%%%%%%%%%%%%%%%%%%%%%%%%%%%%%%%
\section{Global existence and nonnegativity of solutions to the adapted Hu-Cai model}\label{sec:AnalysisD}

\begin{theorem}
	Let $0 < \kappa-\gamma\leq 1$ and fix $T>0$.
	The system \eqref{eq:pineq}, \eqref{eq:auxineqproduction} subject to the initial datum \eqref{pin-IC}--\eqref{aux-IC}
	has a solution $X_{ij}\in C^1(0,T)$, $a_i\in C^1(0,T)$, satisfying
	%\( \label{nonneg}
	$X_{ij}(t) \geq 0, a_i(t) > 0$  for all $t\in[0,T)$ and $i,j\in \Vset$.
	%\)
	Moreover, if $S_i=0$ for all $i\in\Vset$ in \eqref{eq:auxineqproduction}, then $a_i$ is uniformly globally bounded,
	i.e., there exists a constant $\alpha>0$ such that
	\begin{align} \label{global bound}
	a_i(t) \leq \alpha\qquad \text{for all~} t\in[0,\infty) \text{ and } i\in \Vset.
	\end{align}
\end{theorem}

\begin{proof}
	\textbf{Nonnegativity for $X_{ij}$.}
	%Since $\mathcal{P}_{ij}(0)=0$ implies $\mathcal{P}_{ij}(t)=0$ for all $t\geq 0$, we can restrict ourselves to $\mathcal{P}_{ij}(0)>0$.
	With \eqref{eq:pineq} we have
	$\tot{X_{ij}}{t} \geq  -\sigma  \tau X_{ij},$ 
	as long as the solution exists. Consequently, $X_{ij}(0) \geq 0$ implies
	$X_{ij}(t) \geq 0$ on the interval of existence.
	
	\textbf{Boundedness for $|a_{i}|$.}
	Let us denote the adjacency matrix of the graph $G=(\Vset,\Eset)$ by $\mathbb{A}\in \R^{n\times n}$, i.e.\ its entries are given by
	\begin{align}\label{adjM}
	\mathbb{A}_{ij}=\begin{cases}
	0 & \text{if }(i,j)\notin \Eset, \\
	1 & \text{if }(i,j)\in \Eset.
	\end{cases}
	\end{align} 
	For the solutions $a_i$ of the auxin equation \eqref{eq:auxineqproduction} on their joint interval of existence we have
	\begin{align*}
	\frac{1}{2}\frac{\di}{\di t}\sum_{i=1}^N a_i^2 &= \sum_{i=1}^N S_i a_i - \sum_{i=1}^N I_i a_i^2
	%  \sum_{i=1}^N a_i \frac{\di a_i}{\di t}\\&=
	+ \delta\sum_{i=1}^N  \sum_{j=1}^N \mathbb{A}_{ij} X_{ij}a_i\bl a_j-a_i\br\\
	& \leq \sum_{i=1}^N S_i a_i - \frac{\delta}{2} \sum_{i=1}^N  \sum_{j=1}^N \mathbb{A}_{ij} X_{ij}\bl a_i-a_j\br^2\,,
	%\frac{1}{2}\delta\sum_{i=1}^N  \sum_{j=1}^N \mathbb{A}_{ij} \mathcal{P}_{ij}\bl a_i\bl a_j-a_i\br+a_j\bl a_i-a_j\br\br\\&=
	\end{align*}
	where we used the nonnegativity of $I_i$ in the estimate and the usual symmetrization trick
	(recall that both $\mathbb{A}_{ij}$ and $X_{ij}$ are symmetric).
	Now, due to the nonnegativity of $X_{ij}$, we have
	\begin{align*}
	\frac{1}{2}\frac{\di}{\di t}\sum_{i=1}^N a_i^2 \leq \sum_{i=1}^N S_i a_i
	\leq \left( \sum_{i=1}^N S_i^2 \right)^{1/2} \left( \sum_{i=1}^N a_i^2 \right)^{1/2},
	\end{align*}
	implying at most quadratic growth of $a_i^2$ in time, i.e., at most linear growth of $|a_i|=|a_i|(t)$.
	Clearly, if $S_i=0$ for all $i\in\Vset$, then we have the uniform bound \eqref{global bound} with
	$$\alpha :=\sqrt{\sum_{i=1}^N a_i(0)^2}.$$
	
	\textbf{Boundedness for $X_{ij}$.}
	Due to the nonnegativity of $X_{ij}$ we have
	\[
	\tot{X_{ij}}{t} \leq \sigma \frac{|\mathcal{Q}_{ij}|^{\kappa}}{X_{ij}^{\gamma}} L_{ij},
	\]
	and the boundedness of $|a_i|$ on bounded time intervals implies
	\[
	|\mathcal{Q}_{ij}|^\kappa = \left| X_{ij}\frac{a_j-a_i}{L_{ij}} \right|^\kappa \leq C |X_{ij}|^\kappa
	\]
	for a suitable constant $C>0$. Hence,
	\begin{align*}
	\tot{X_{ij}}{t} \leq C X_{ij}^{\kappa-\gamma},
	\end{align*}
	and, therefore, for $0 < \kappa-\gamma < 1$, $X_{ij}=X_{ij}(t)$ grows at most algebraically in time,
	while for $\kappa-\gamma = 1$ the growth is at most exponential.
	
	\textbf{Positivity for $a_i$.}
	According to the assumption, there exists $\underline{a} >0$ such that $a_i(0) \geq \underline{a}$ for all $i\in\Vset$.
	Let us assume that $t_0<+\infty$ is the first instant when any of the curves $a_i=a_i(t)$ hits zero.
	Due to continuity, we have $t_0>0$, and, clearly, $a_i(t) > 0$ for $t\in [0,t_0)$ for all $i\in\Vset$.
	With the nonnegativity of the sources $S_i\geq 0$, \eqref{eq:auxineqproduction} implies
	\[
	\tot{a_i}{t} \geq - I_ia_i + \delta\sum_{j\in \Neigh(i)} X_{ij}\frac{a_j-a_i}{L_{ij}}\qquad \text{for~}i\in \Vset,\, t> 0,
	\]
	and with the nonnegativity of $X_{ij}$ we have
	\[
	\tot{a_i}{t} \geq - I_ia_i - \delta \left( \sum_{j\in \Neigh(i)} \frac{X_{ij}}{L_{ij}} \right) a_i \qquad \text{for~}i\in \Vset,\, t\in (0,t_0).
	\]
	Finally, since $X_{ij}=X_{ij}(t)$ grow at most exponentially in time, there exist constants $C$, $\lambda>0$
	independent of $t_0$ such that
	\[
	\tot{a_i}{t} \geq - C e^{\lambda t} a_i \qquad \text{for~}i\in \Vset,\, t\in (0,t_0),
	\]
	implying
	\[
	a_i(t) \geq a_i(0) \exp \left( \frac{\lambda}{C} (1-\exp(\lambda t)) \right).
	\]
	Therefore, $a_i(t_0) > 0$ for all $i\in\Vset$, a contradiction to the assumption $t_0 <+\infty$.
\end{proof}

Note that under the relaxed initial condition 
\begin{align}\label{eq:auxin:relaxedic}
a_{i}(0) = a_{i}^0 \geq 0\qquad \text{for all~}i \in \Vset
\end{align}
with an initial auxin concentration $\sum_{i\in\Vset} a_i(0)>0$ some cells may get no auxin over time. If $a_i(0)=0$ for some $i\in\Vset$, it follows from \eqref{eq:auxineqproduction} that cell $i$ gets no auxin as long as its neighboring cells have zero auxin. However, if $a_i(0)=0$ for some $i\in\Vset$ and $a_j(0)>0$ for some $j\in\mathcal{N}(i)$, then \eqref{eq:auxineqproduction} implies that $$\left.\tot{a_i(t)}{t}\right|_{t=0}\begin{cases} >0 & X_{ij}(0)>0, a_j(0)>0,\\ =0 & \text{otherwise.}\end{cases}$$  In particular, the relaxed initial condition \eqref{eq:auxin:relaxedic}  guarantees the nonnegativity for $a_i$.

\section{Murray's law}\label{sec:Murray}
In this section we demonstrate the validity of the Murray's law \cite{Murray1, Murray2}
for the steady states of %both the original model of Hu and Cai \eqref{eq:GF-orig}, \eqref{eq:kirchhoff}
the auxin transport activity model \eqref{eq:pineq}, \eqref{eq:auxineqproduction}.
Murray's law is a basic physical principle for transportation networks which
predicts the thickness or conductivity of branches,
such that the cost for transport and maintenance of the transport medium is minimized.
This law is observed in the vascular and respiratory systems of animals, xylem in plants,
and the respiratory system of insects \cite{Sherman}.

The stationary version of the auxin transport activity model \eqref{eq:pineq}, \eqref{eq:auxineqproduction} consists of the algebraic system
\begin{align}
\delta\sum_{j\in N(i)} \mathcal{Q}_{ji} &= S_i - I_i a_i  \qquad \text{for all~}i\in \Vset,   \label{AP-steady1} \\
\left( \frac{|\mathcal{Q}_{ij}|^{\kappa}}{X_{ij}^{\gamma+1}} - \tau \right) X_{ij} &= 0  \qquad \text{for all~}(i,j)\in \Eset.  \label{AP-steady2}
\end{align}
Noting that $\mathcal{Q}_{ij}=0$ if $X_{ij}=0$,
\eqref{AP-steady2} implies
\begin{align}  \label{AP-steady3}
|\mathcal{Q}_{ij}|^{\kappa} = \tau X_{ij}^{\gamma+1} \qquad \text{for all~}(i,j)\in \Eset.
\end{align}
Then, we rewrite \eqref{AP-steady1}  in the form
\[
\delta \sum_{j\in N^+(i)} |\mathcal{Q}_{ij}| + S_i  - I_i a_i = \delta \sum_{j\in N^-(i)} |\mathcal{Q}_{ij}| \qquad \text{for all~}i\in \Vset
\]
with
\[
N^+(i) := \{ j\in N(i); \; \mathcal{Q}_{ij} > 0 \},\qquad
N^-(i) := \{ j\in N(i); \; \mathcal{Q}_{ij} < 0 \}.
\]
Using \eqref{AP-steady3}, we have
\[
\delta \sum_{j\in N^+(i)} (\tau X_{ij}^{\gamma+1})^{1/\kappa} + S_i  - I_i a_i
= \delta \sum_{j\in N^-(i)}  (\tau X_{ij}^{\gamma+1})^{1/\kappa} \qquad \text{for all~}i\in \Vset.
\]
In particular, when all $I_i = 0$, we obtain the generalized Murray's law
\[
\delta \sum_{j\in N^+(i)} (\tau X_{ij}^{\gamma+1})^{1/\kappa} + S_i
= \delta \sum_{j\in N^-(i)}  (\tau X_{ij}^{\gamma+1})^{1/\kappa} \qquad \text{for all~}i\in \Vset.
\]

\section{Numerical simulation}\label{sec:num}

In this section, we provide numerical results for the discrete model \eqref{eq:auxineqproduction}--\eqref{eq:pineq}. Since the  problem is stiff, implicit formulas are necessary and we consider a multi-step solver based on the numerical differentiation formulas of orders 1 to 5 \cite{matlab}.  

We consider a planar graph $G=(\Vset,\Eset)$, whose vertices and edges define a diamond shaped geometry embedded in the two-dimensional domain $\Omega=(-0.5,2)\times (-1.5,0.5)$ with $|V|=81$ vertices and $|E|=208$ edges. Let $(x^i,y^i)$ denote the position of vertex $i\in V$. We assume that the source terms $S_i\geq 0$ are positive on the subset of vertices
\begin{align*}
V^+:=\{i\in V;~ x^i\leq -0.4\},
\end{align*}
and vanish on its complement $V\backslash \Vset^+$, 
\begin{align*}
S_i:=\begin{cases}\xi_S, &i\in \Vset^+,\\0,& i\in V\backslash V^+,\end{cases}
\end{align*}
where $\xi_S:=100$, implying that we have a single source in the top corner of the diamond. The decay terms $I_i,~ i\in V,$ are assumed to positive on the complement $V\backslash V^+$,
\begin{align*}
I_i:=\begin{cases}0, &i\in \Vset^+,\\\xi_I,& i\in V\backslash V^+,\end{cases}
\end{align*}
where $\xi_I:=1$.
Note that in terms of the distribution of source and sink terms, we consider the same situation as in \cite{HaKrMa}.
We prescribe the initial condition $\bar{X}_{ij}:=1$ for every $(i,j)\in\Eset$ and $a_i:=1$ for all $i\in \Vset$, unless stated otherwise. Besides, we consider $\delta:=1$, $\sigma:=1$, $\kappa:=2$, $\gamma:=0.5$ and $\tau:=1$ in the numerical simulations, if not stated otherwise.

In the sequel, we present the stationary solutions obtained by solving the system \eqref{eq:auxineqproduction}--\eqref{eq:pineq}. We plot the value of the  transport activity $X_{ij}$ for every edge $(i,j)\in\Eset$ in terms of its width and color. The auxin concentration in each cell $i\in \Vset$ is indicated by the color of that cell.

\begin{figure}[htbp]
	\centering
	\subfloat[$\epsilon=0.5$] {\includegraphics[width=0.24\textwidth]{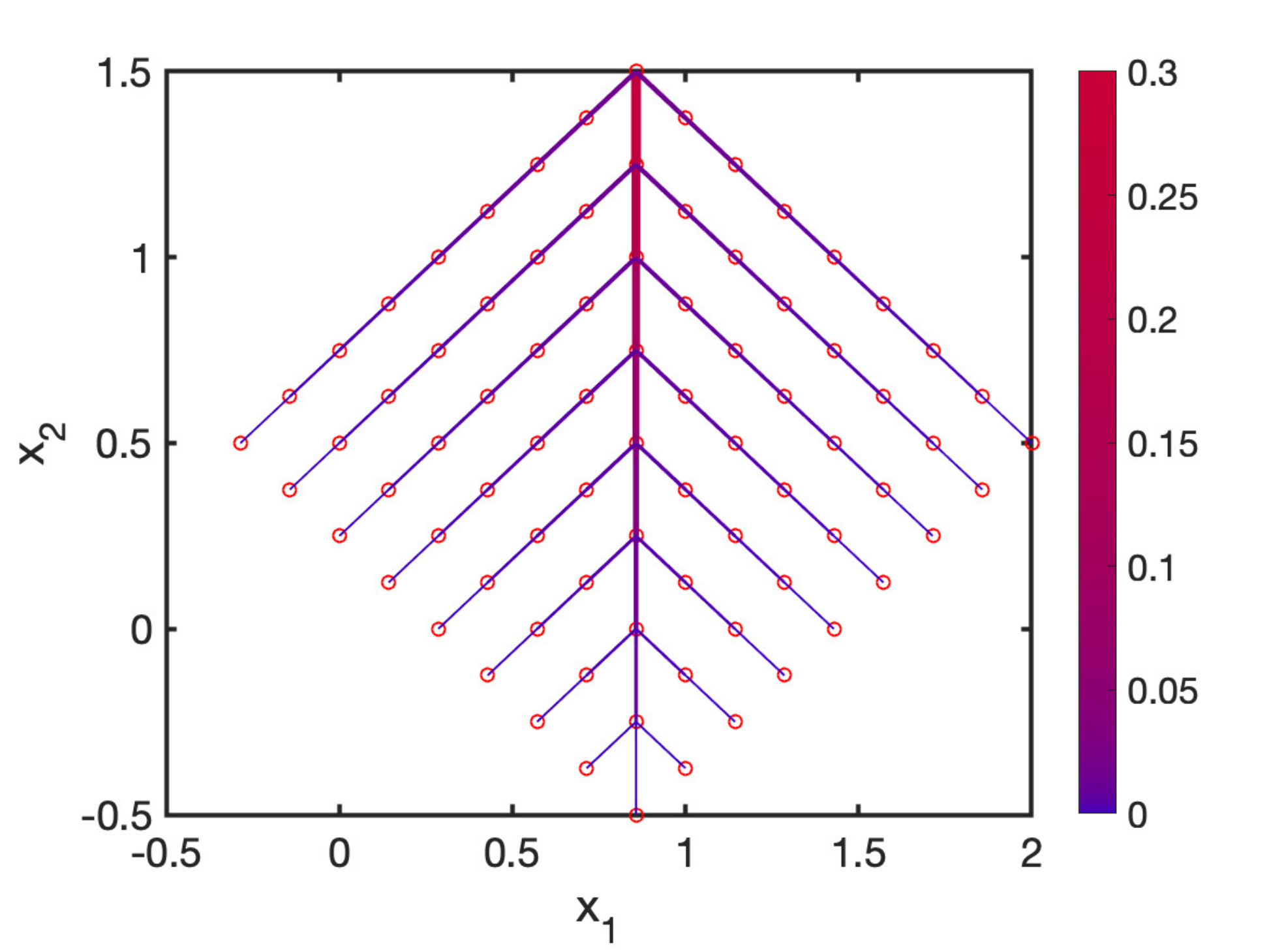}}
	\subfloat[$\epsilon=1$] {\includegraphics[width=0.24\textwidth]{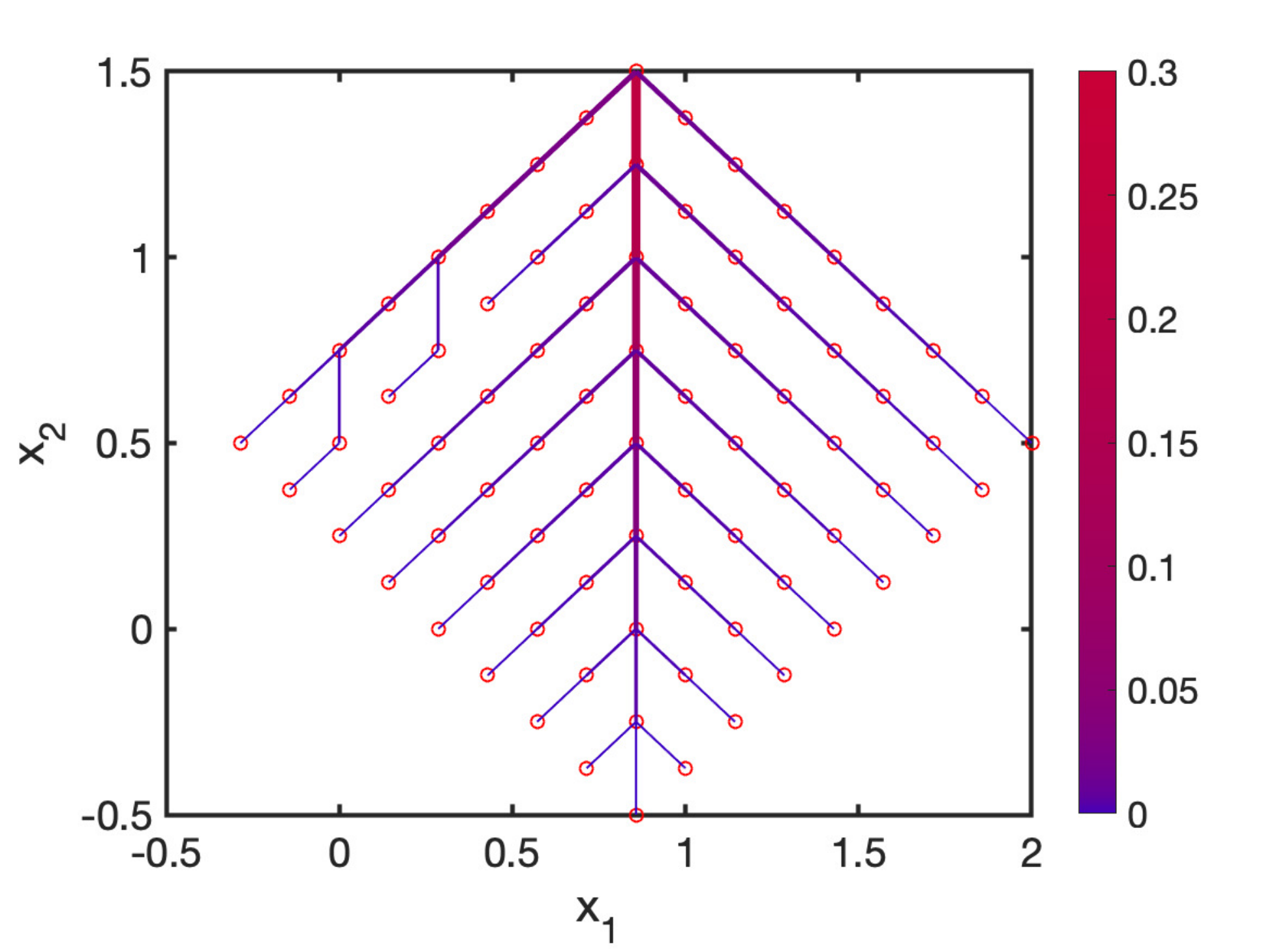}}
	\subfloat[$\epsilon=5$] {\includegraphics[width=0.24\textwidth]{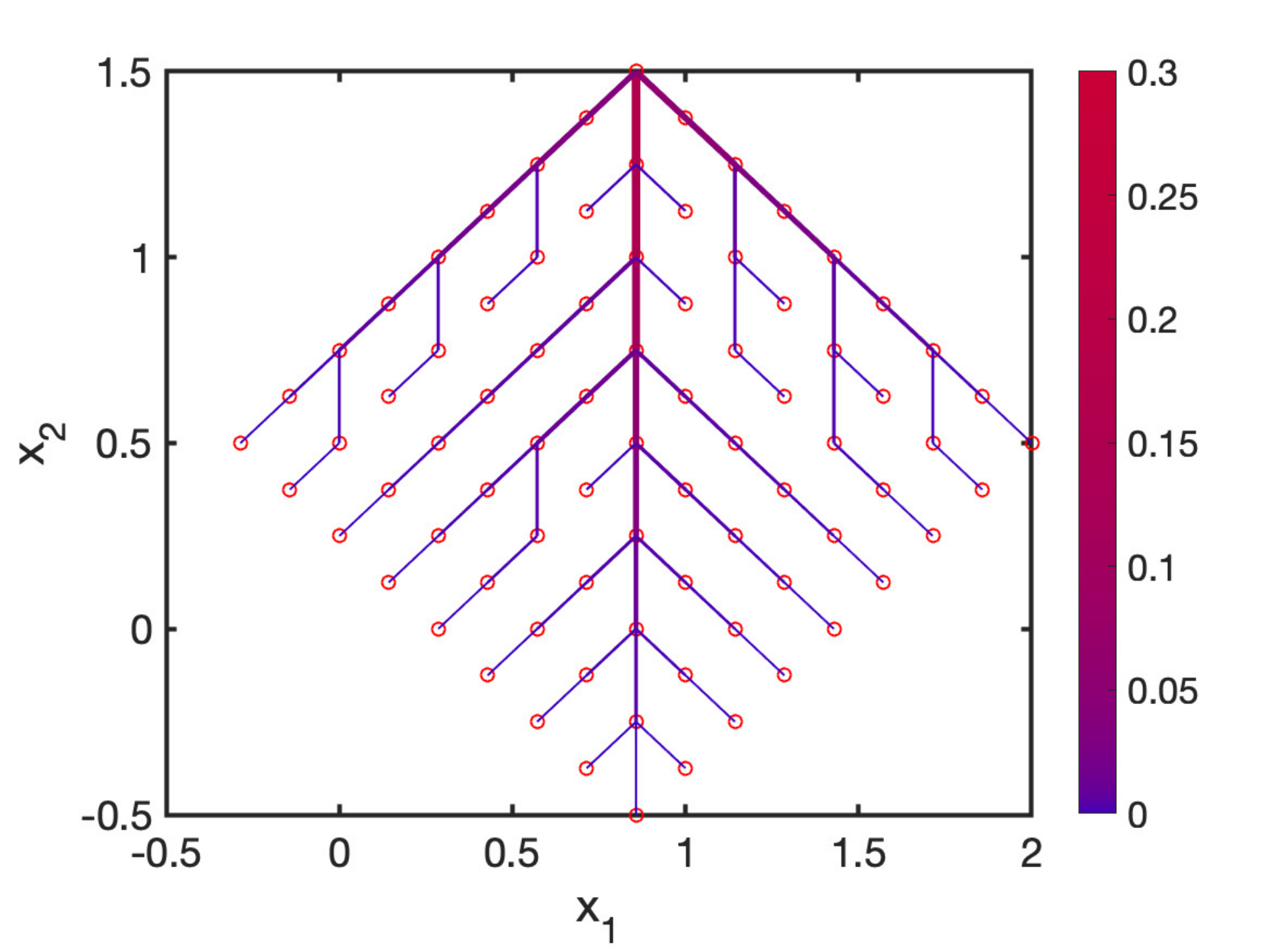}}
	\subfloat[$\epsilon=10$] {\includegraphics[width=0.24\textwidth]{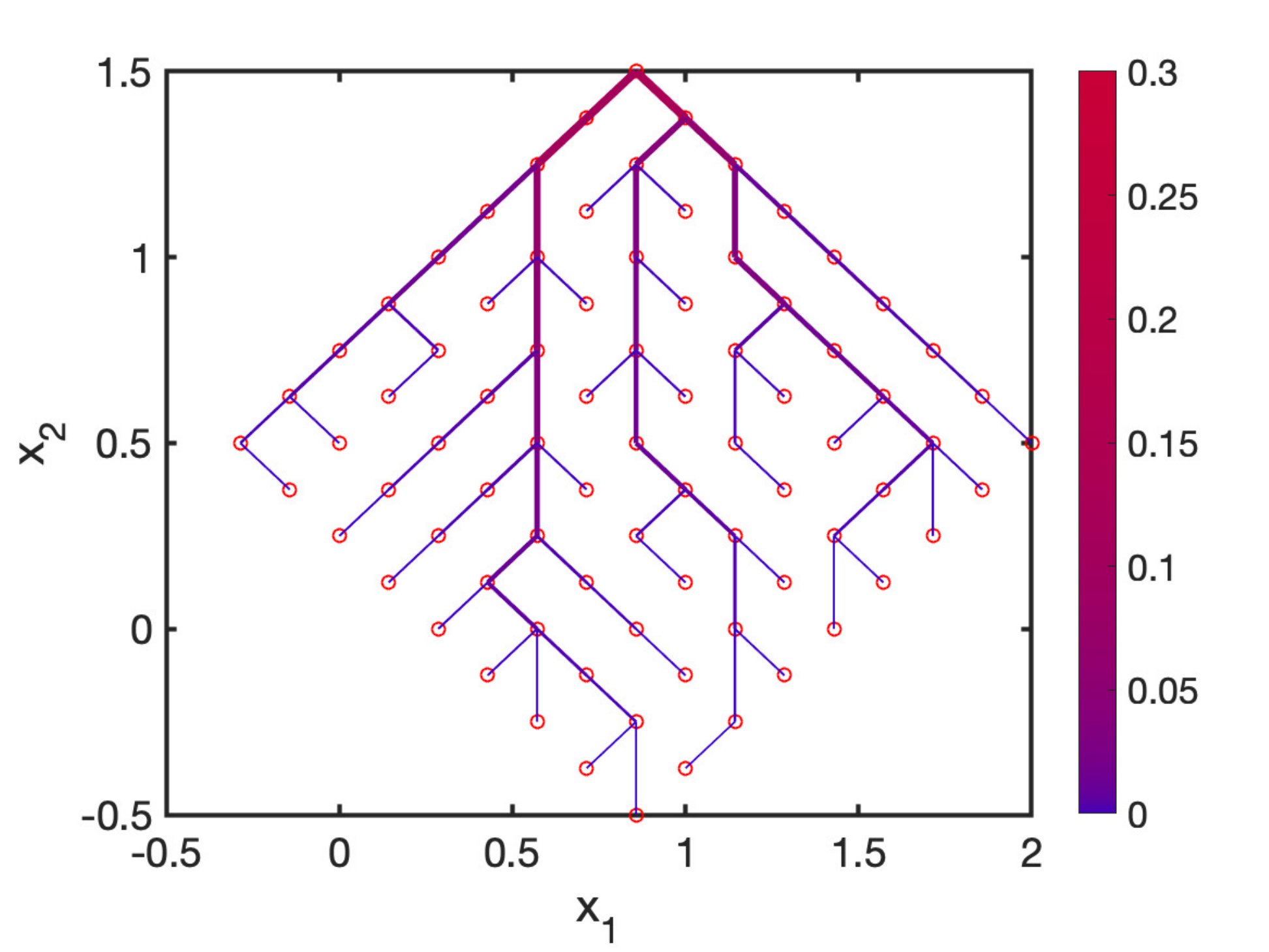}}
	\caption{Steady states for transport activity for perturbation $\epsilon\mathcal{U}(0,1)$ of the initial transport activity   $\bar{X}_{ij}$ with initial data $\bar{X}_{ij}, \bar{a}_{i}$.}\label{fig:pinperturbation}
\end{figure}

\begin{figure}[htbp]
	\centering
	%	\subfloat[$\xi_S=1$] {\includegraphics[width=0.24\textwidth]{Gamma50e-2_InitSingleSourceStrength_1Aux.eps}}
	\subfloat[$\xi_S=10$] {\includegraphics[width=0.24\textwidth]{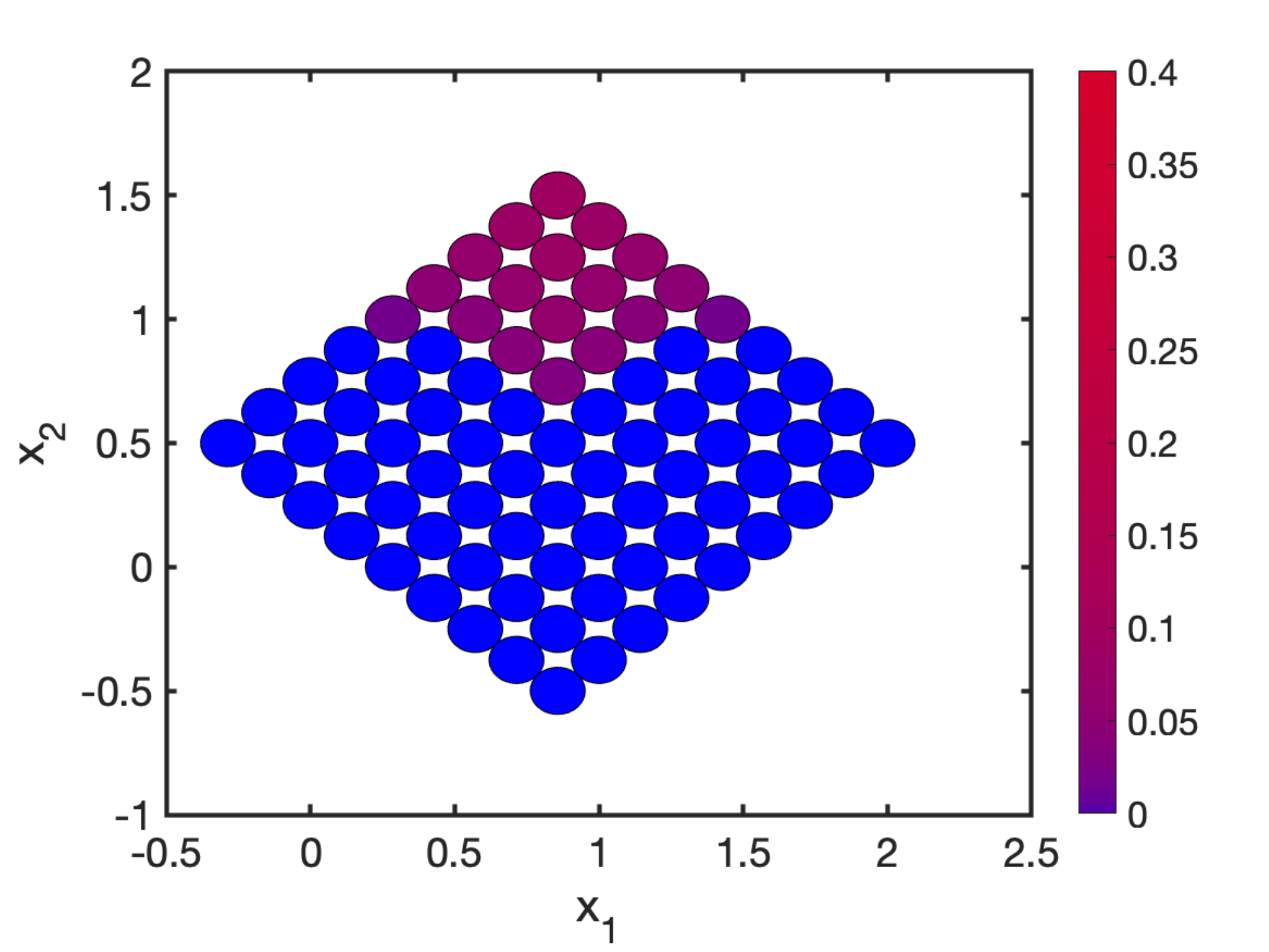}}
	\subfloat[$\xi_S=50$] {\includegraphics[width=0.24\textwidth]{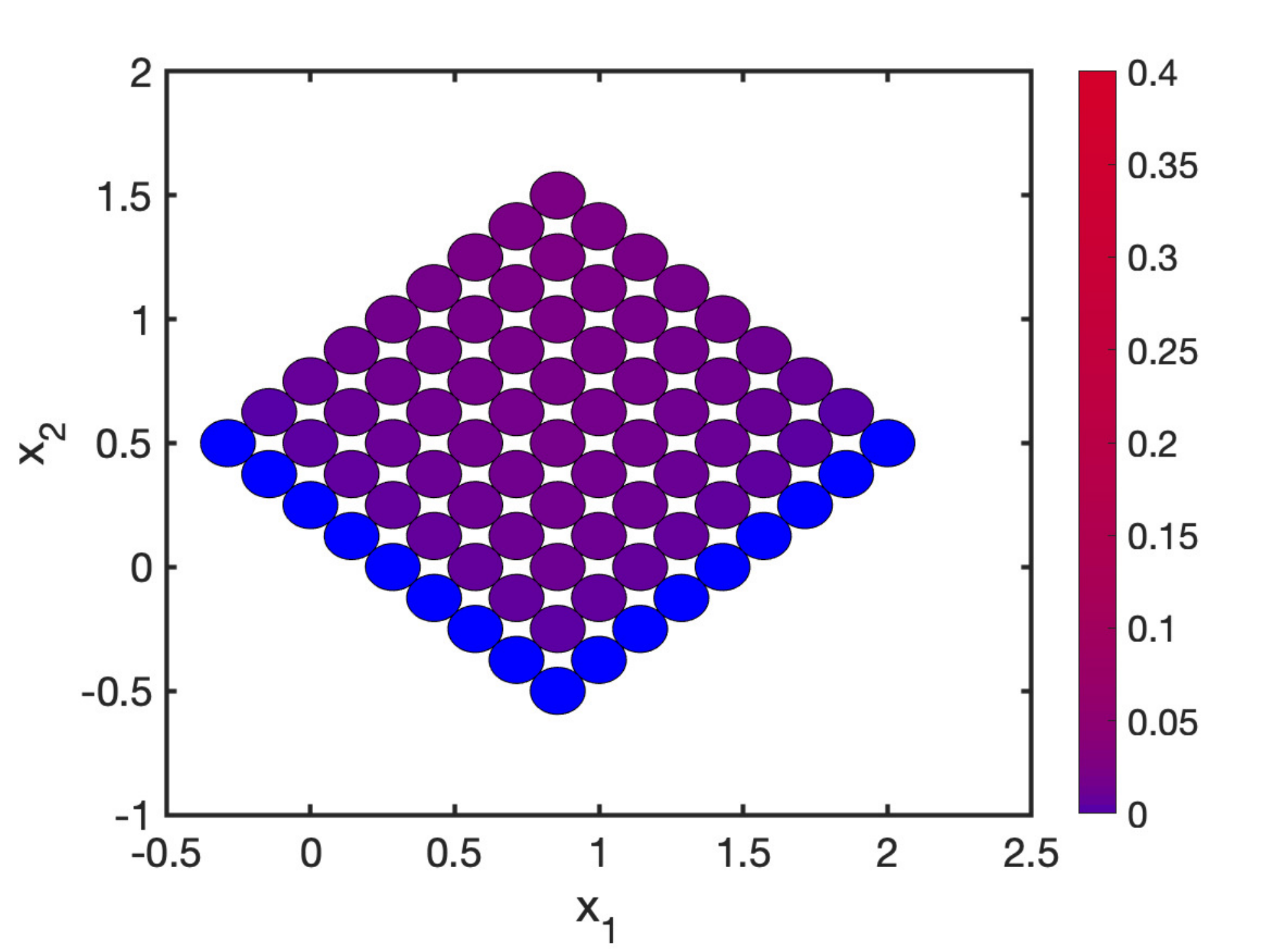}}
	\subfloat[$\xi_S=100$] {\includegraphics[width=0.24\textwidth]{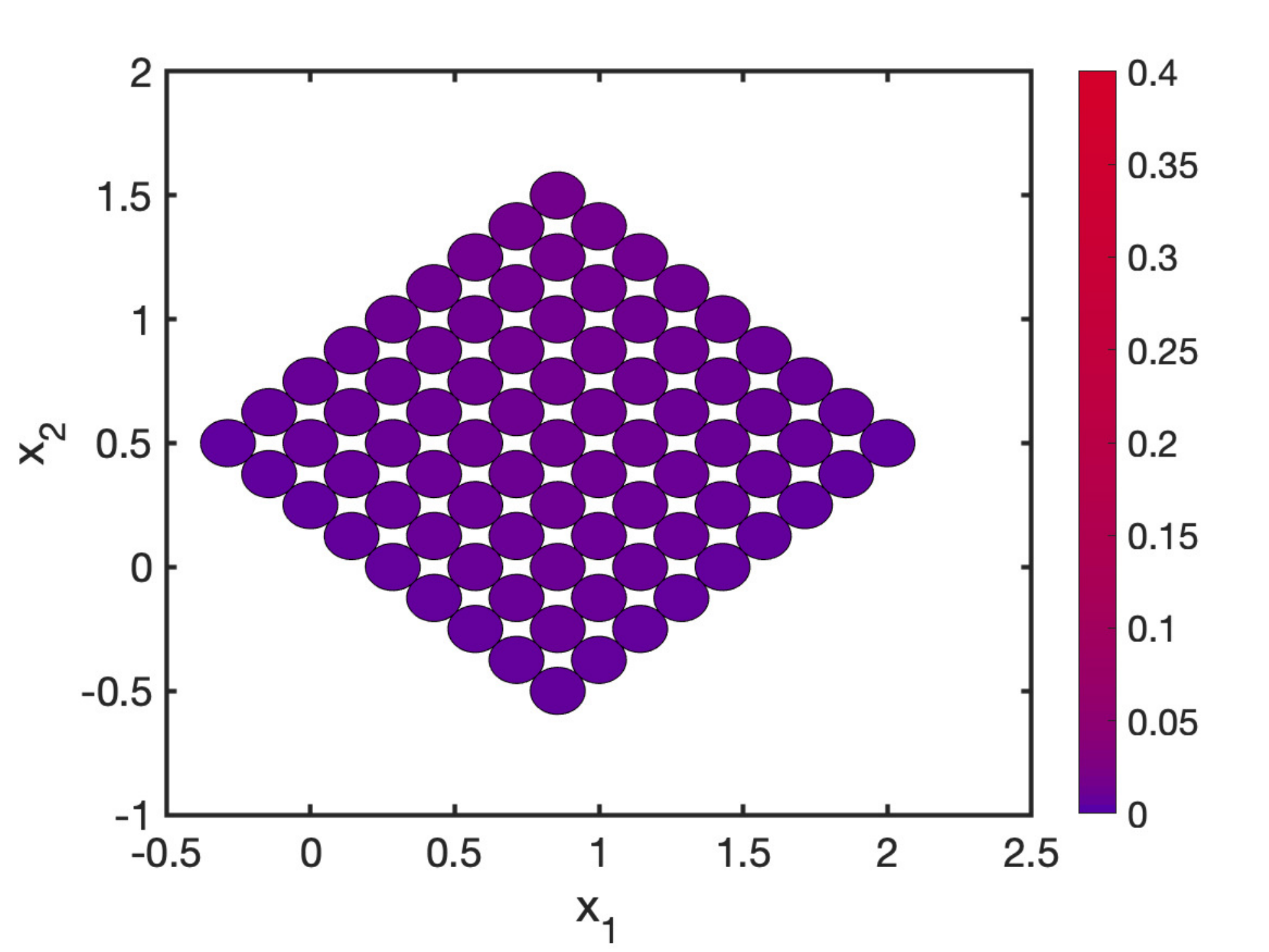}}
	\subfloat[$\xi_S=200$] {\includegraphics[width=0.24\textwidth]{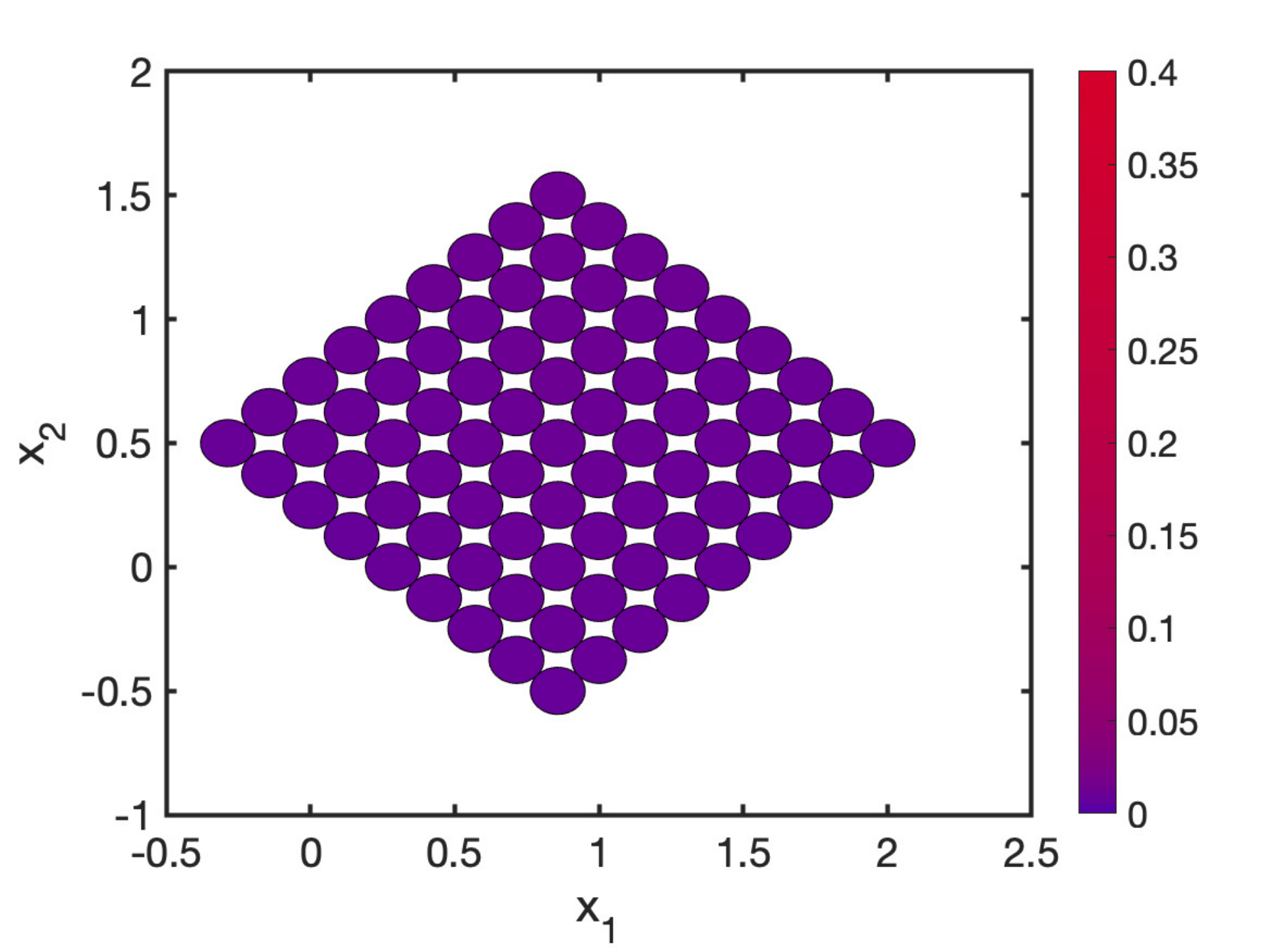}}

	%	\subfloat[$\xi_S=1$] {\includegraphics[width=0.24\textwidth]{Gamma50e-2_InitSingleSourceStrength_1PIN.eps}}
	\subfloat[$\xi_S=10$] {\includegraphics[width=0.24\textwidth]{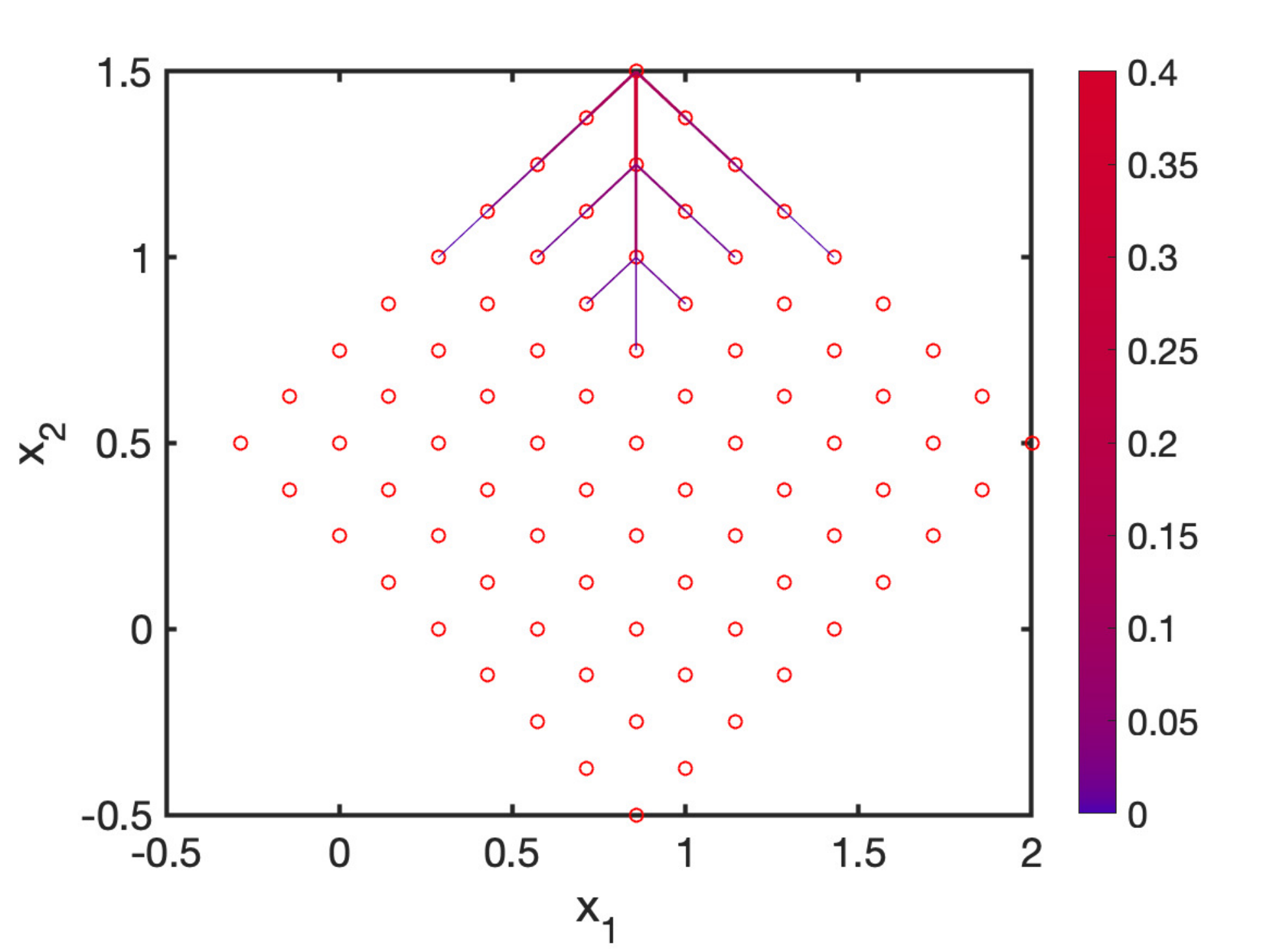}}
	\subfloat[$\xi_S=50$] {\includegraphics[width=0.24\textwidth]{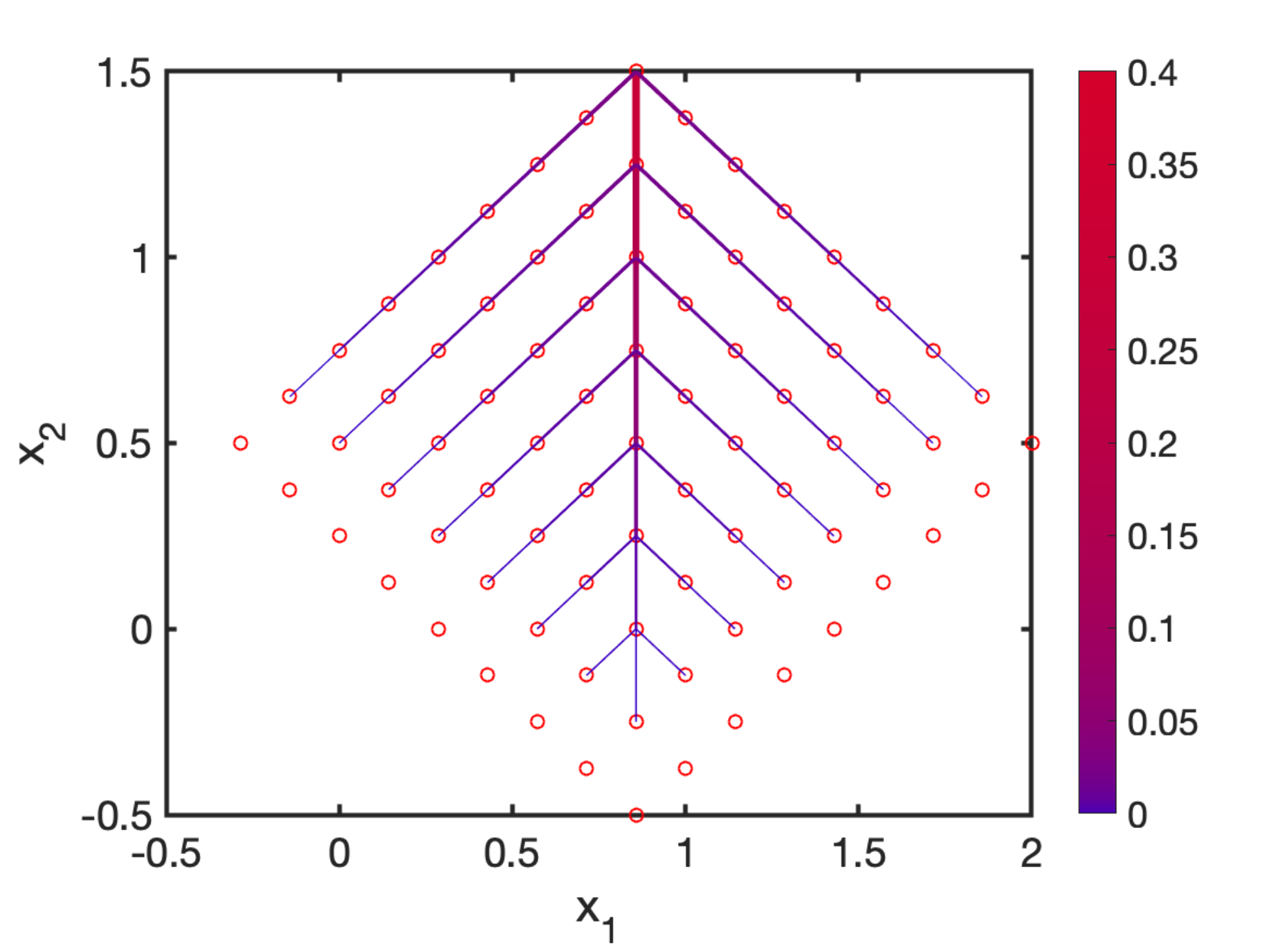}}
	\subfloat[$\xi_S=100$] {\includegraphics[width=0.24\textwidth]{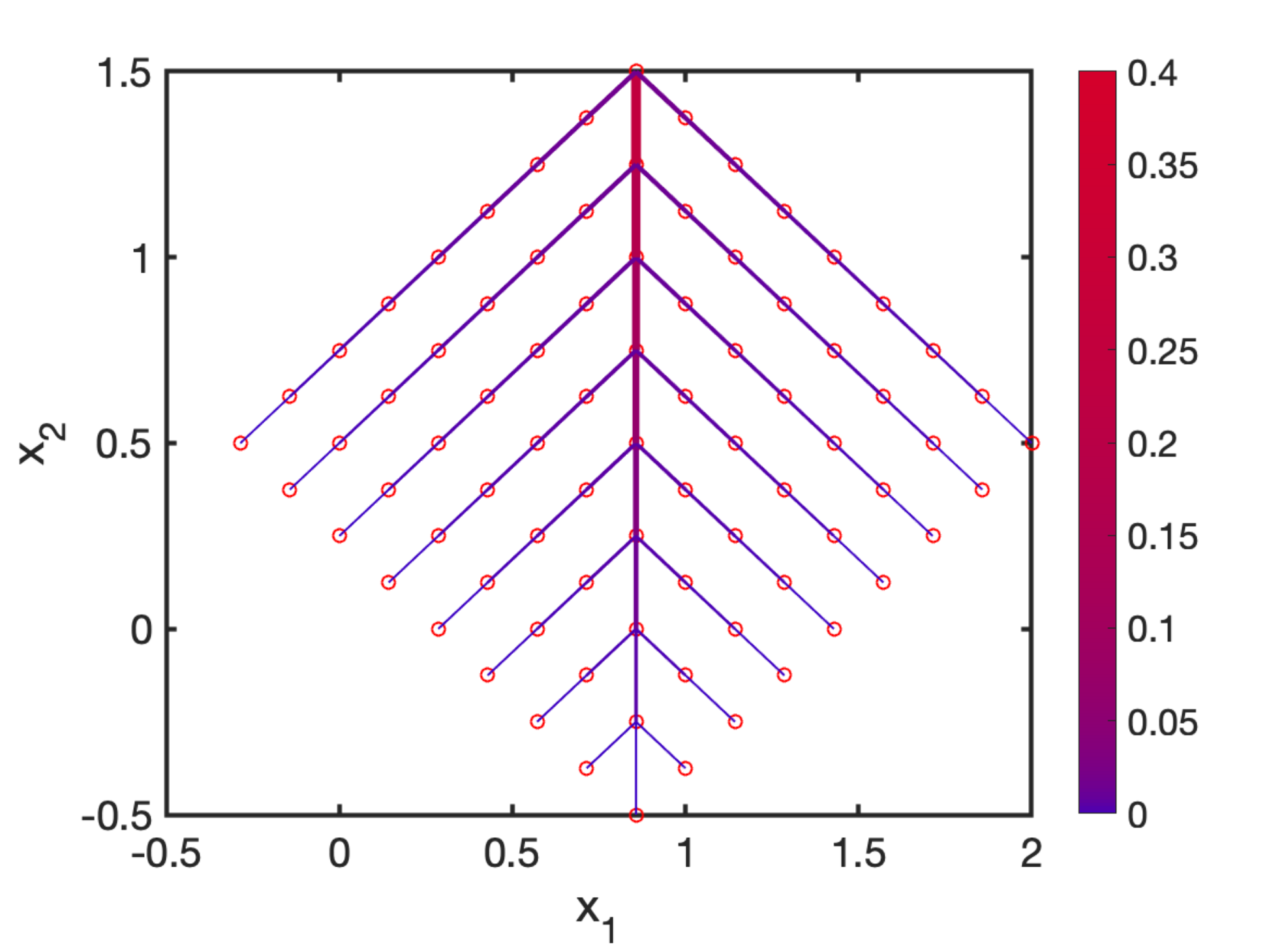}\label{fig:sourcestrengthorig}}
	\subfloat[$\xi_S=200$] {\includegraphics[width=0.24\textwidth]{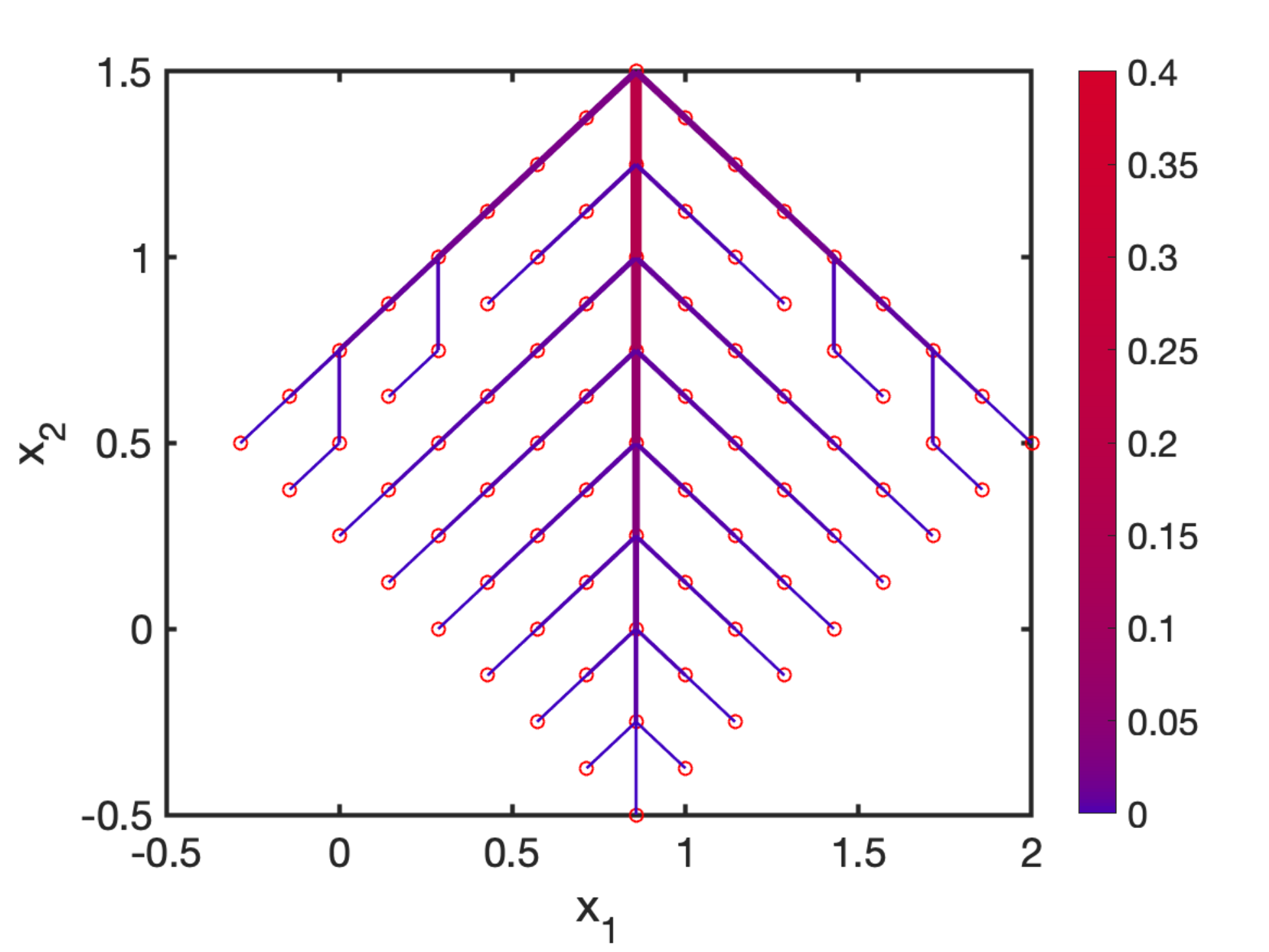}}
	\caption{Steady states  for auxin concentration and transport activity for different background source strengths $\xi_S$ with initial data $\bar{X}_{ij}, \bar{a}_{i}$.}\label{fig:sourcestrength}
\end{figure}

In Figure \ref{fig:pinperturbation}, we show the stationary transport activity for perturbed initial data $\bar{X}_{ij}$, i.e., we consider $\bar{X}_{ij}+\epsilon\mathcal{U}(0,1)$ instead of $\bar{X}_{ij}$ as initial data, where $\mathcal{U}(0,1)$ denotes  a uniformly distributed random variable on $[0,1]$. In particular, the resulting network is stable under small perturbation. This can be seen by comparing the results with Figure \subref*{fig:sourcestrengthorig} where the same parameters without perturbation are considered. 
%Similarly, we perturb the initial data for the auxin concentration in Figure \ref{fig:auxinperturbation}, i.e., we consider $\bar{a}_{i}+\epsilon\mathcal{U}(0,1)$ for all $i\in\Vset$ as initial data and show the resulting stationary PIN concentration. 
The perturbations of the initial data result in more complex steady states compared to the steady states obtained from unperturbed initial data.

%\begin{figure}[htbp]
%	\centering
%	\subfloat[$\epsilon=0.5$] {\includegraphics[width=0.24\textwidth]{Gamma50e-2_InitPIN_05PIN.eps}}
%	\subfloat[$\epsilon=1$] {\includegraphics[width=0.24\textwidth]{Gamma50e-2_InitPIN_1PIN.eps}}
%	\subfloat[$\epsilon=5$] {\includegraphics[width=0.24\textwidth]{Gamma50e-2_InitPIN_5PIN.eps}}
%	\subfloat[$\epsilon=10$] {\includegraphics[width=0.24\textwidth]{Gamma50e-2_InitPIN_10PIN.eps}}
%	\caption{Steady states for PIN concentration for perturbation $\epsilon\mathcal{U}(0,1)$ of the initial auxin concentration  $\bar{\mathcal{P}}_{ij}$ with initial data $\bar{\mathcal{P}}_{ij}, \bar{a}_{i}$.}\label{fig:auxinperturbation}
%\end{figure}

In Figure \ref{fig:sourcestrength} we vary the strength $\xi_S$ of the source in the top corner of the diamond.  As $\xi_S$ increases,   auxin is transported over a larger area, resulting in lower auxin levels and transport activity close to the source in the top corner of the diamond. Note that the area of large auxin levels and transport activities coincide in the steady states. Further note that not the entire graph is covered with auxin  for $\xi_S\in\{10,50\}$ and the resulting pattern is symmetric due to symmetric initial data for the auxin levels and the transport activity.

\begin{figure}[htbp]
	\centering
	\subfloat[$\xi_S=100$] {\includegraphics[width=0.24\textwidth]{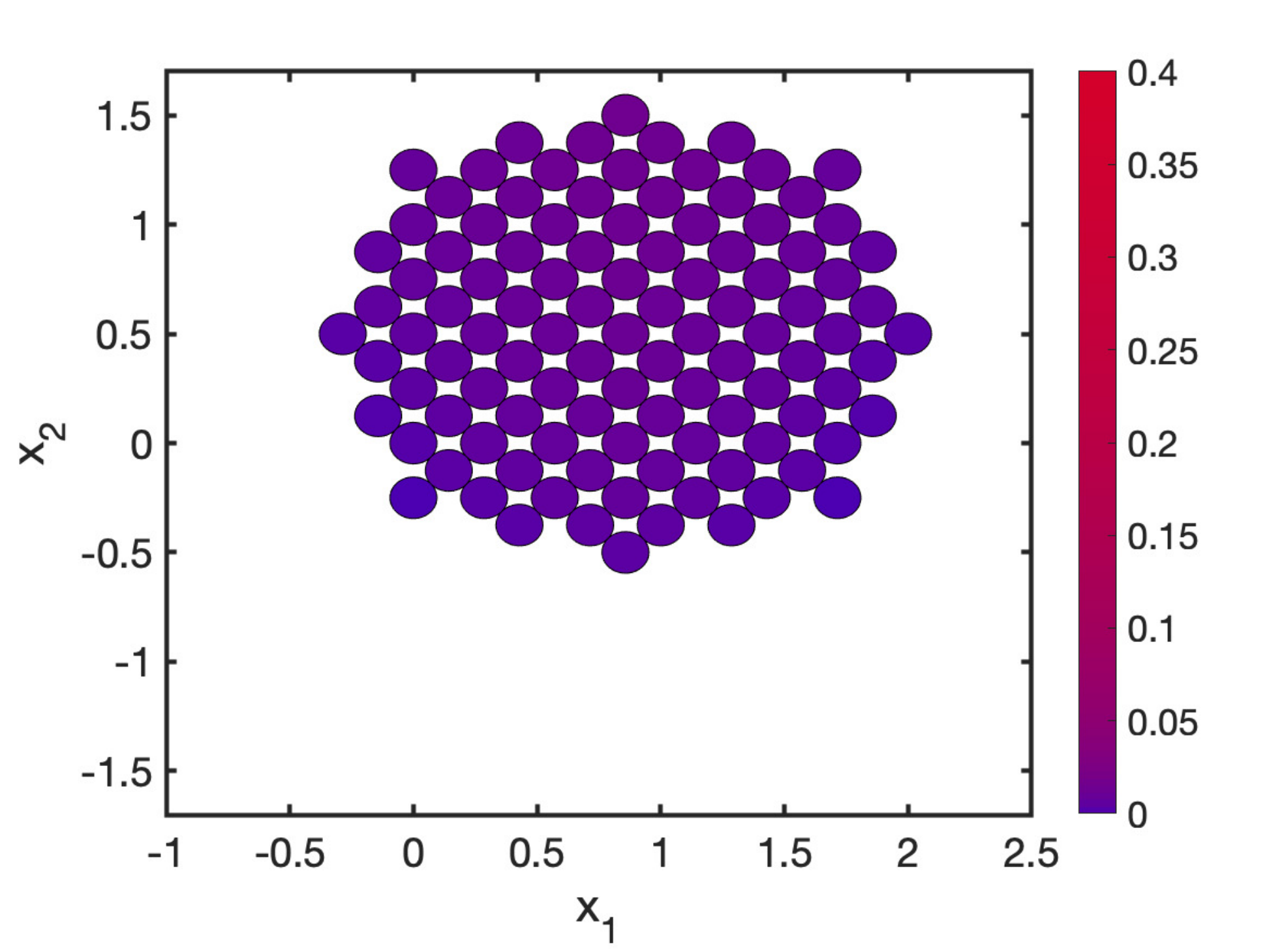}}
	\subfloat[$\xi_S=200$] {\includegraphics[width=0.24\textwidth]{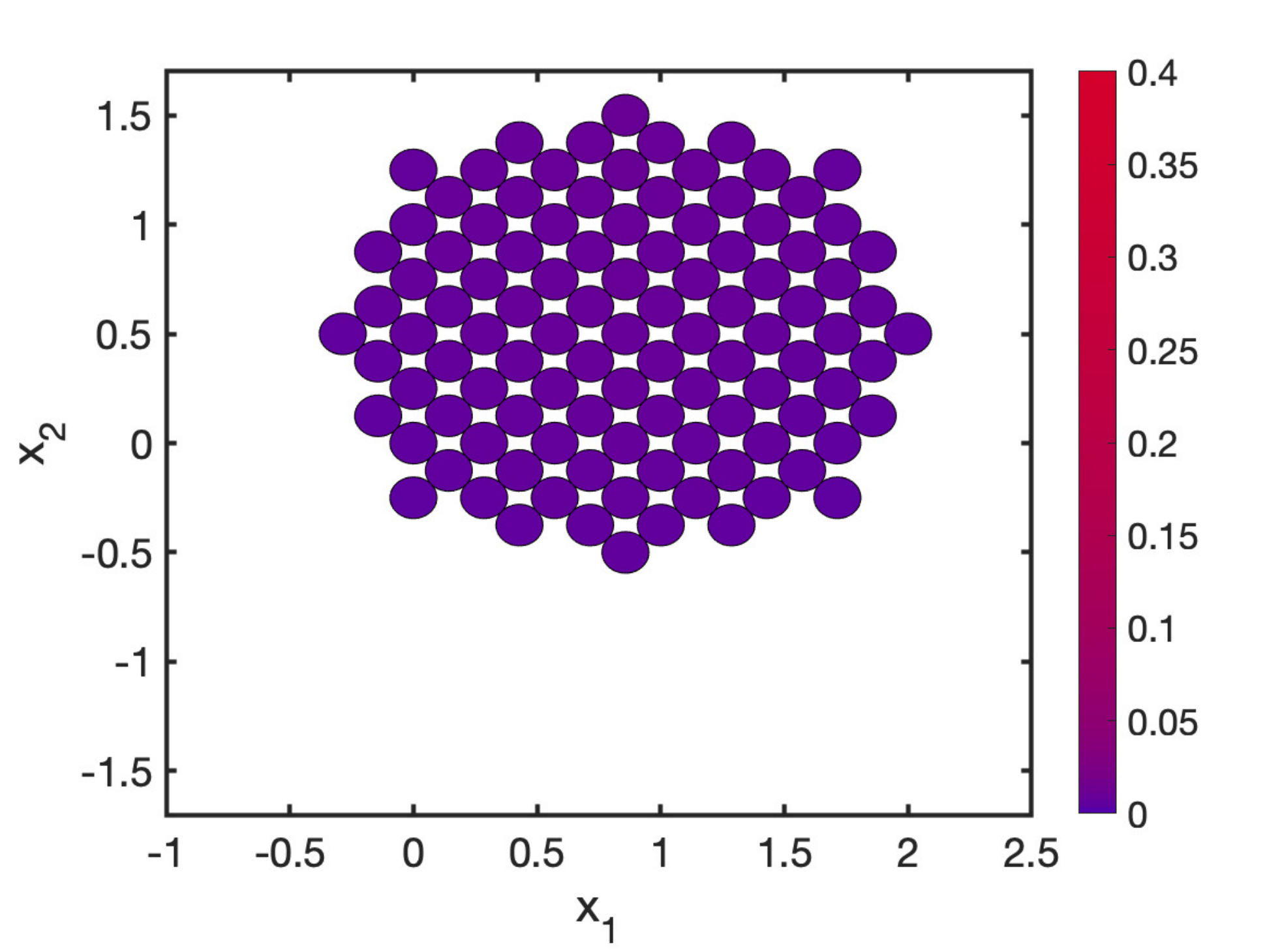}}
	\subfloat[$\xi_S=100$] {\includegraphics[width=0.24\textwidth]{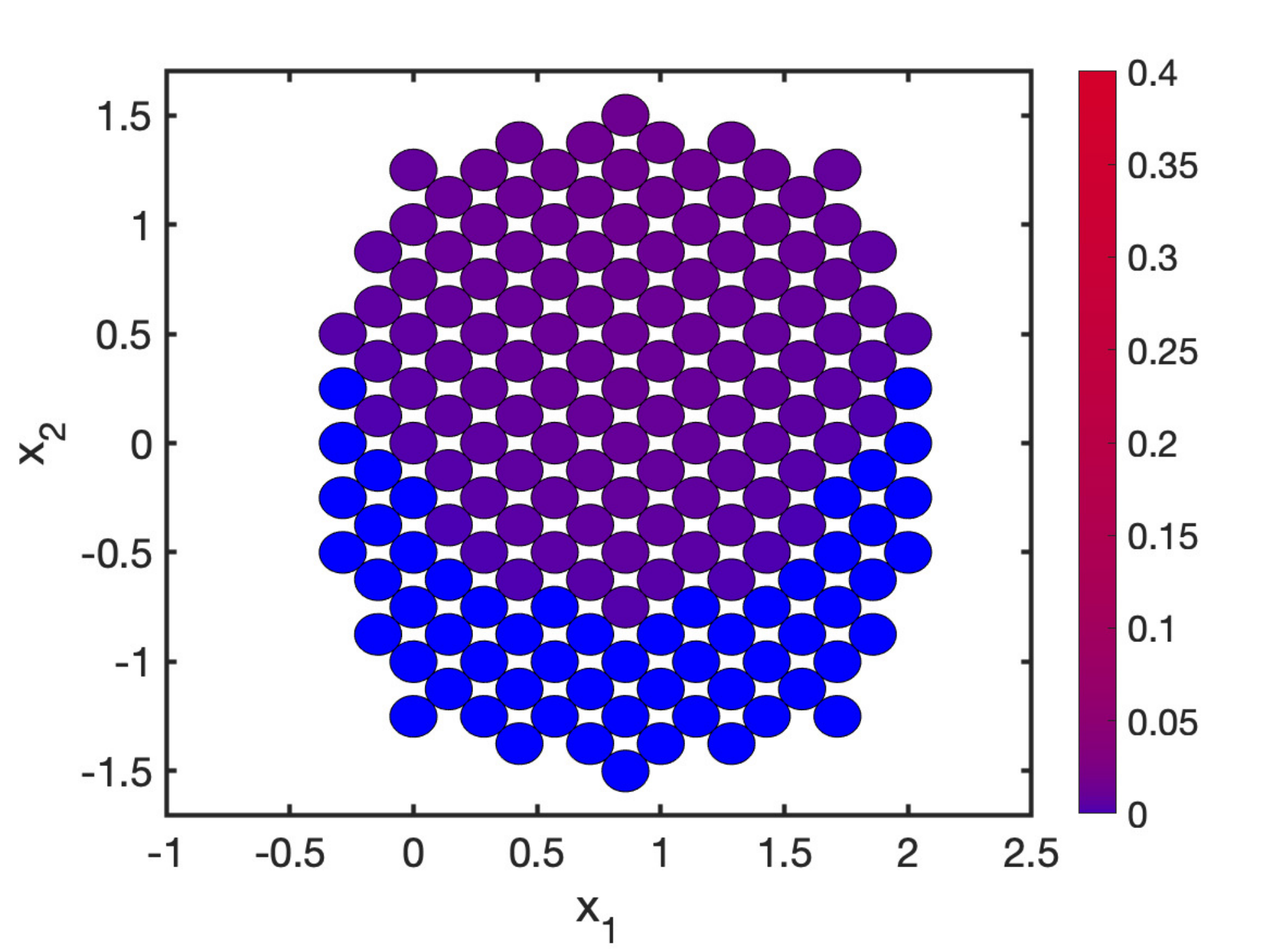}}
	\subfloat[$\xi_S=200$] {\includegraphics[width=0.24\textwidth]{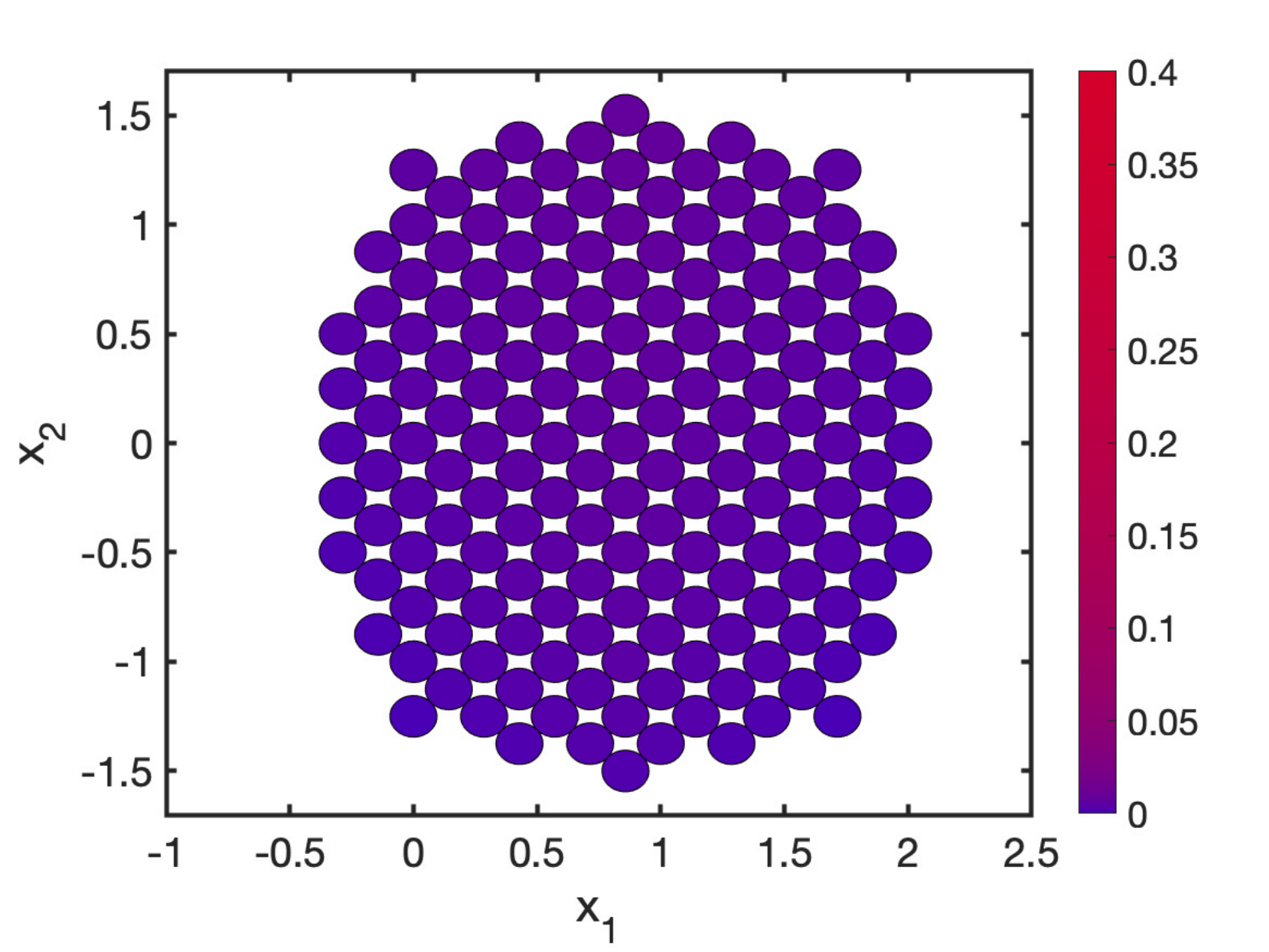}}\\
	\subfloat[$\xi_S=100$] {\includegraphics[width=0.24\textwidth]{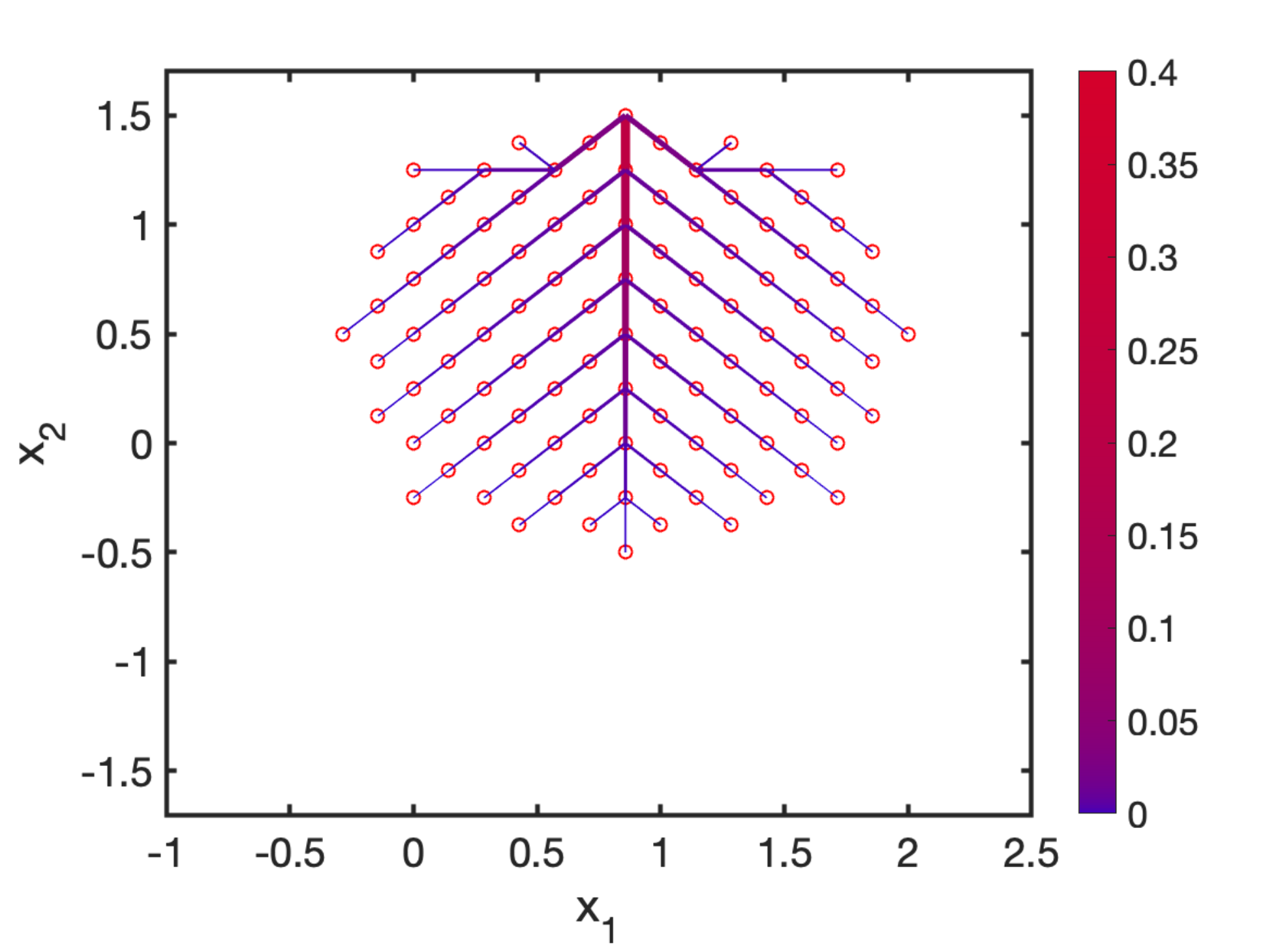}}
	\subfloat[$\xi_S=200$] {\includegraphics[width=0.24\textwidth]{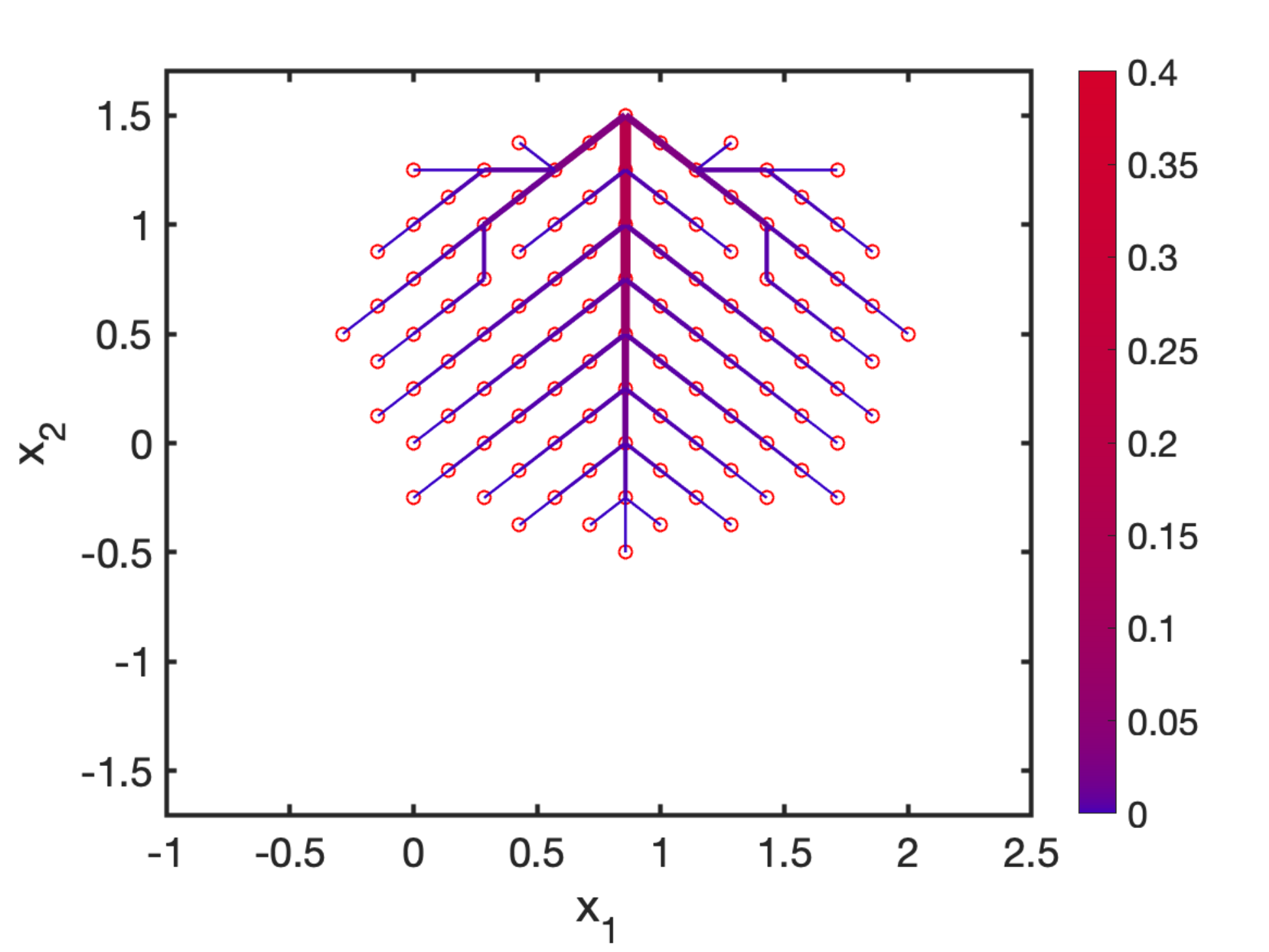}}
	\subfloat[$\xi_S=100$] {\includegraphics[width=0.24\textwidth]{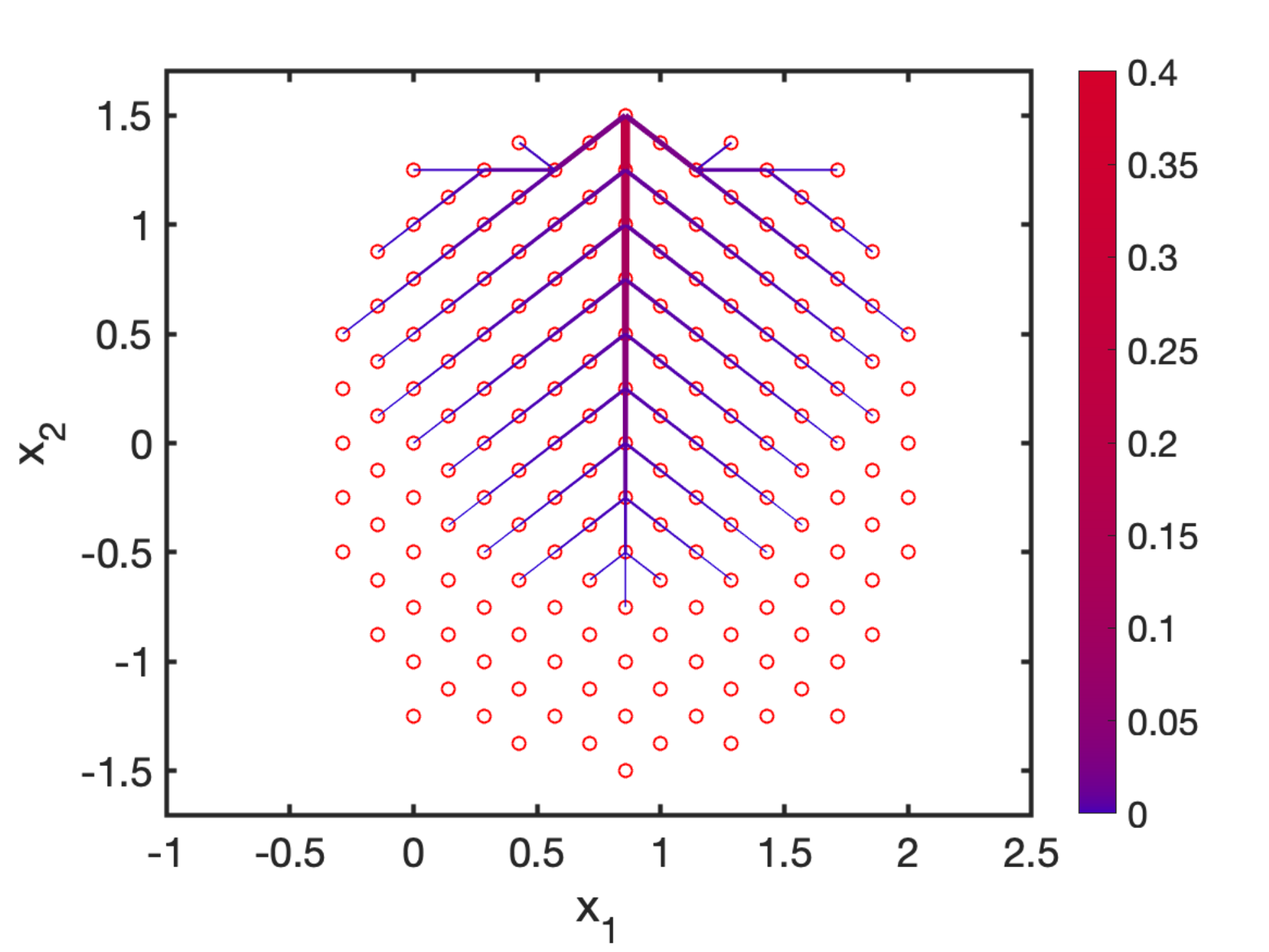}}
	\subfloat[$\xi_S=200$] {\includegraphics[width=0.24\textwidth]{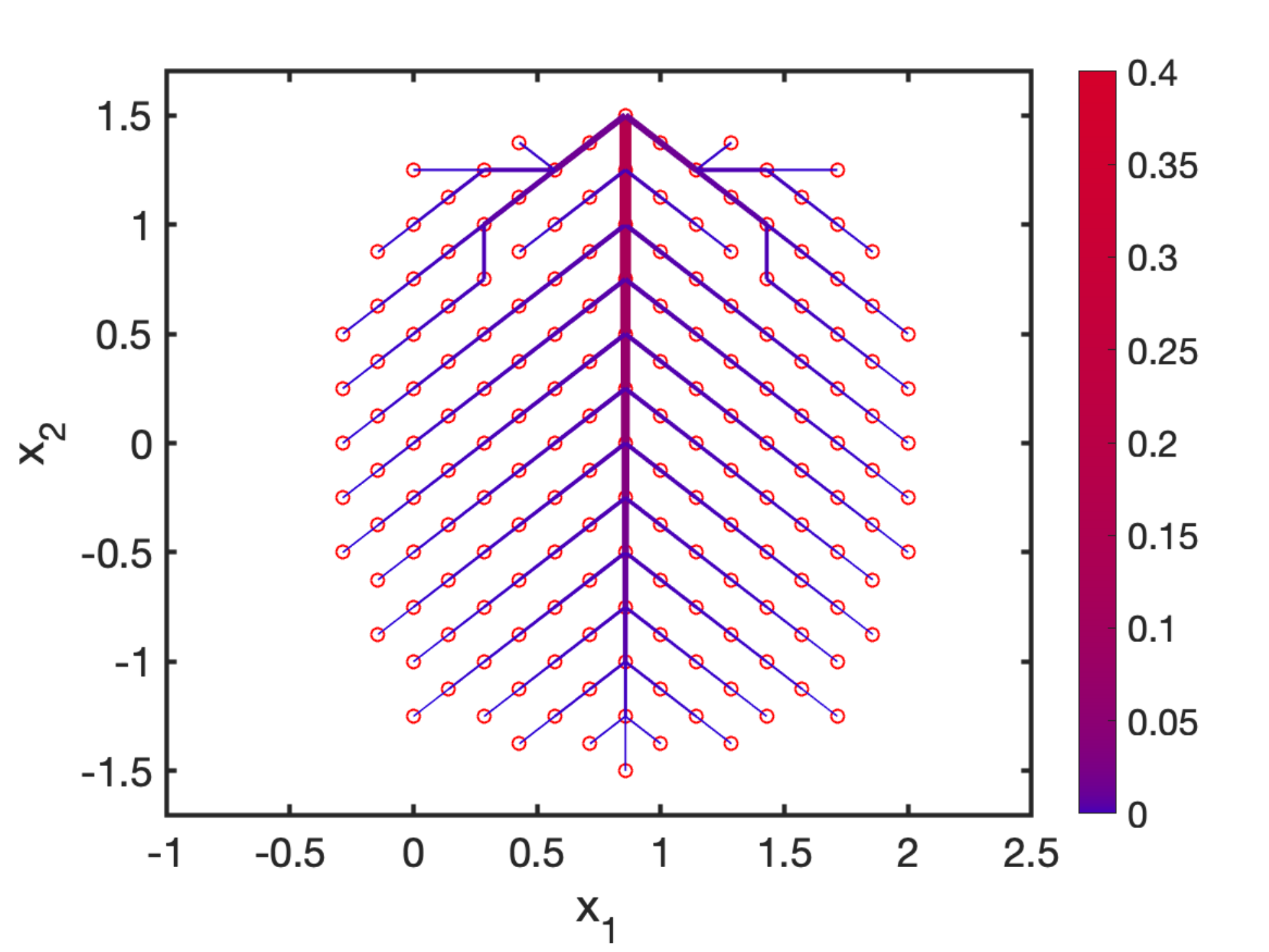}}
	\caption{Steady states  for auxin concentration and transport activity for different background source strengths $\xi_S$ and different grid shapes (round, oval) with initial data $\bar{X}_{ij}, \bar{a}_{i}$.}\label{fig:sourcestrengthgrids}
\end{figure}

In Figure \ref{fig:sourcestrengthgrids} we consider different grids (round, oval). As in Figure \ref{fig:sourcestrength} we vary the  strength $\xi_S$ of the source in the top middle corner of these grids. The resulting pattern formation  for round and oval grids is very similar to the patterns obtained with the same source strengths in Figure~\ref{fig:sourcestrength} for the diamond grid. In particular, this demonstrates the robustness of the model to variations of the underlying grid. Note that due to the larger size of the oval grid compared to the other considered grids, a stronger source is required for obtaining stationary patterns covering the entire simulation domain. 

%Similarly, we vary the background sink strength $\xi_I$ in Figure \ref{fig:sinkstrength}. While for  sink strengths sufficiently small, i.e., $\xi_I\in\{0.5,1\}$, the auxin and PINs are concentrated on the entire graph, the area of auxin and PIN concentration shrinks as the parameter $\xi_I$ increases. As before, the area of auxin and PIN concentration is of the same size and their absolute values increase as $\xi_I$ increases.
%
%\begin{figure}[htbp]
%	\centering
%	\subfloat[$\xi_I=0.5$] {\includegraphics[width=0.24\textwidth]{Gamma50e-2_InitBackgroundSinkStrength_05Aux.eps}}
%	\subfloat[$\xi_I=1$] {\includegraphics[width=0.24\textwidth]{Gamma50e-2_InitBackgroundSinkStrength_1Aux.eps}}
%	\subfloat[$\xi_I=5$] {\includegraphics[width=0.24\textwidth]{Gamma50e-2_InitBackgroundSinkStrength_5Aux.eps}}
%	\subfloat[$\xi_I=10$] {\includegraphics[width=0.24\textwidth]{Gamma50e-2_InitBackgroundSinkStrength_10Aux.eps}}
%	
%	\subfloat[$\xi_I=0.5$] {\includegraphics[width=0.24\textwidth]{Gamma50e-2_InitBackgroundSinkStrength_05PIN.eps}}
%	\subfloat[$\xi_I=1$] {\includegraphics[width=0.24\textwidth]{Gamma50e-2_InitBackgroundSinkStrength_1PIN.eps}}
%	\subfloat[$\xi_I=5$] {\includegraphics[width=0.24\textwidth]{Gamma50e-2_InitBackgroundSinkStrength_5PIN.eps}}
%	\subfloat[$\xi_I=10$] {\includegraphics[width=0.24\textwidth]{Gamma50e-2_InitBackgroundSinkStrength_10PIN.eps}}
%	\caption{Steady states  for auxin concentration and PIN concentration for different background sink strengths $\xi_I$ with initial data $\bar{\mathcal{P}}_{ij}, \bar{a}_{i}$.}\label{fig:sinkstrength}
%\end{figure}

\begin{figure}[htbp]
	\centering
	\subfloat[$\xi_I^C=10$] {\includegraphics[width=0.24\textwidth]{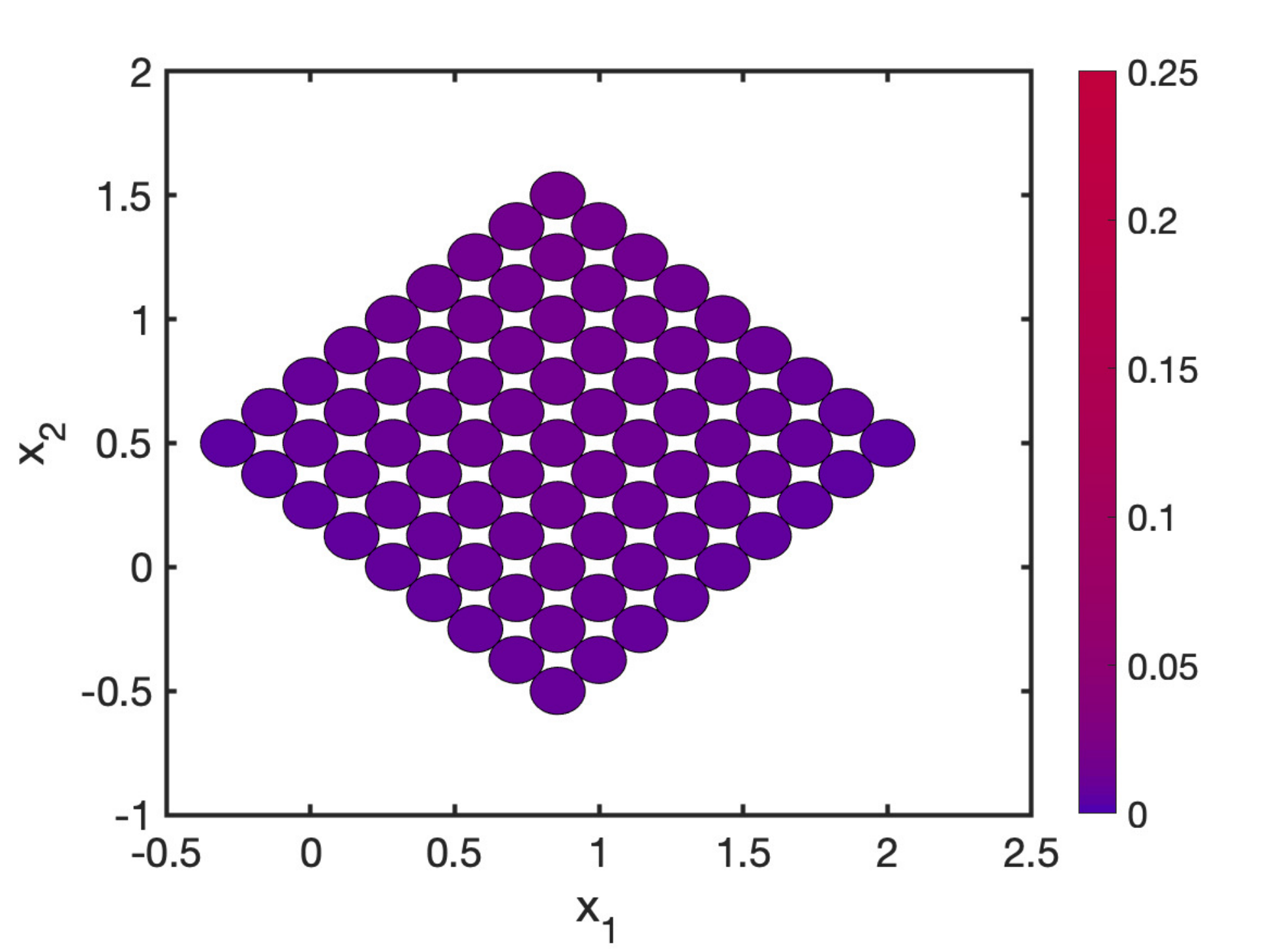}}
	\subfloat[$\xi_I^C=50$] {\includegraphics[width=0.24\textwidth]{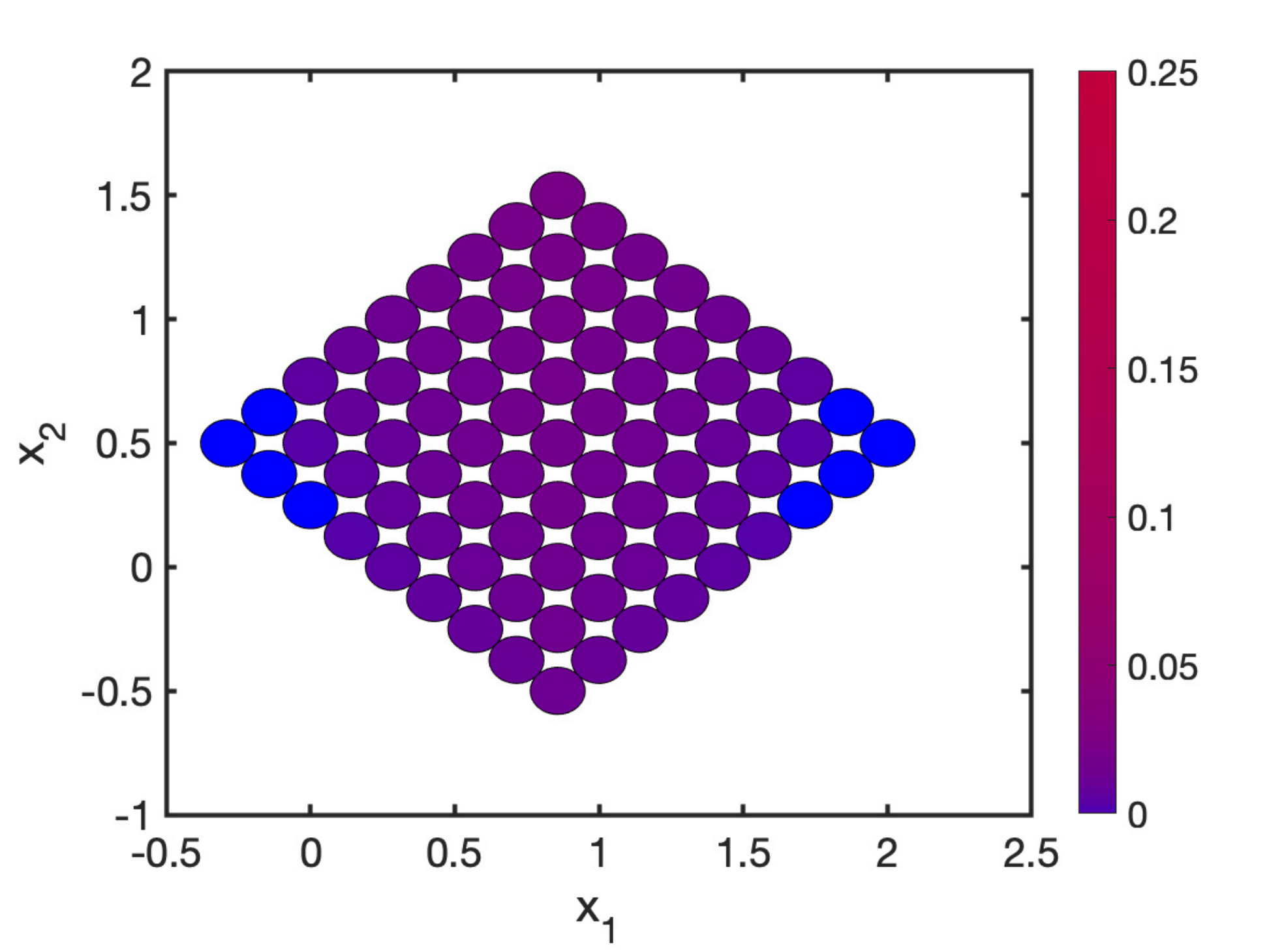}}
	\subfloat[$\xi_I^C=100$] {\includegraphics[width=0.24\textwidth]{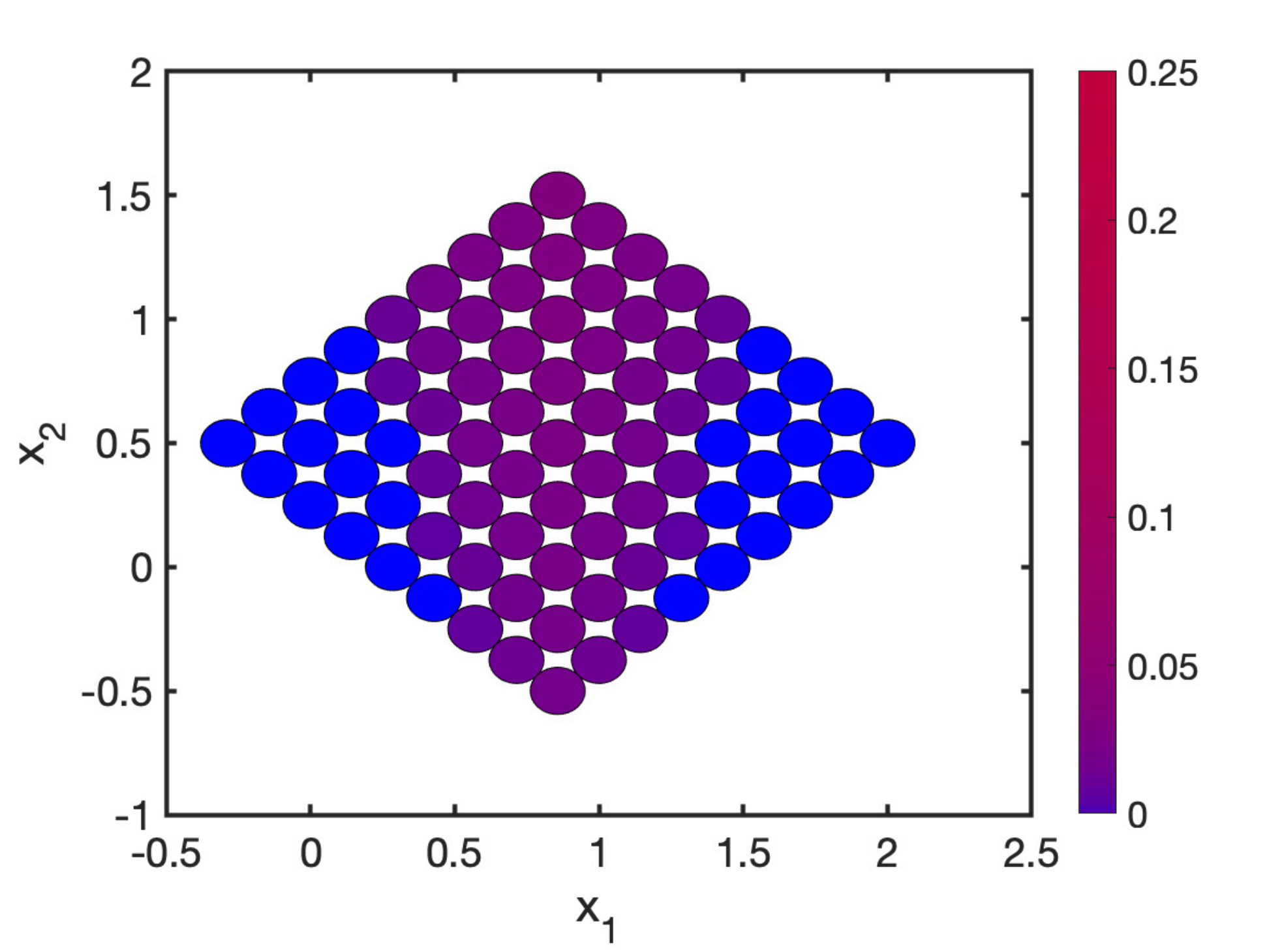}}
	\subfloat[$\xi_I^C=5000$] {\includegraphics[width=0.24\textwidth]{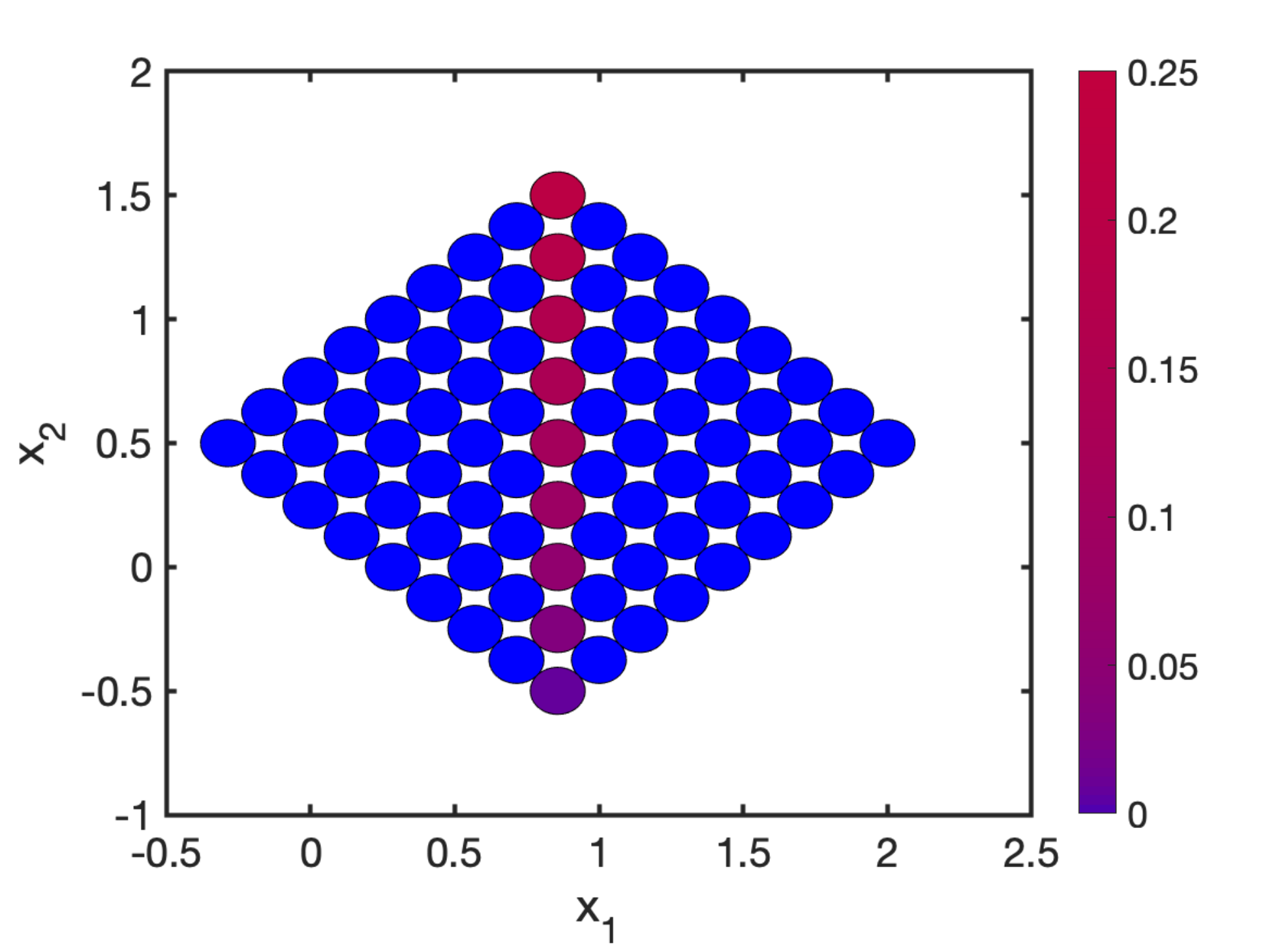}}
	
	\subfloat[$\xi_I^C=10$] {\includegraphics[width=0.24\textwidth]{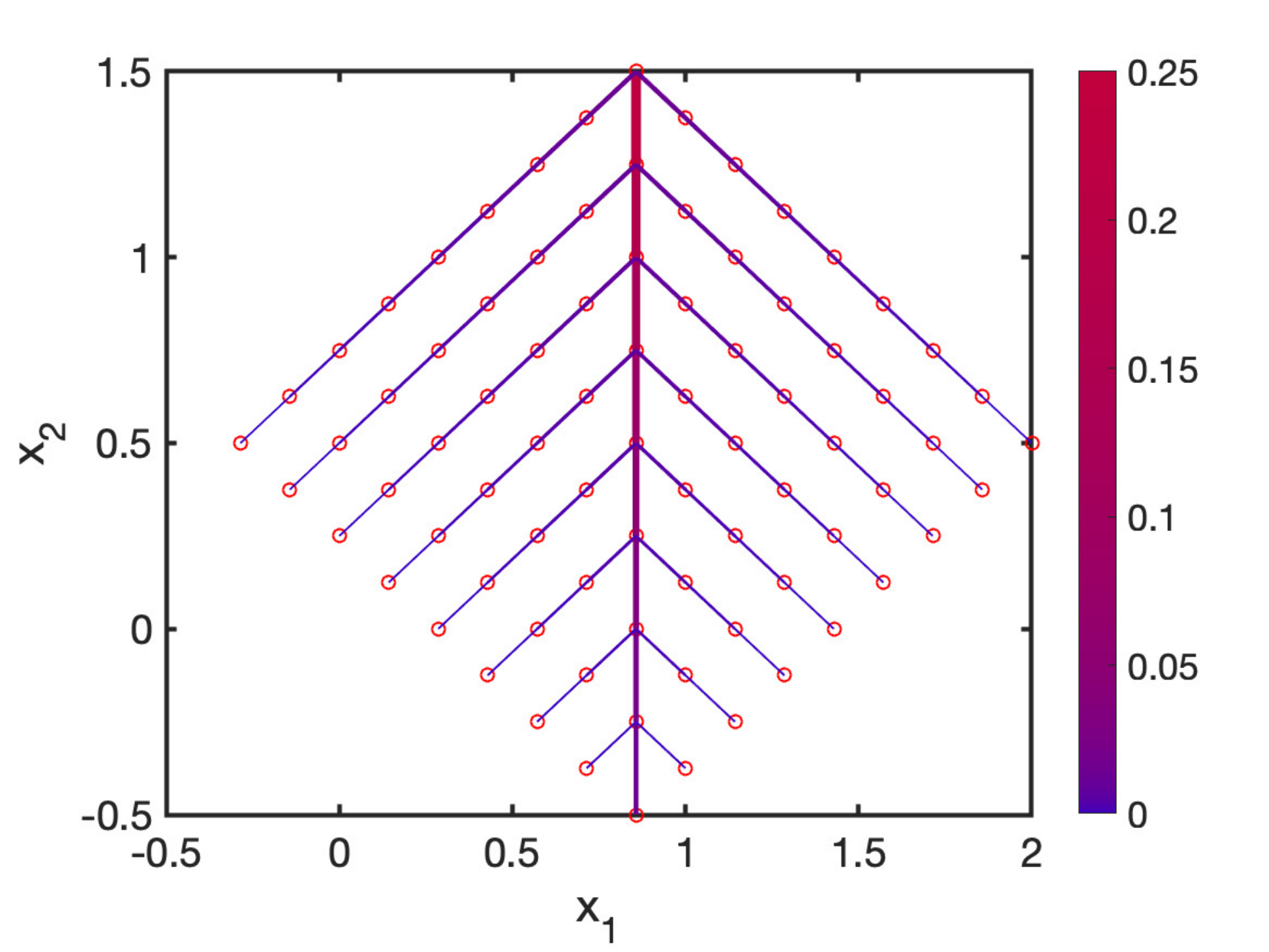}}
	\subfloat[$\xi_I^C=50$] {\includegraphics[width=0.24\textwidth]{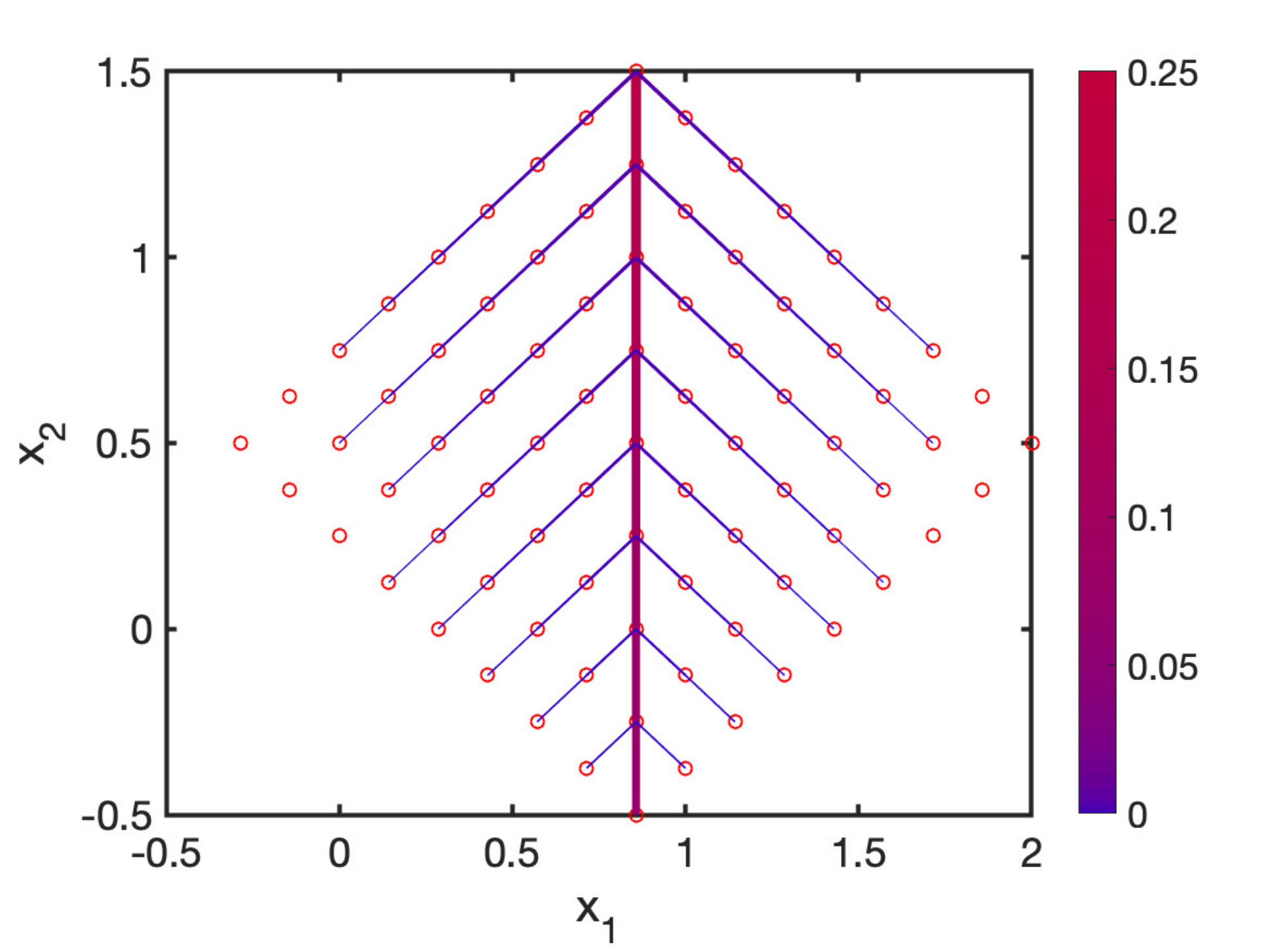}}
	\subfloat[$\xi_I^C=100$] {\includegraphics[width=0.24\textwidth]{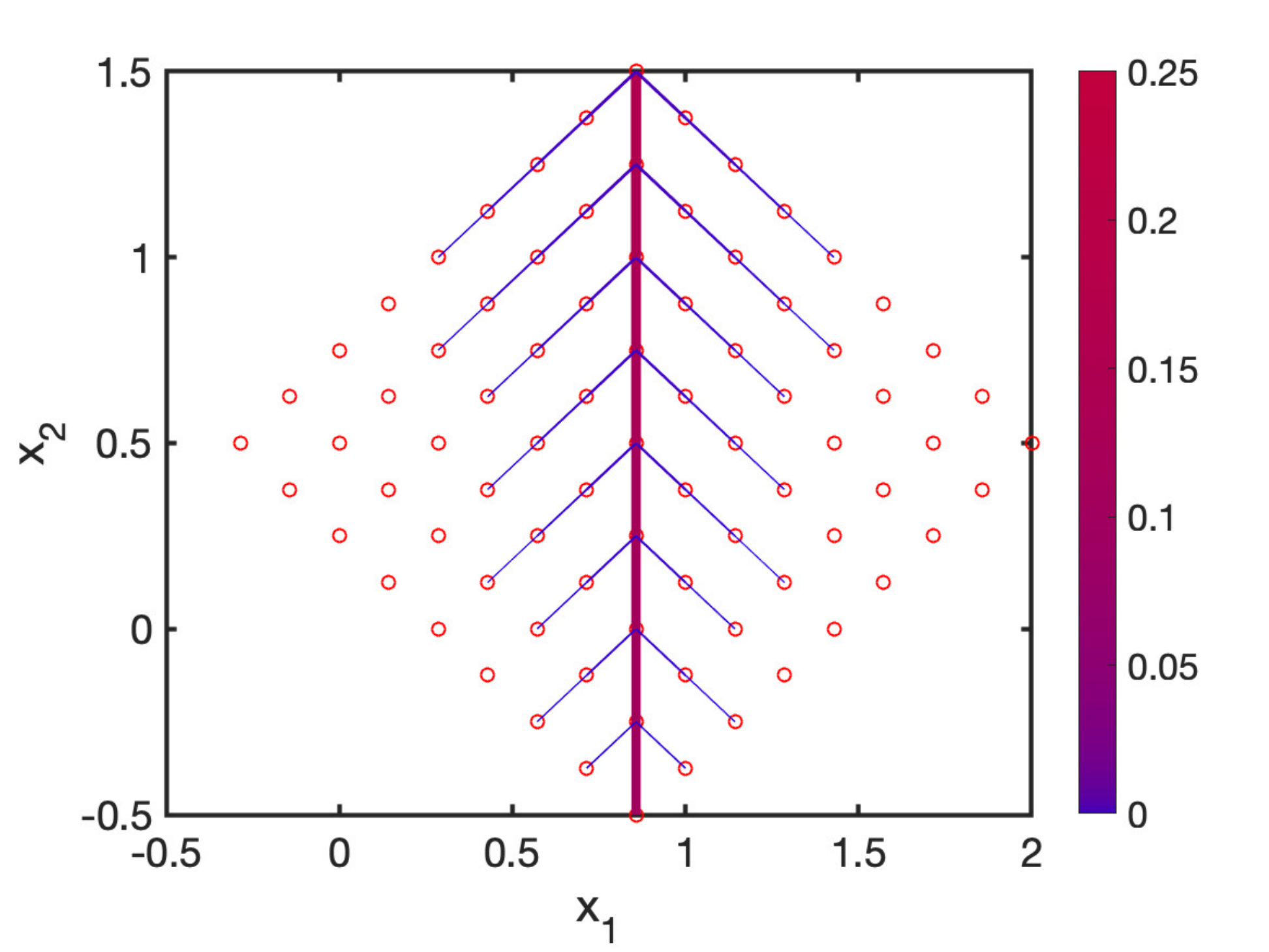}}
	\subfloat[$\xi_I^C=5000$] {\includegraphics[width=0.24\textwidth]{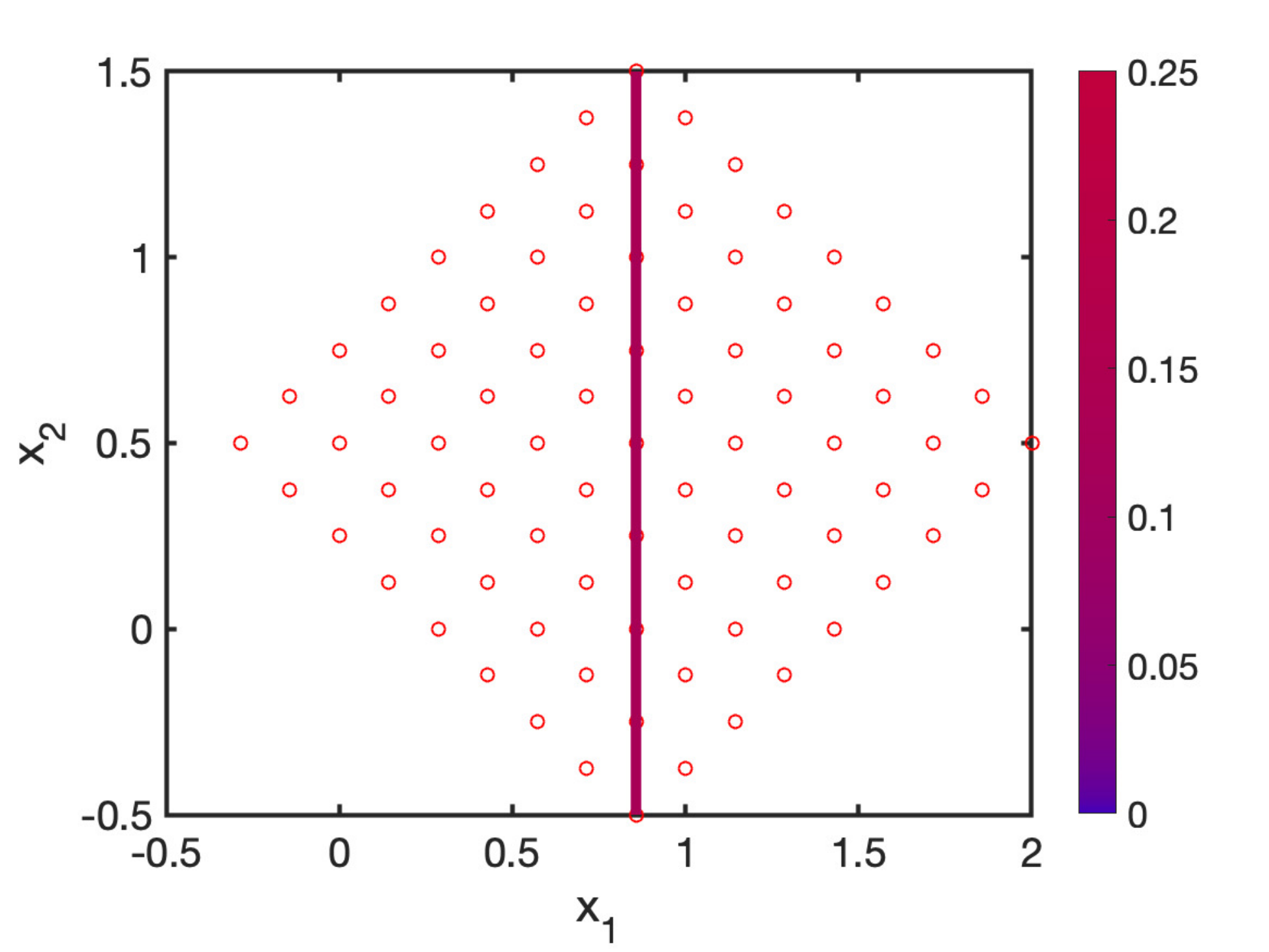}}
	\caption{Steady states  for auxin concentration and transport activity for different sink strengths $\xi_I^C$ with initial data $\bar{X}_{ij}, \bar{a}_{i}$.}\label{fig:sinkstrengthcorner}
\end{figure}

In Figure \ref{fig:sinkstrengthcorner}, we vary the strength of the sink in the bottom corner, denoted by $\xi_I^C$, while keeping the values of $I_i$ for all other vertices $i\in\Vset$ as before. Similarly as for the  variation of $\xi_I$, the area of the network decreases as $\xi_I^C$ increases for both auxin levels and transport activity. In this case, however, it decreases outside a neighborhood of the line connecting the source in the top corner and the increasing sink of size $\xi_I^C$ in the bottom corner. In particular, the network structure for large $\xi_S^C$ is given by a high auxin levels and transport activity along the line of cells, connecting the source in the top corner with the strong sink in the bottom corner. Moreover, this variation of the size of the source $\xi_S$ in Figures \ref{fig:sourcestrength} and \ref{fig:sourcestrengthgrids}, as well as,  of the sinks $\xi_I$ and $\xi_I^C$ in Figure \ref{fig:sinkstrengthcorner}  illustrate how crucial the choice of sources and sinks for the resulting pattern formation is.

\begin{figure}[htbp]
	\centering
	
	\subfloat[$\delta=0.1$] {\includegraphics[width=0.24\textwidth]{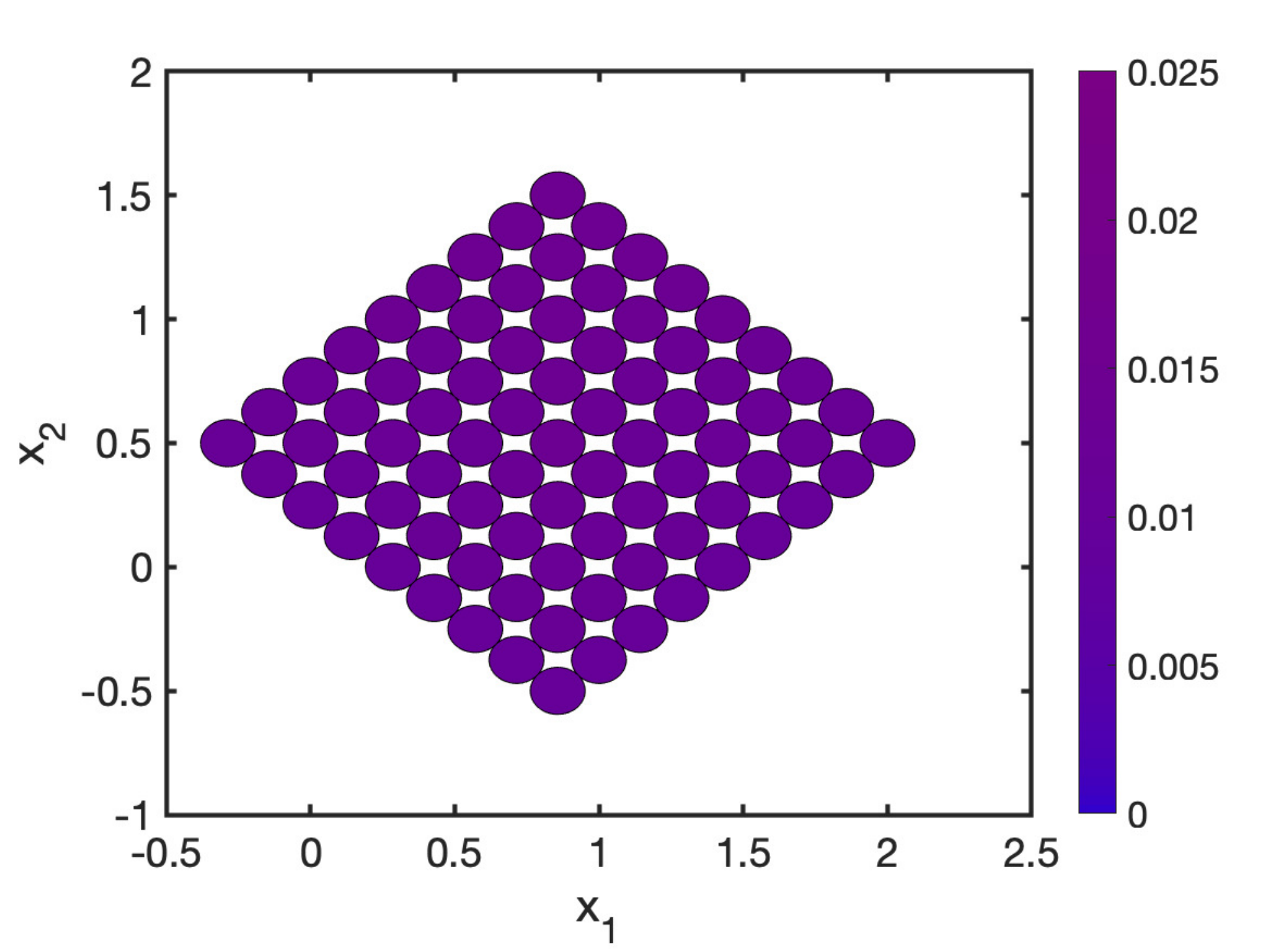}}
	\subfloat[$\delta=0.5$] {\includegraphics[width=0.24\textwidth]{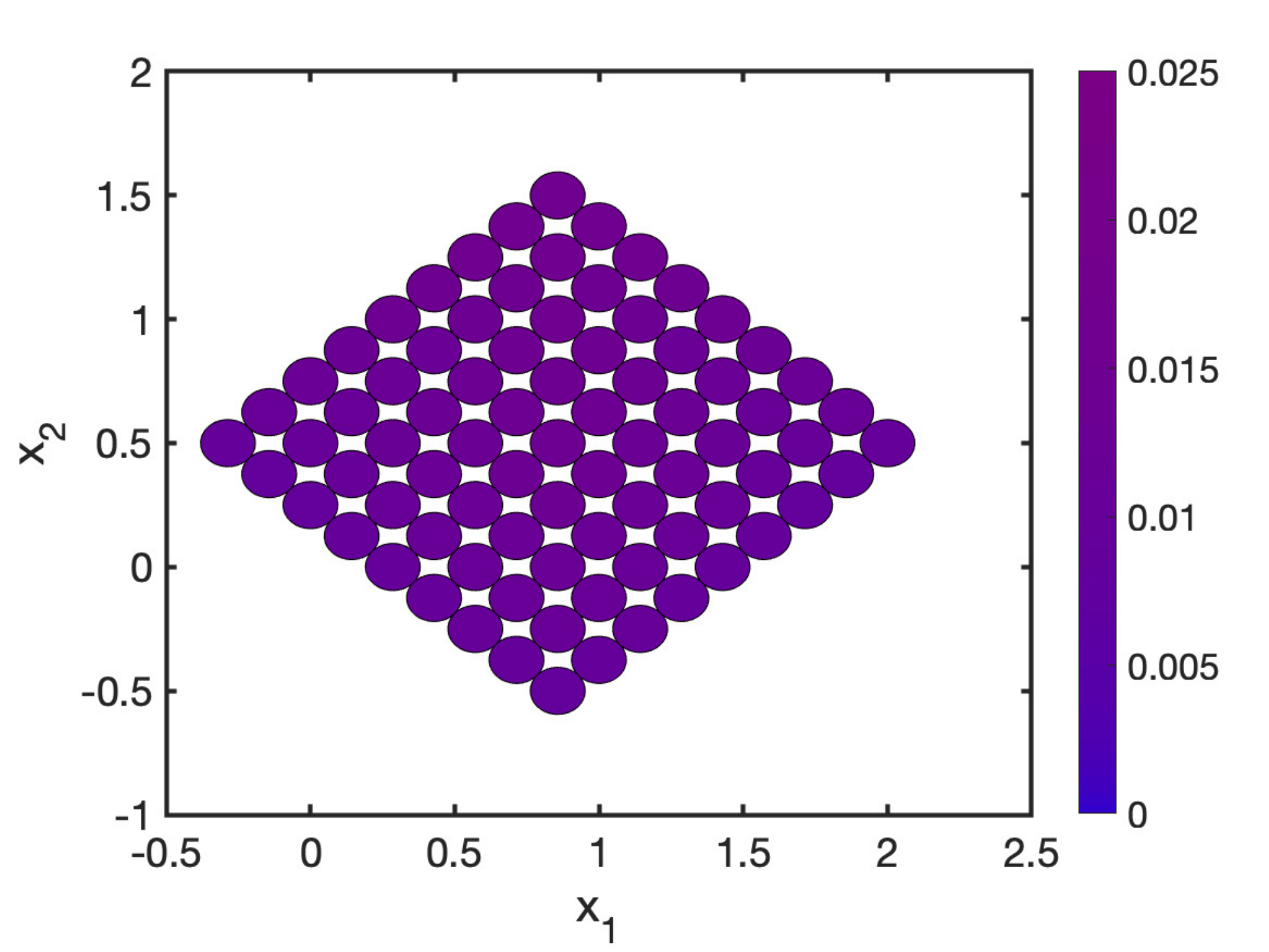}}
	\subfloat[$\delta=2$] {\includegraphics[width=0.24\textwidth]{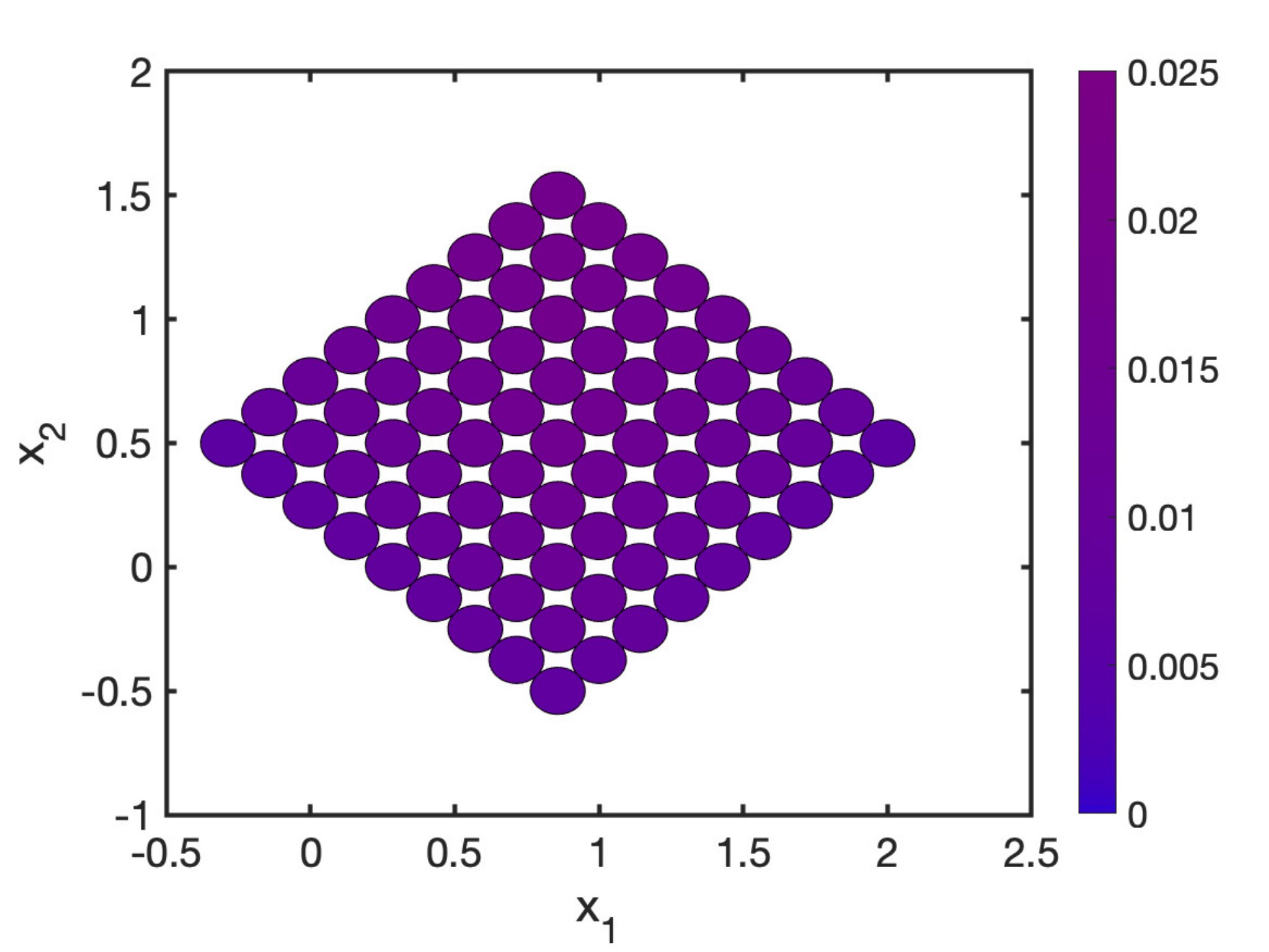}}
	\subfloat[$\delta=10$] {\includegraphics[width=0.24\textwidth]{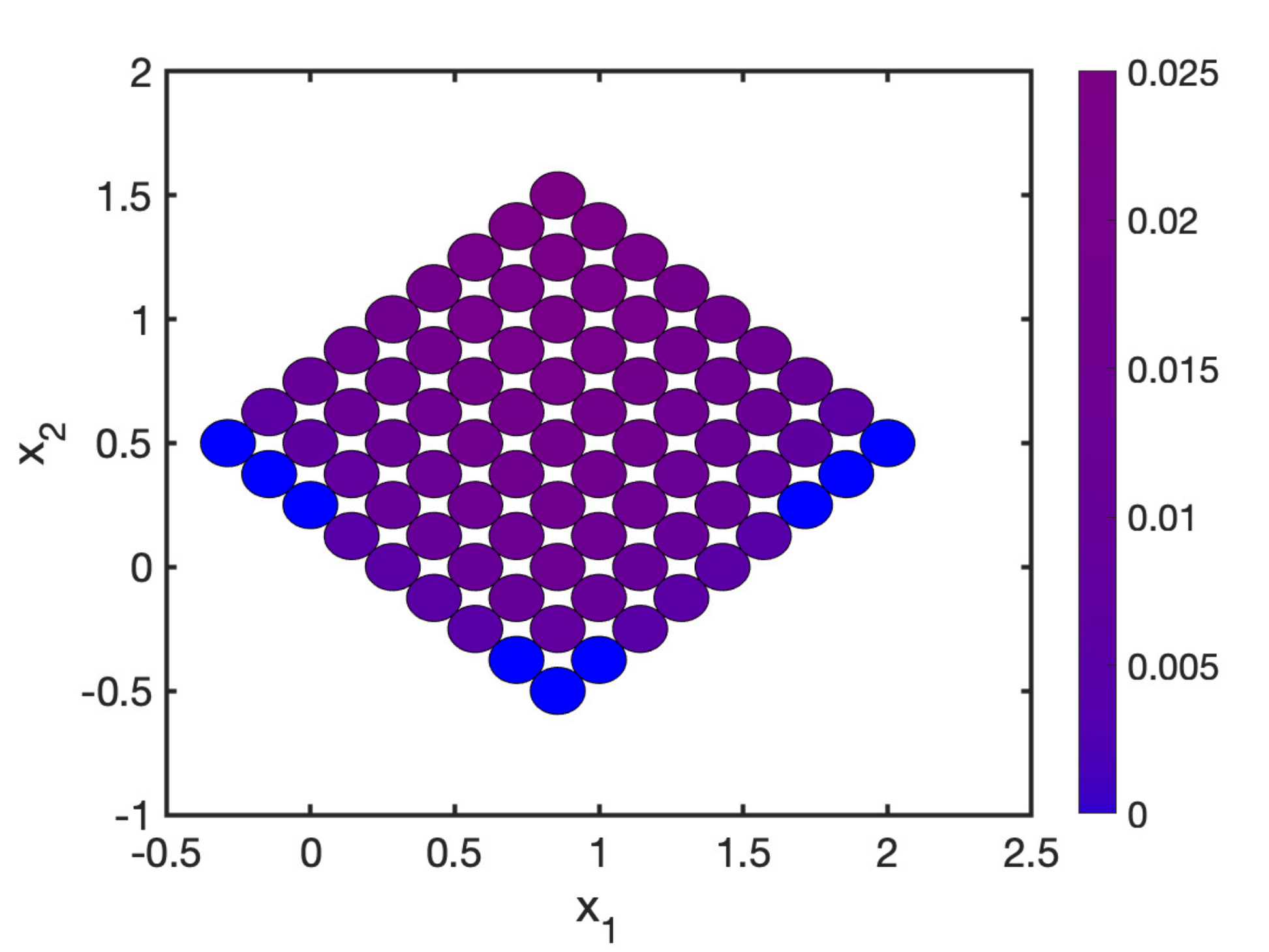}}
	
	\subfloat[$\delta=0.1$] {\includegraphics[width=0.24\textwidth]{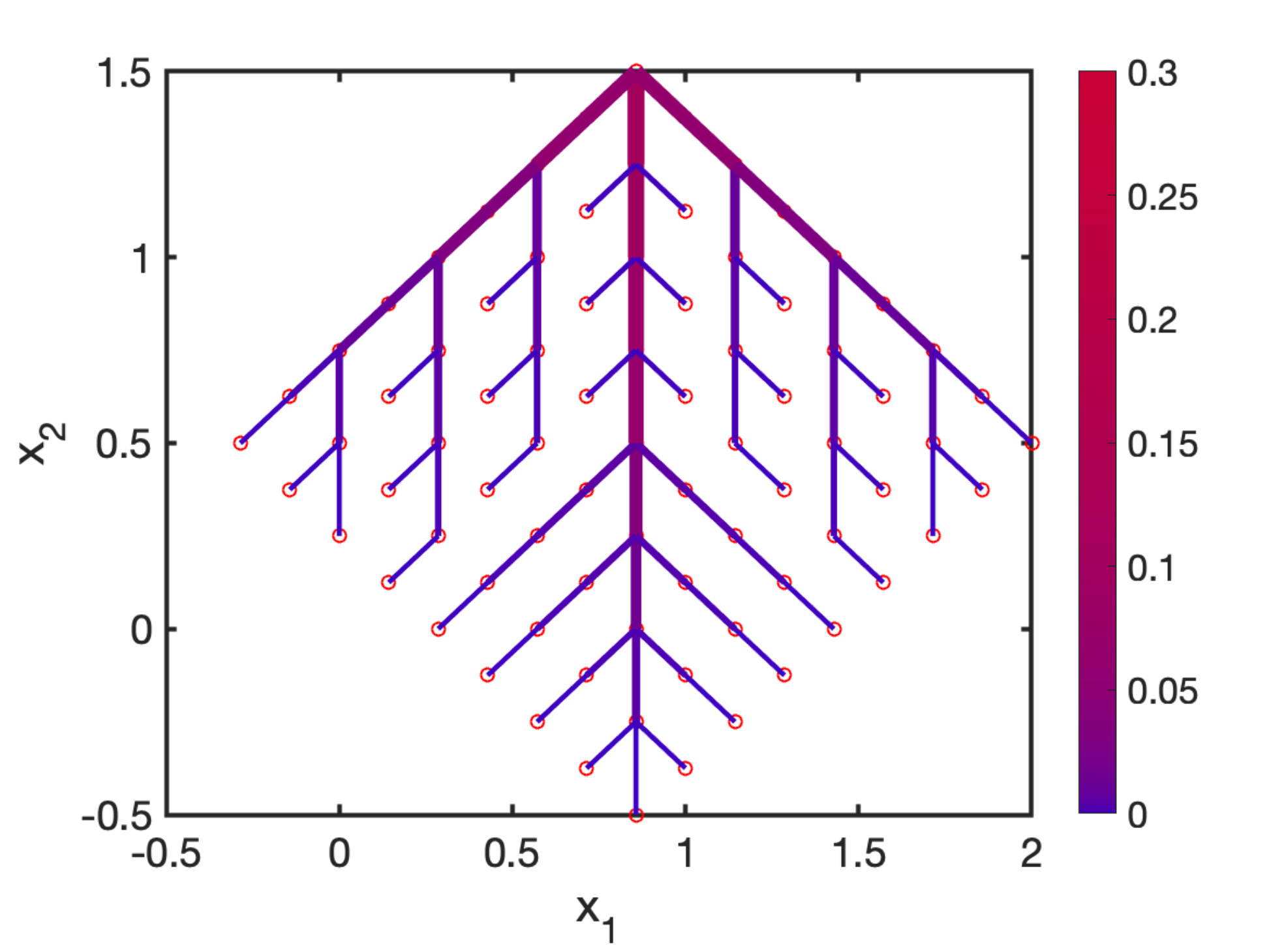}}
	\subfloat[$\delta=0.5$] {\includegraphics[width=0.24\textwidth]{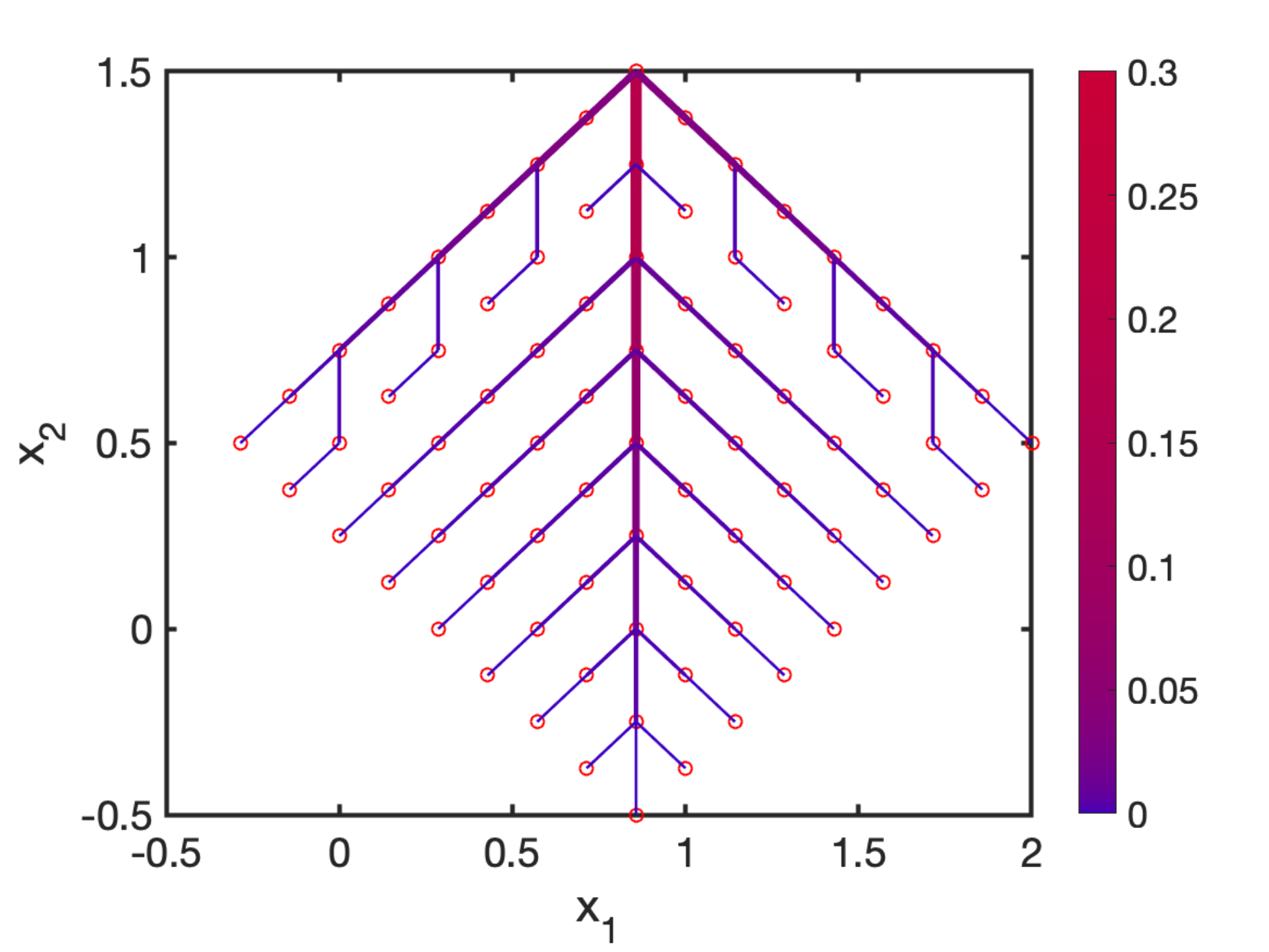}}
	\subfloat[$\delta=2$] {\includegraphics[width=0.24\textwidth]{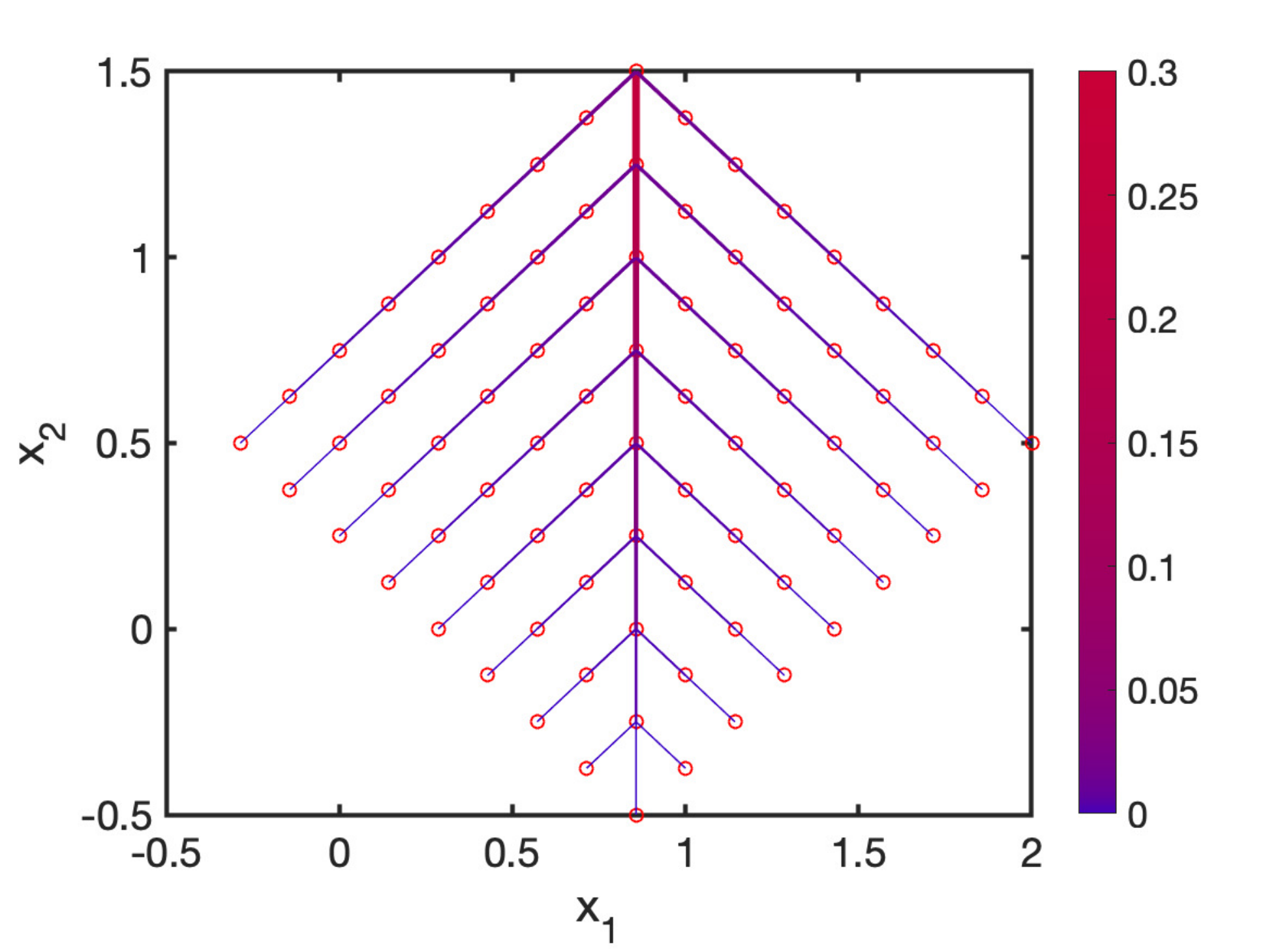}}
	\subfloat[$\delta=10$] {\includegraphics[width=0.24\textwidth]{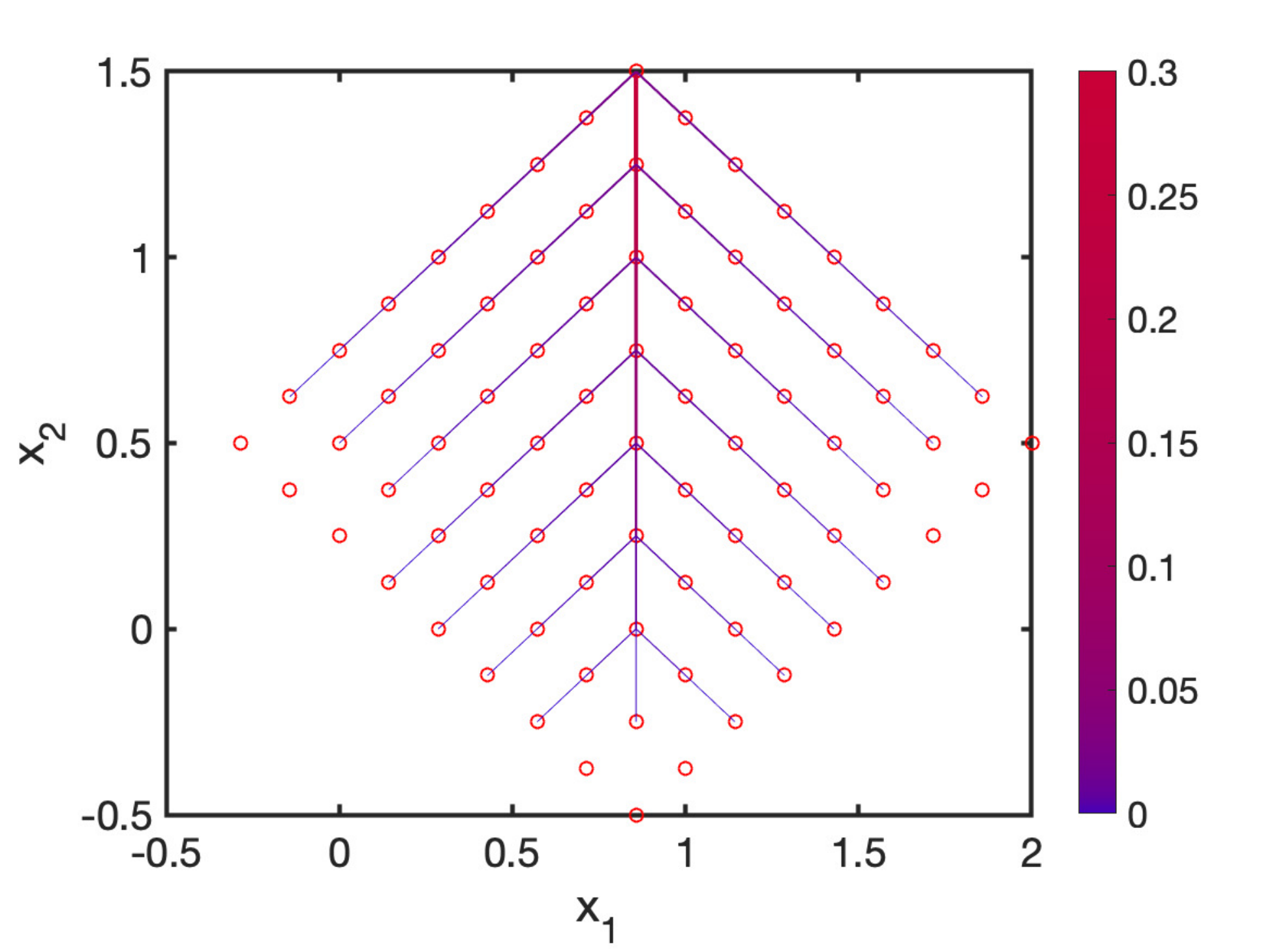}}
	
	\caption{Steady states  for auxin transport activity and auxin levels  for different parameter values $\delta$ with initial data $\bar{X}_{ij}, \bar{a}_{i}$.}\label{fig:deltadependence}
\end{figure}

\begin{figure}[htbp]
	\centering
	
	\subfloat[$\tau=0.5$] {\includegraphics[width=0.24\textwidth]{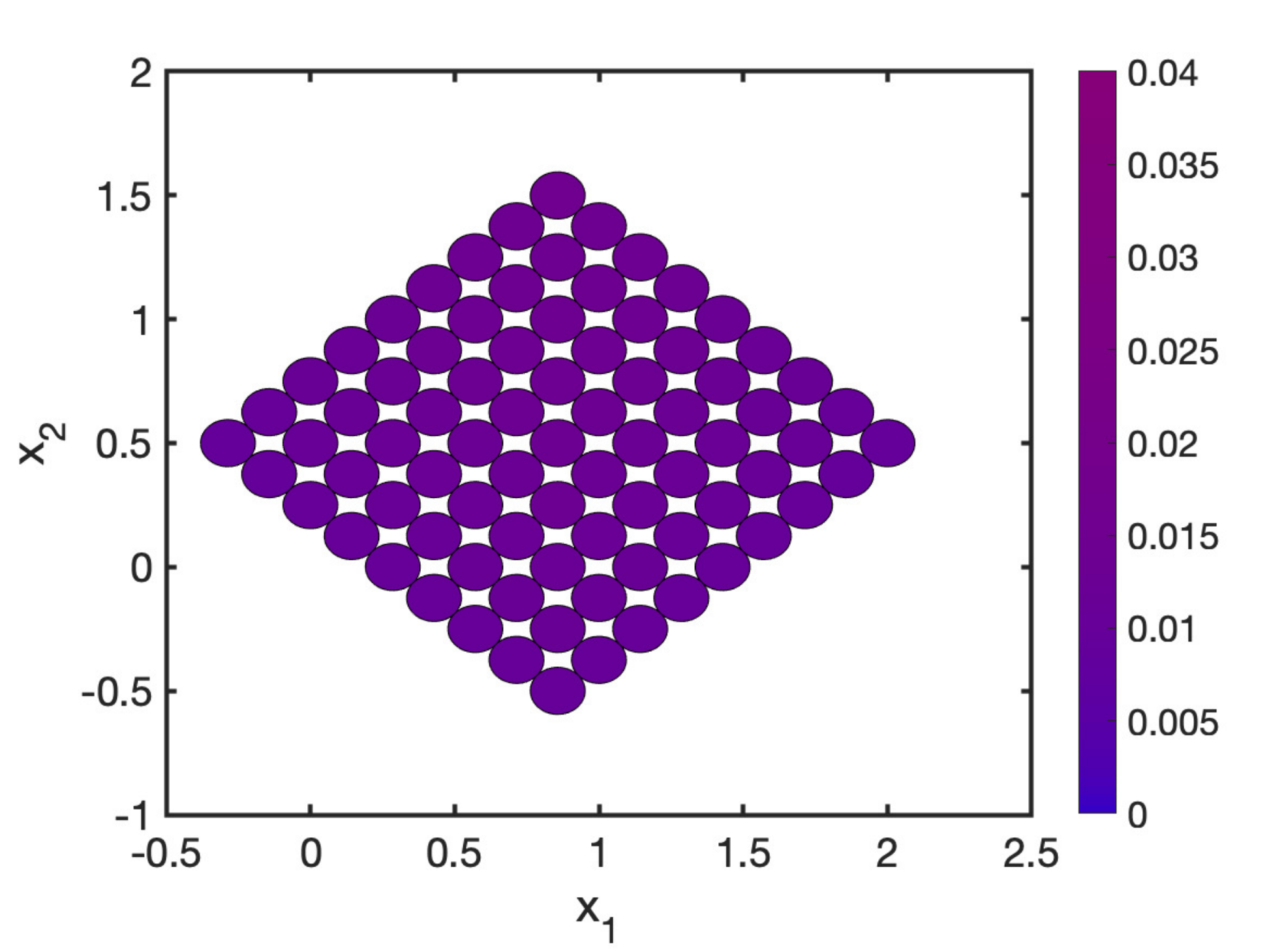}}
	\subfloat[$\tau=2$] {\includegraphics[width=0.24\textwidth]{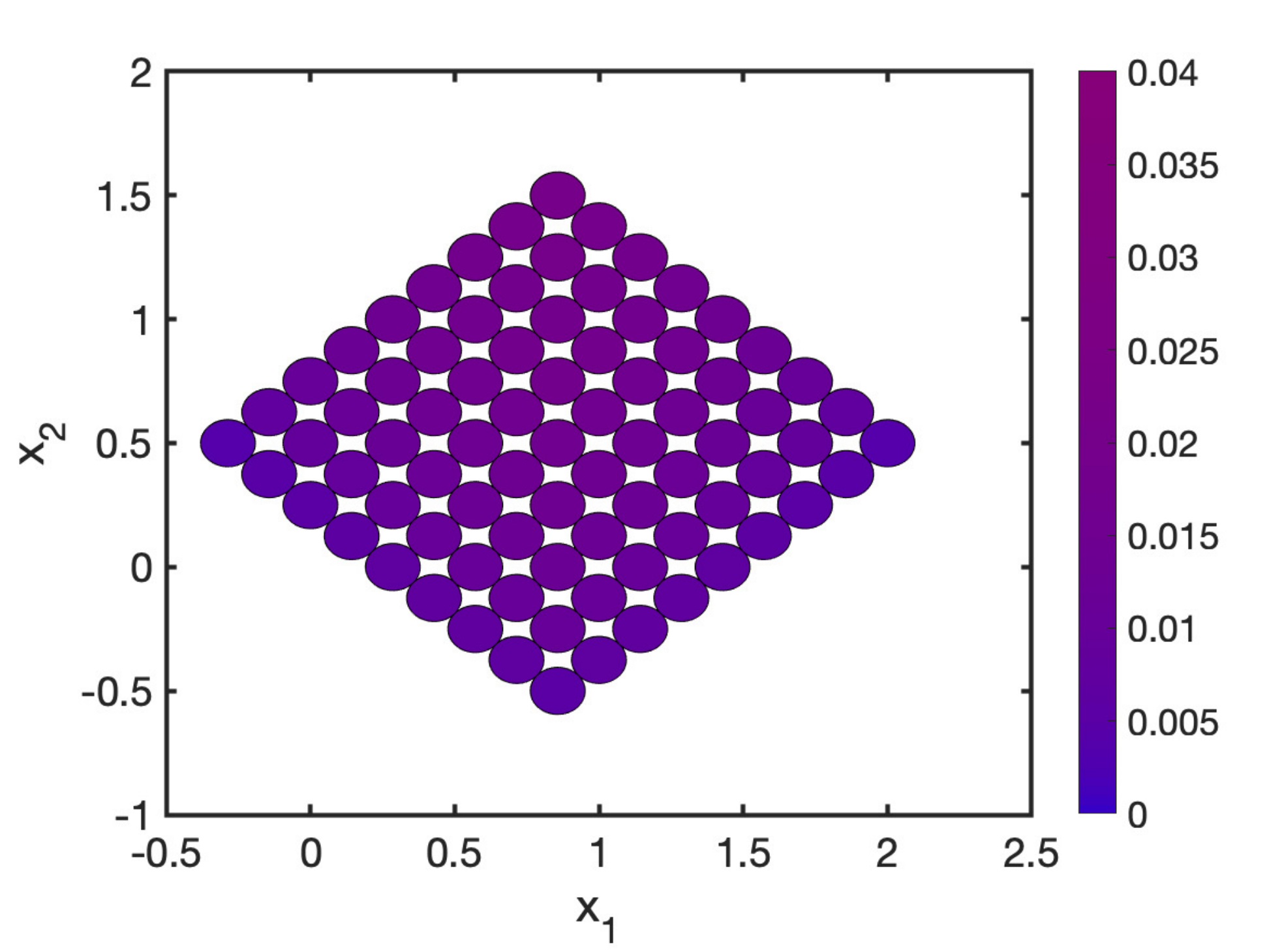}}
	\subfloat[$\tau=5$] {\includegraphics[width=0.24\textwidth]{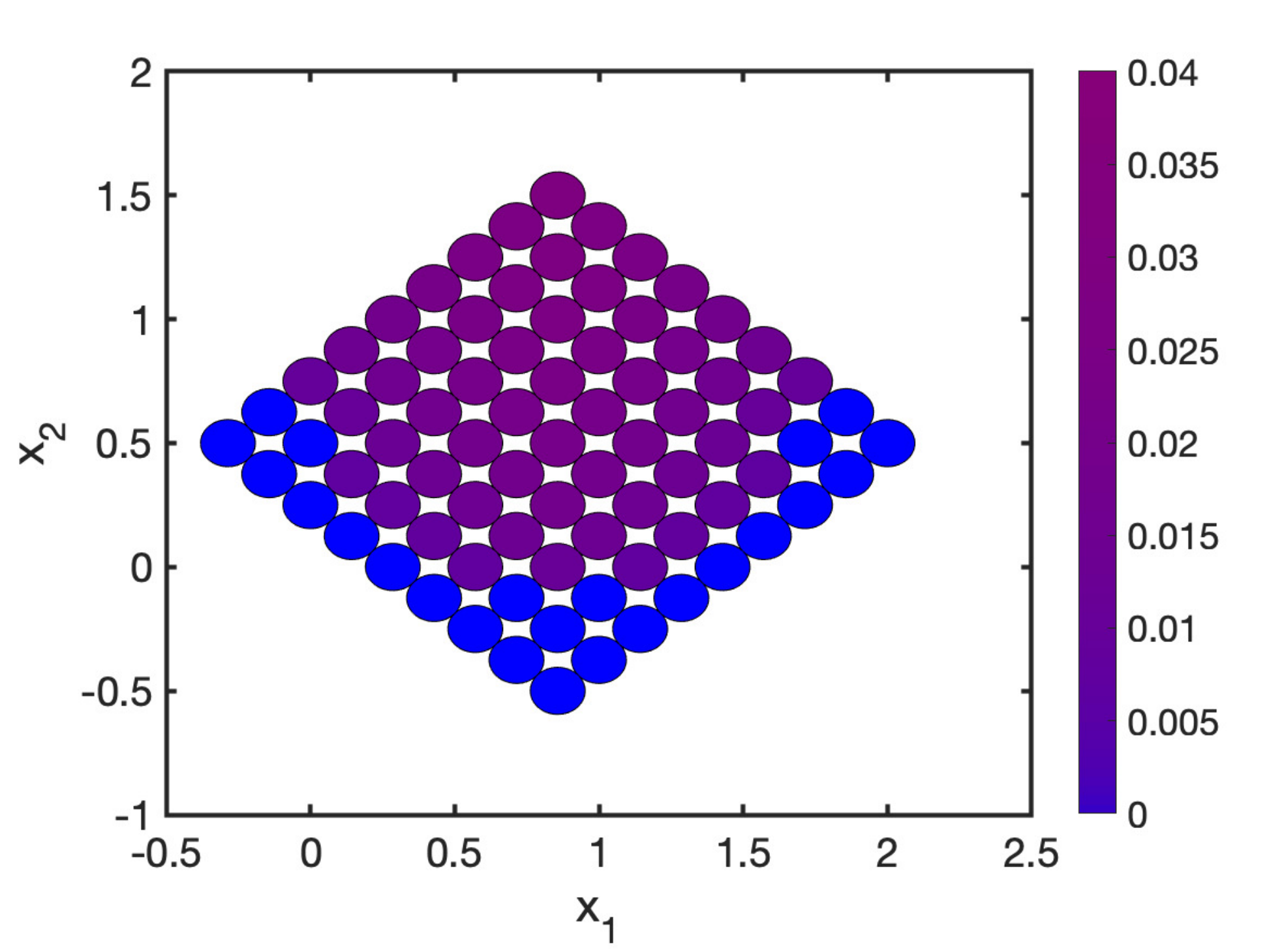}}
	\subfloat[$\tau=10$] {\includegraphics[width=0.24\textwidth]{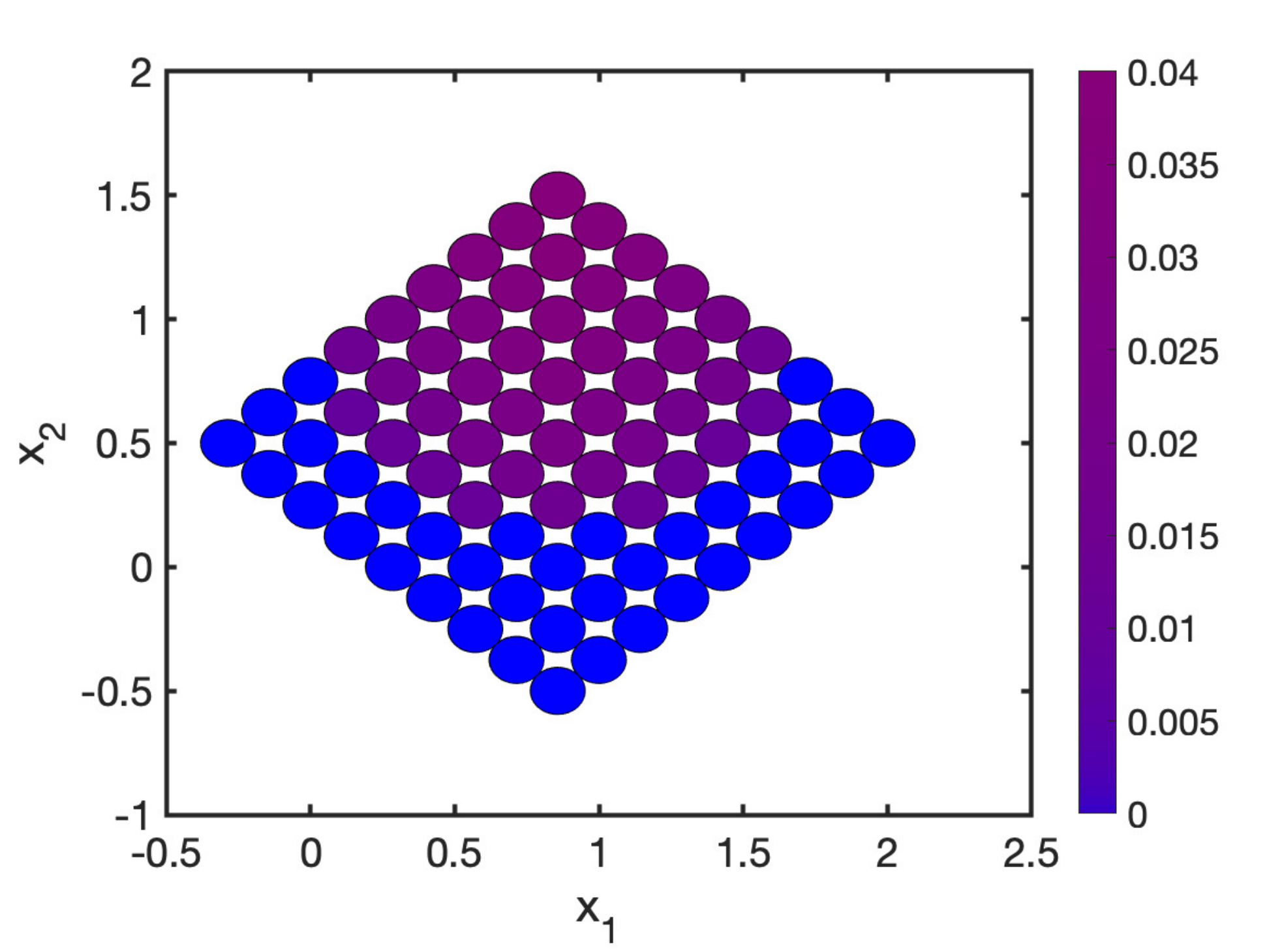}}
	
	\subfloat[$\tau=0.5$] {\includegraphics[width=0.24\textwidth]{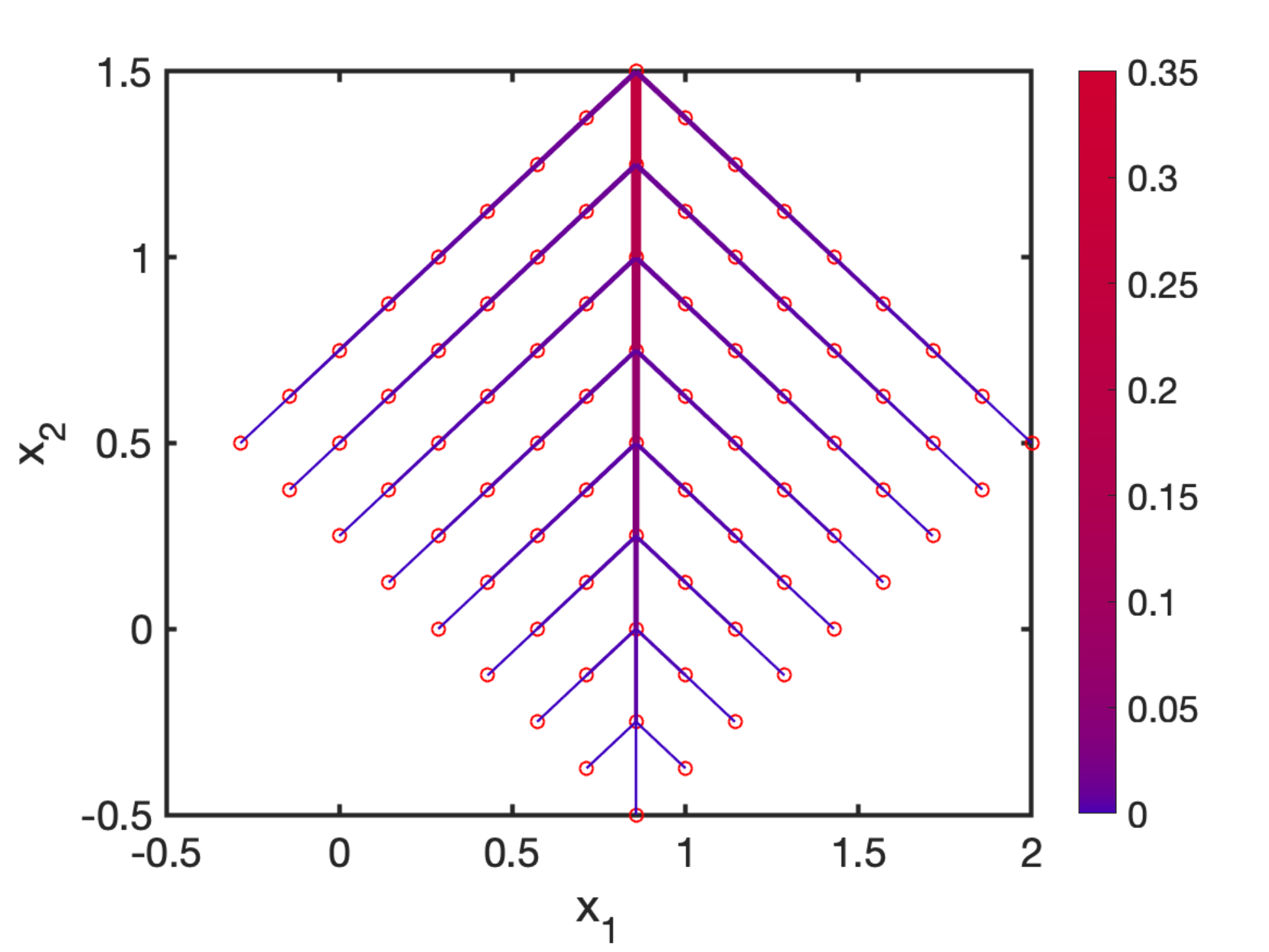}}
	\subfloat[$\tau=2$] {\includegraphics[width=0.24\textwidth]{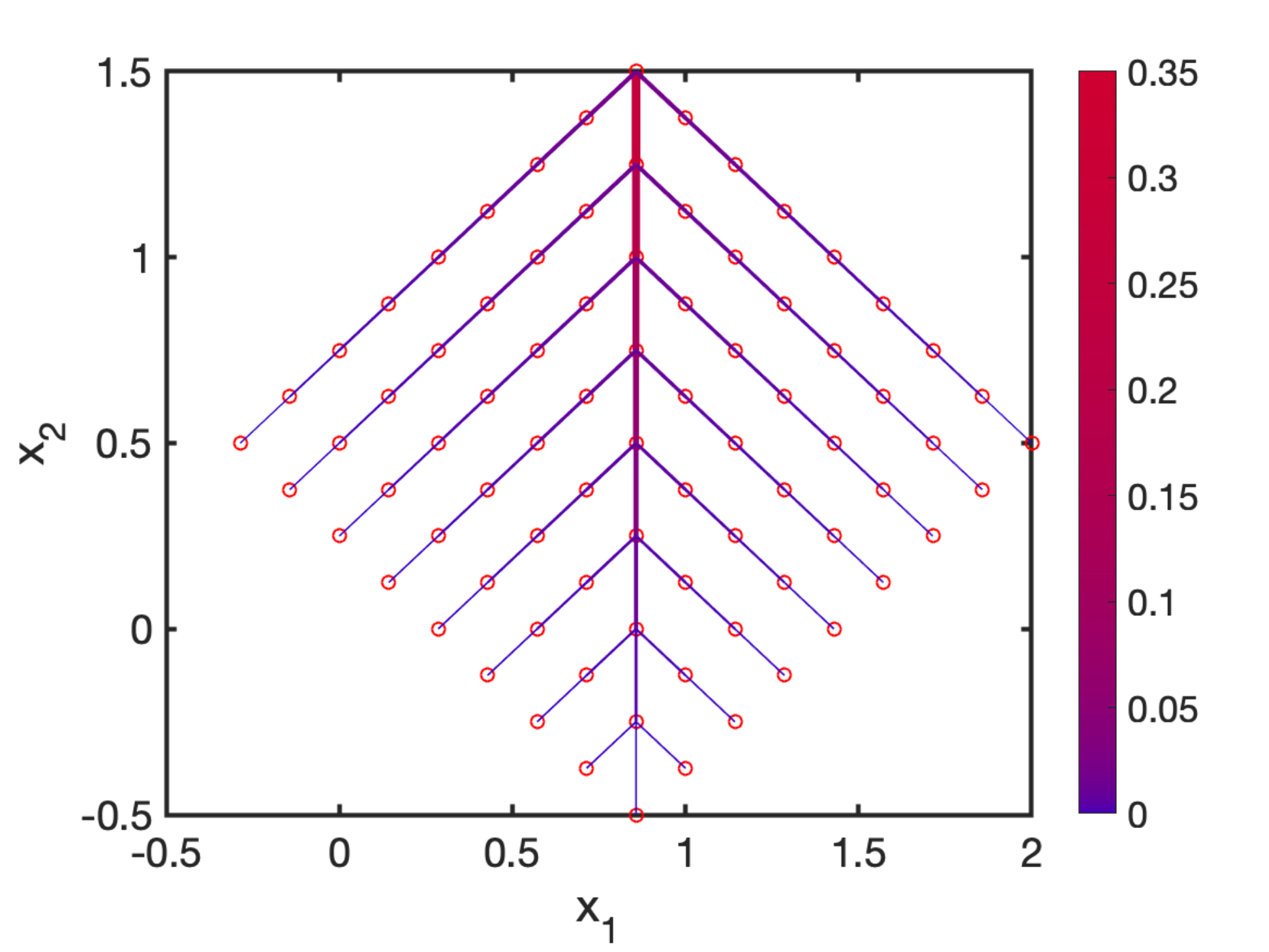}}
	\subfloat[$\tau=5$] {\includegraphics[width=0.24\textwidth]{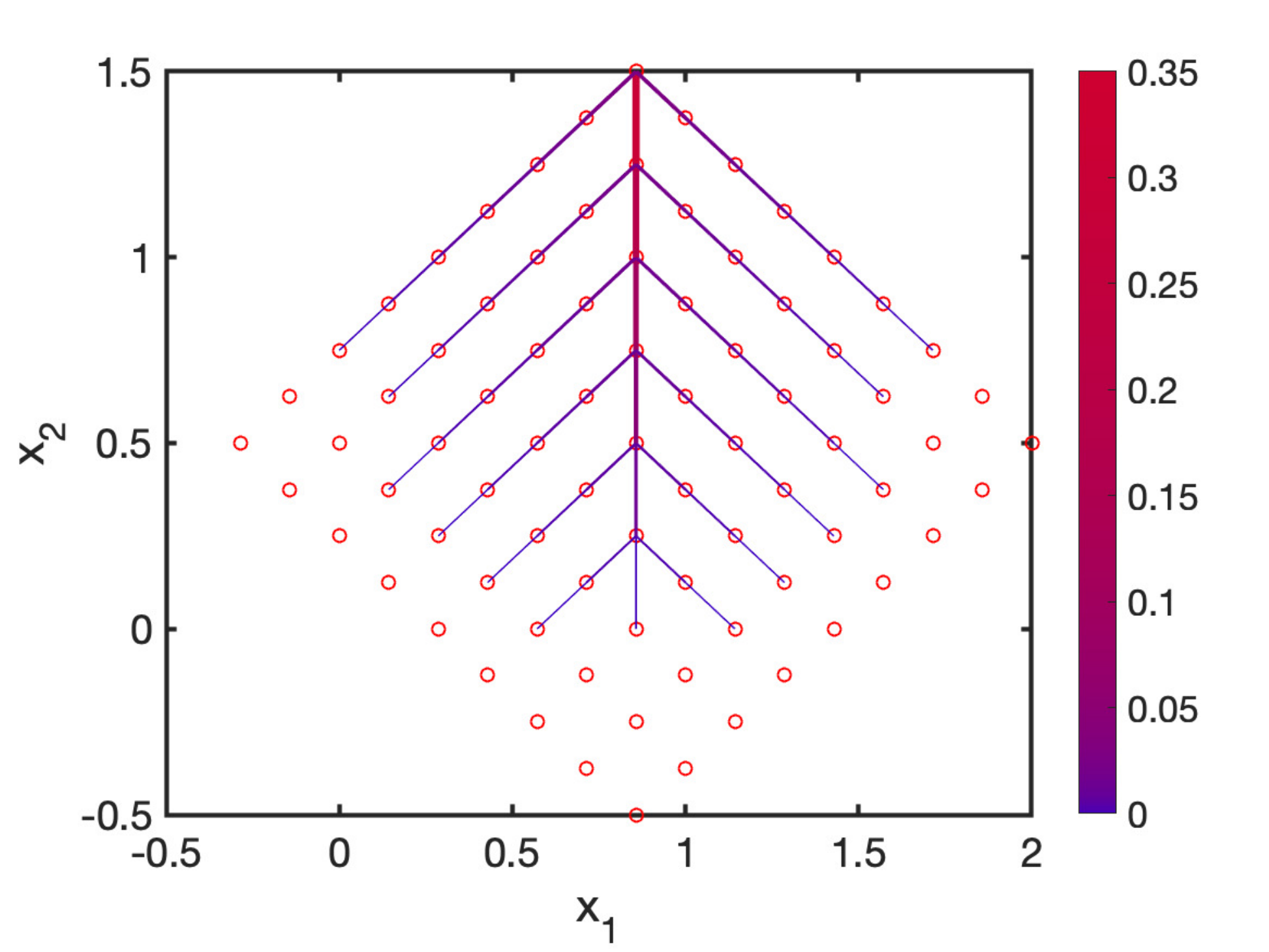}}
	\subfloat[$\tau=10$] {\includegraphics[width=0.24\textwidth]{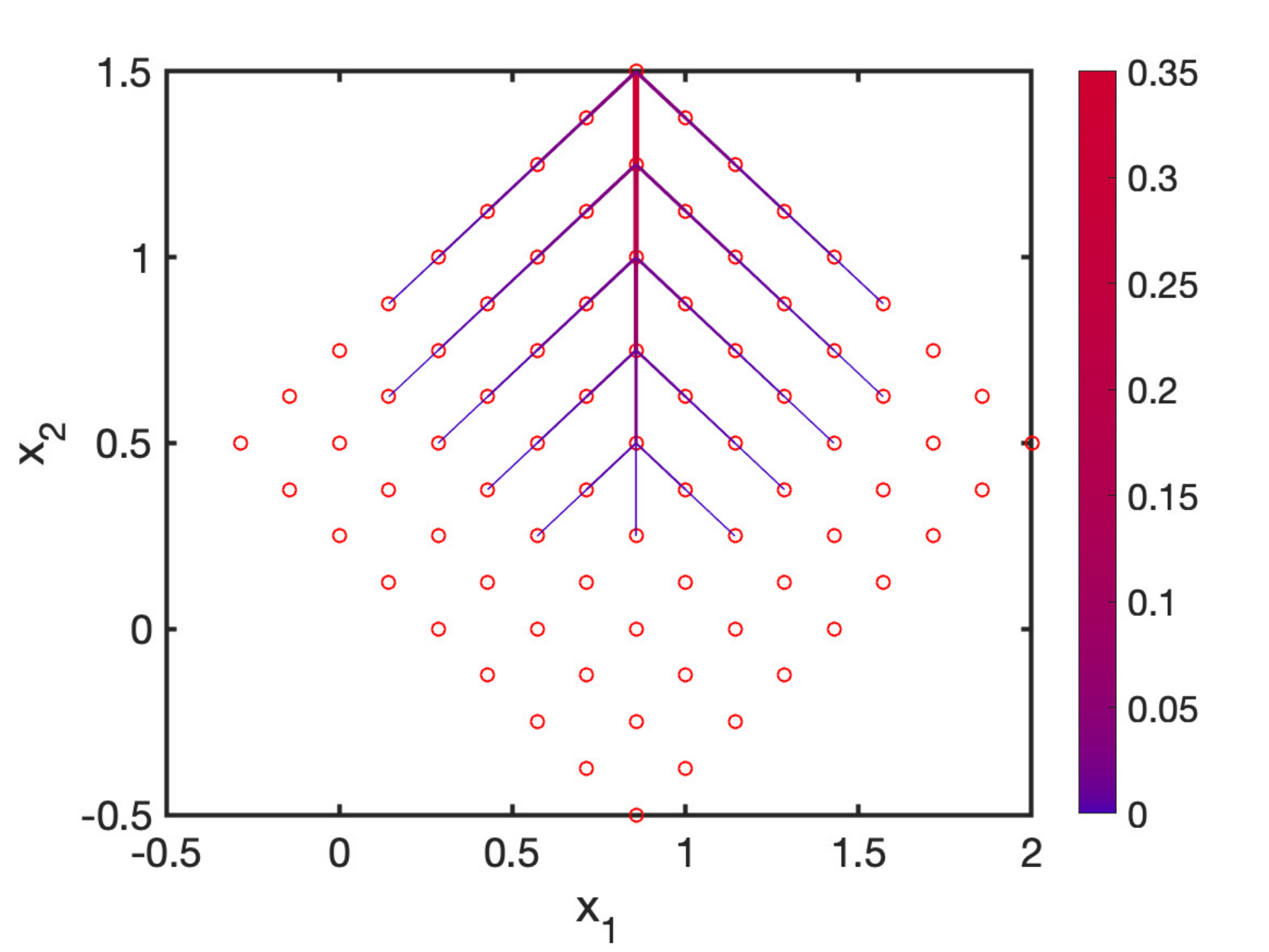}}
	
	\caption{Steady states  for auxin transport activity and auxin levels for different parameter values $\tau$ with initial data $\bar{X}_{ij}, \bar{a}_{i}$.}\label{fig:taudependence} 
\end{figure}

In Figures \ref{fig:deltadependence} and \ref{fig:taudependence}, we investigate the dependence of the stationary states on the model parameters $\delta$ and $\tau$ in \eqref{eq:auxineqproduction}--\eqref{eq:pineq}. For small values of $\delta$, more complex stationary patterns for the transport activity  can be seen in Figure \ref{fig:deltadependence} and auxin is transported over the entire graph. As $\delta$ increases, the auxin levels and the transport activity  increase close to the source, but they are no longer transported over the entire graph. As before, the area covered by  auxin transport activity and auxin levels are of a similar size, i.e., auxin transport activity and auxin levels  are  co-existent. The increase of $\tau$ shows a similar change of the steady states of both the auxin transport activity and auxin levels as the increase of $\delta$.

\begin{figure}[htbp]
	\centering
	\subfloat[$\epsilon=0.5$] {\includegraphics[width=0.24\textwidth]{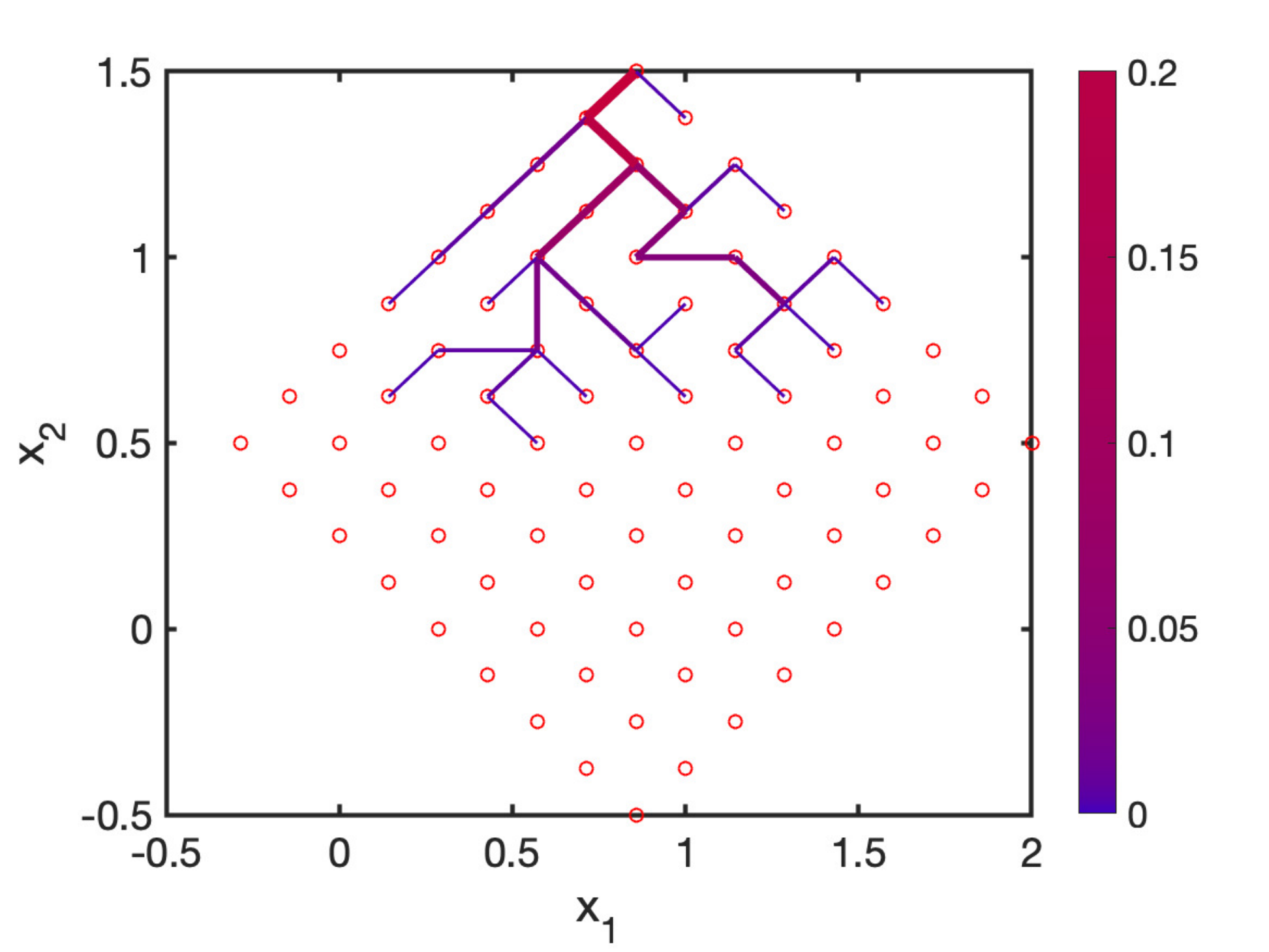}}
	\subfloat[$\epsilon=5$] {\includegraphics[width=0.24\textwidth]{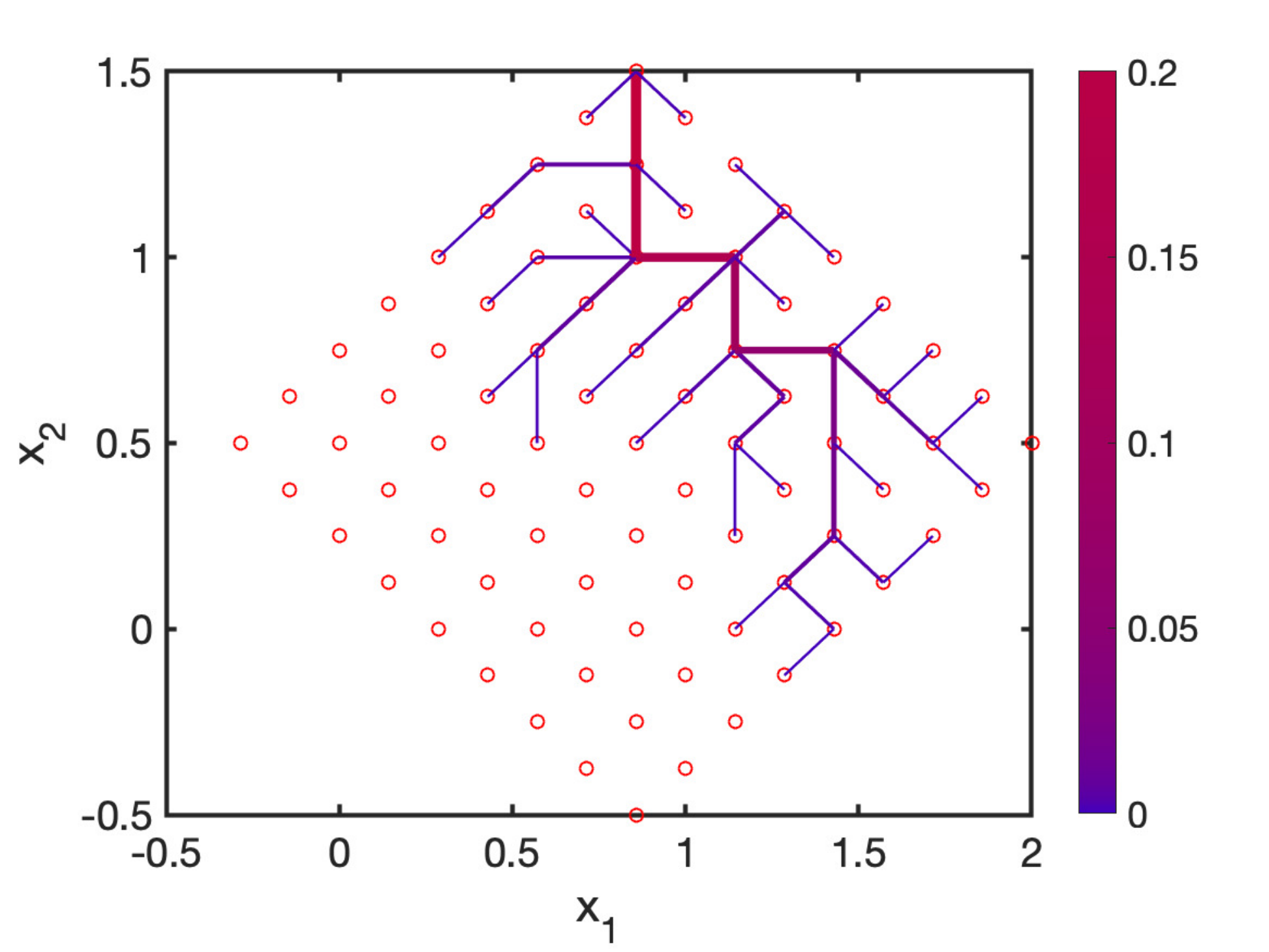}}
	\subfloat[$\epsilon=50$] {\includegraphics[width=0.24\textwidth]{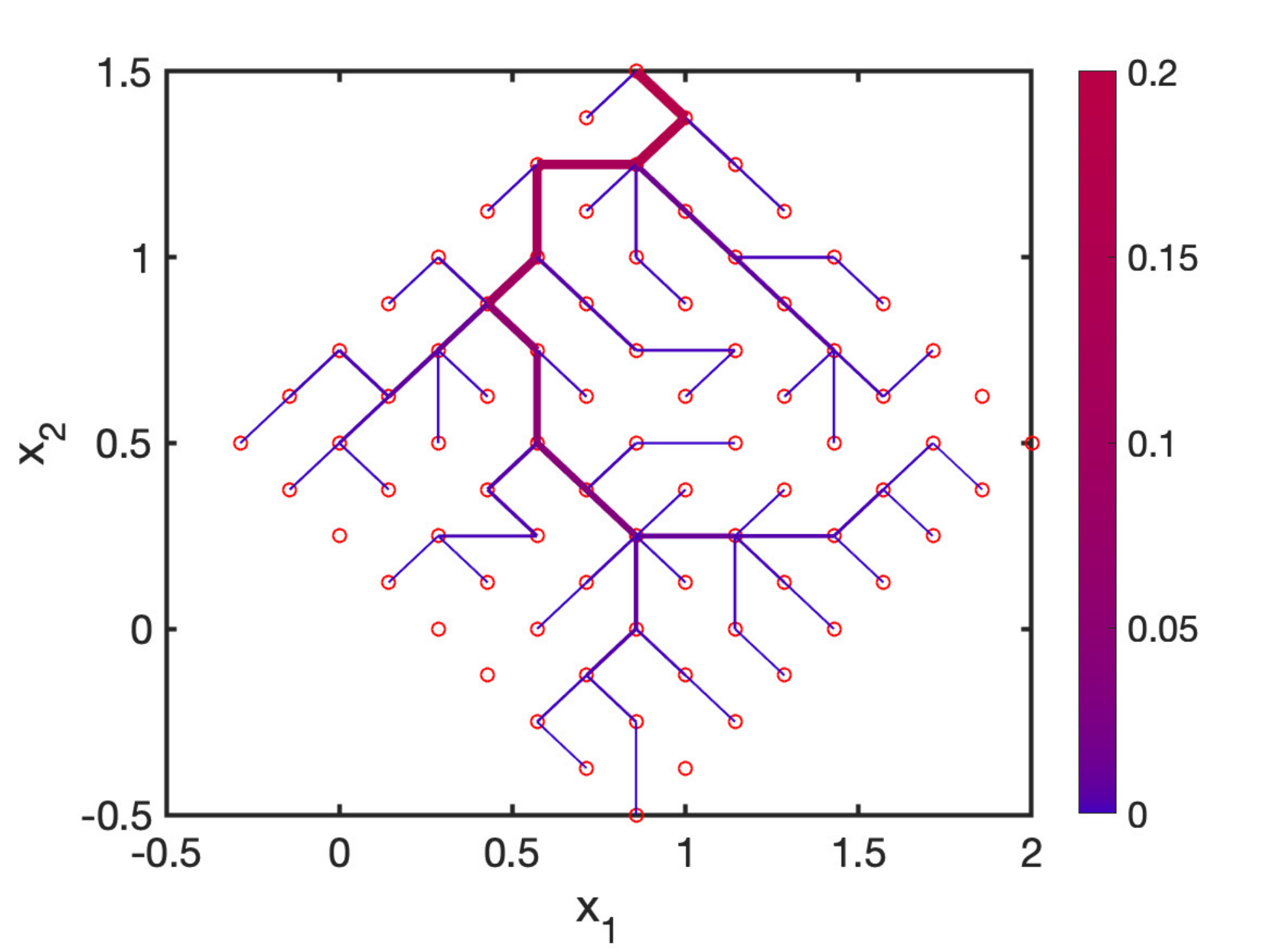}}
	\subfloat[$\epsilon=100$] {\includegraphics[width=0.24\textwidth]{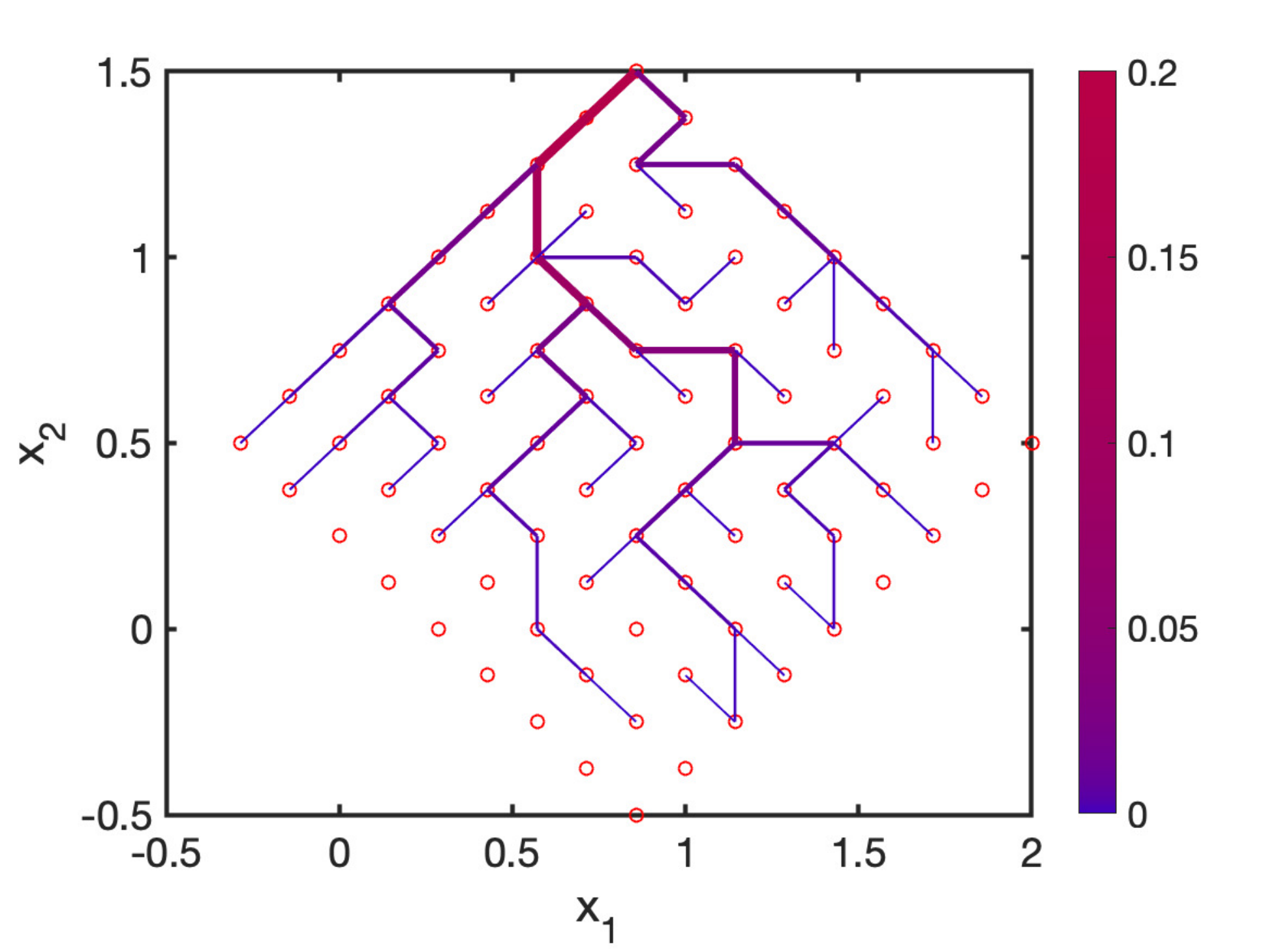}}
	\caption{Steady states for  the transport activity  for initial transport activity $\theta+ 0.00001\epsilon$ where $\theta$ is a random variable with $\theta=1$ with probability $0.2$ and $\theta=0$ with probability $0.8$.}\label{fig:pininitialdatavariationrandom}
\end{figure}

\begin{figure}[htbp]
	\centering
	\subfloat[$\epsilon=0.5$] {\includegraphics[width=0.24\textwidth]{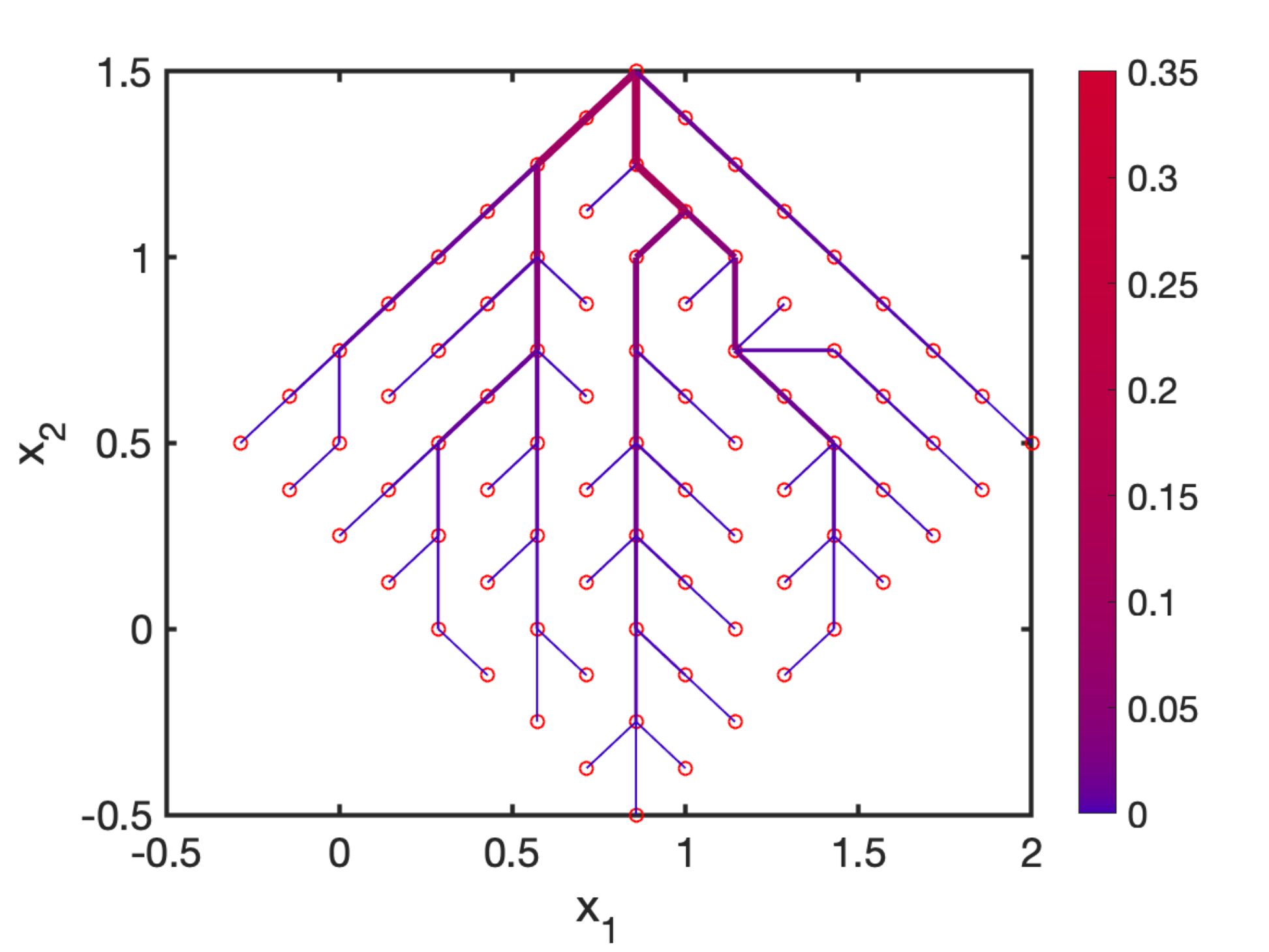}}
	\subfloat[$\epsilon=1$] {\includegraphics[width=0.24\textwidth]{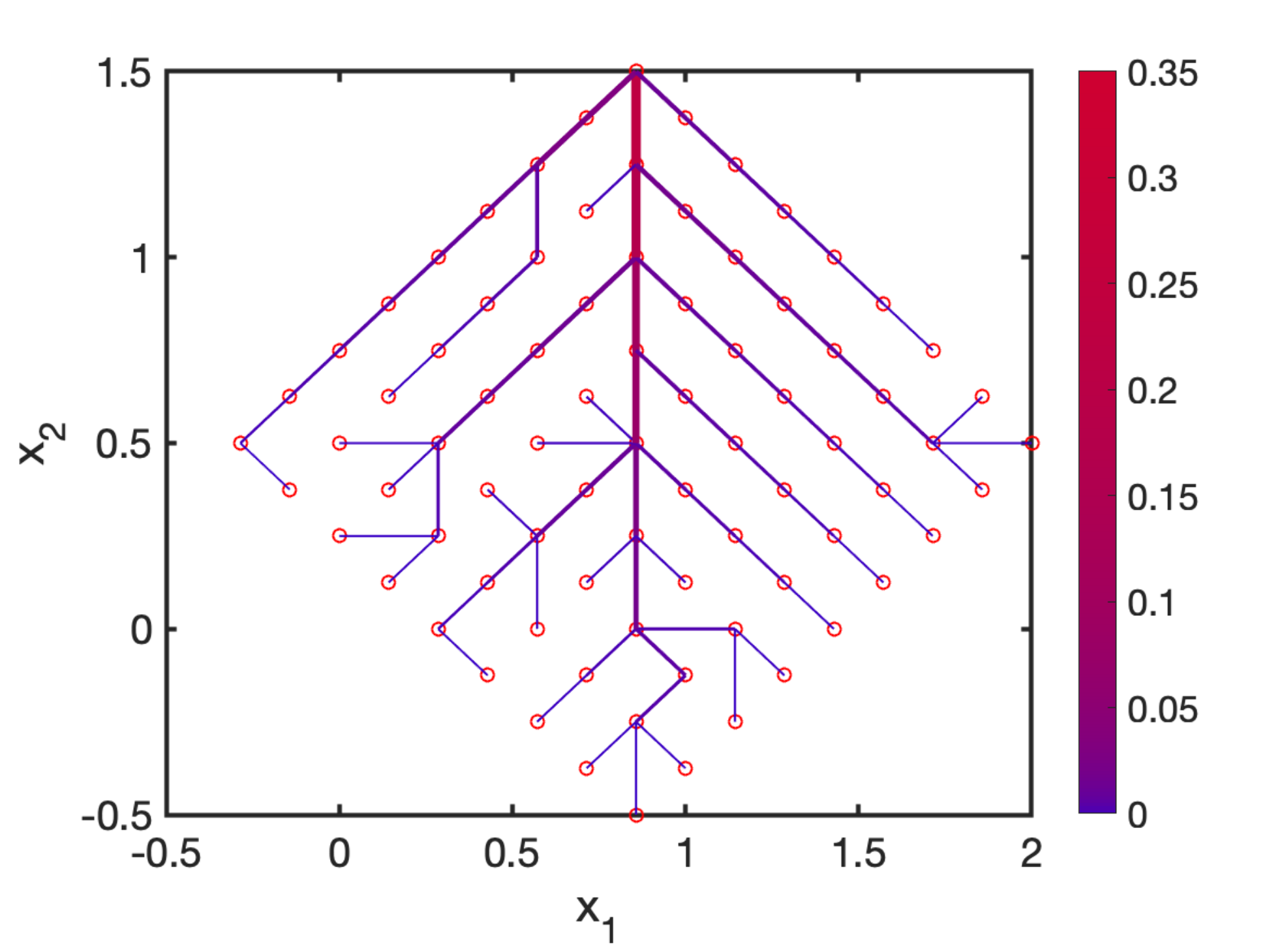}}
	\subfloat[$\epsilon=5$] {\includegraphics[width=0.24\textwidth]{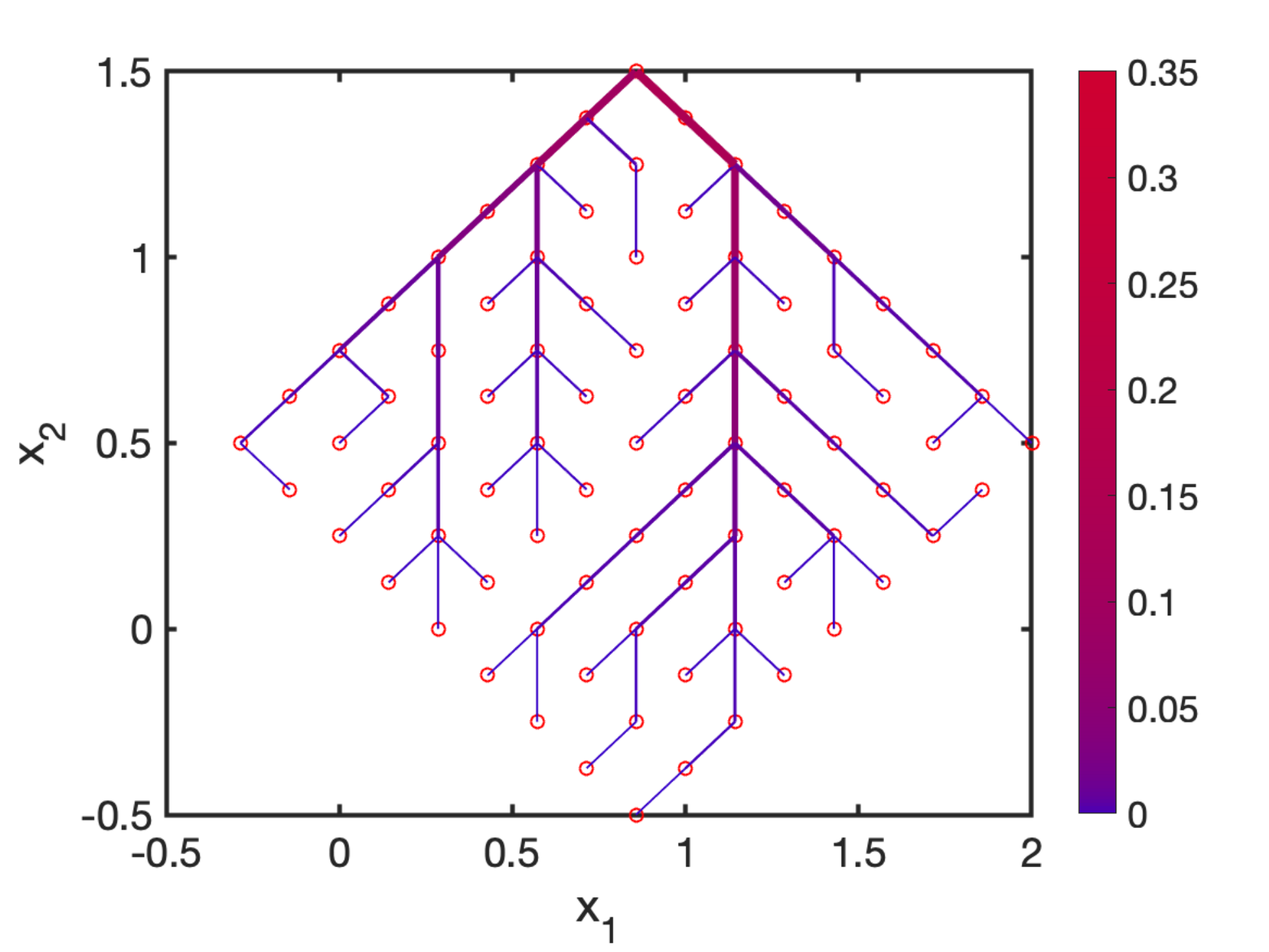}}
	\subfloat[$\epsilon=100$] {\includegraphics[width=0.24\textwidth]{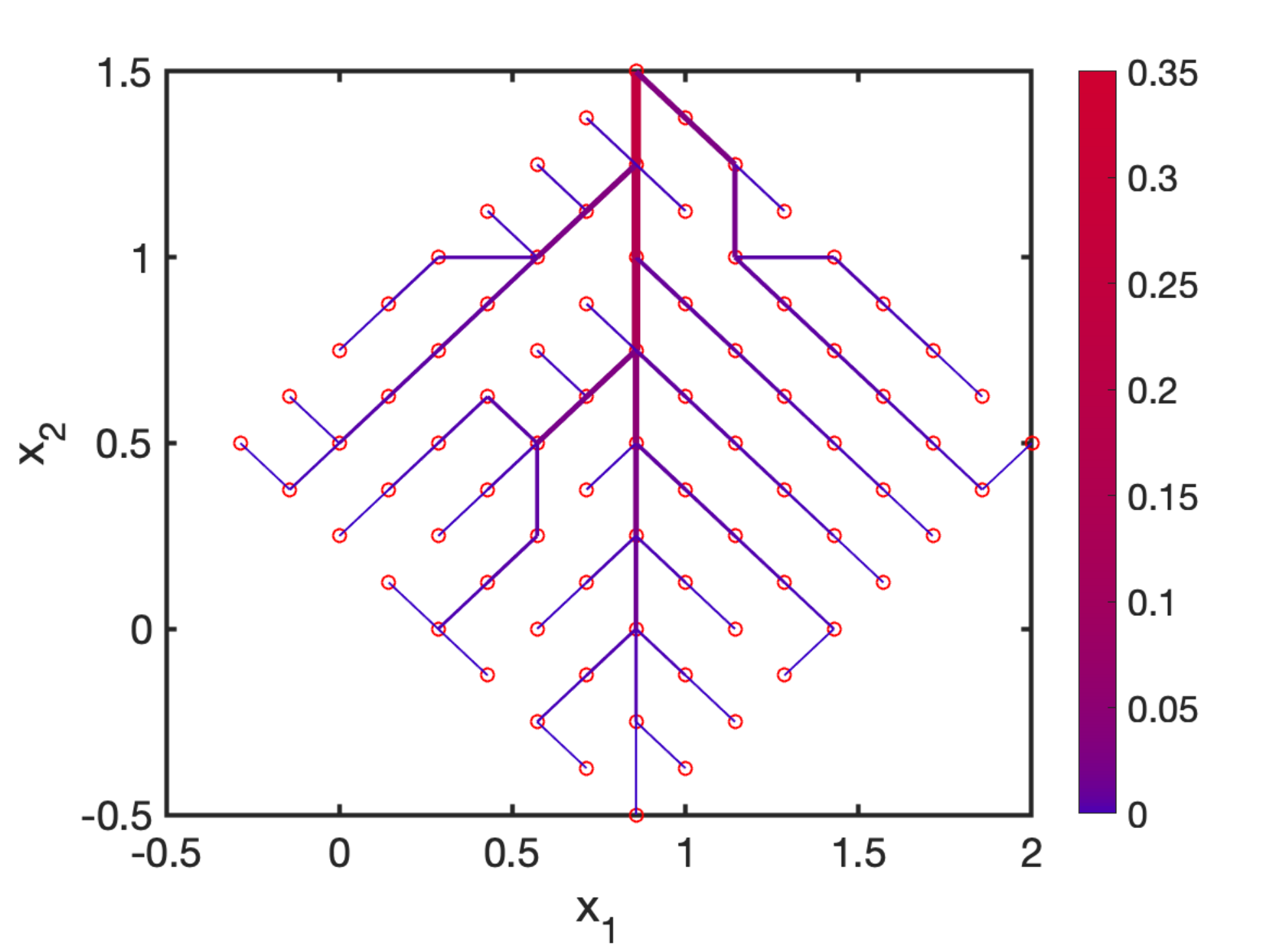}\label{fig:pininitialdatavariation100}}
	\caption{Steady states for transport activity for initial transport activity $\epsilon\mathcal{U}(0,1)$.}\label{fig:pininitialdatavariation}
\end{figure}

\begin{figure}[htbp]
	\centering
	\subfloat[$\epsilon=1$] {\includegraphics[width=0.24\textwidth]{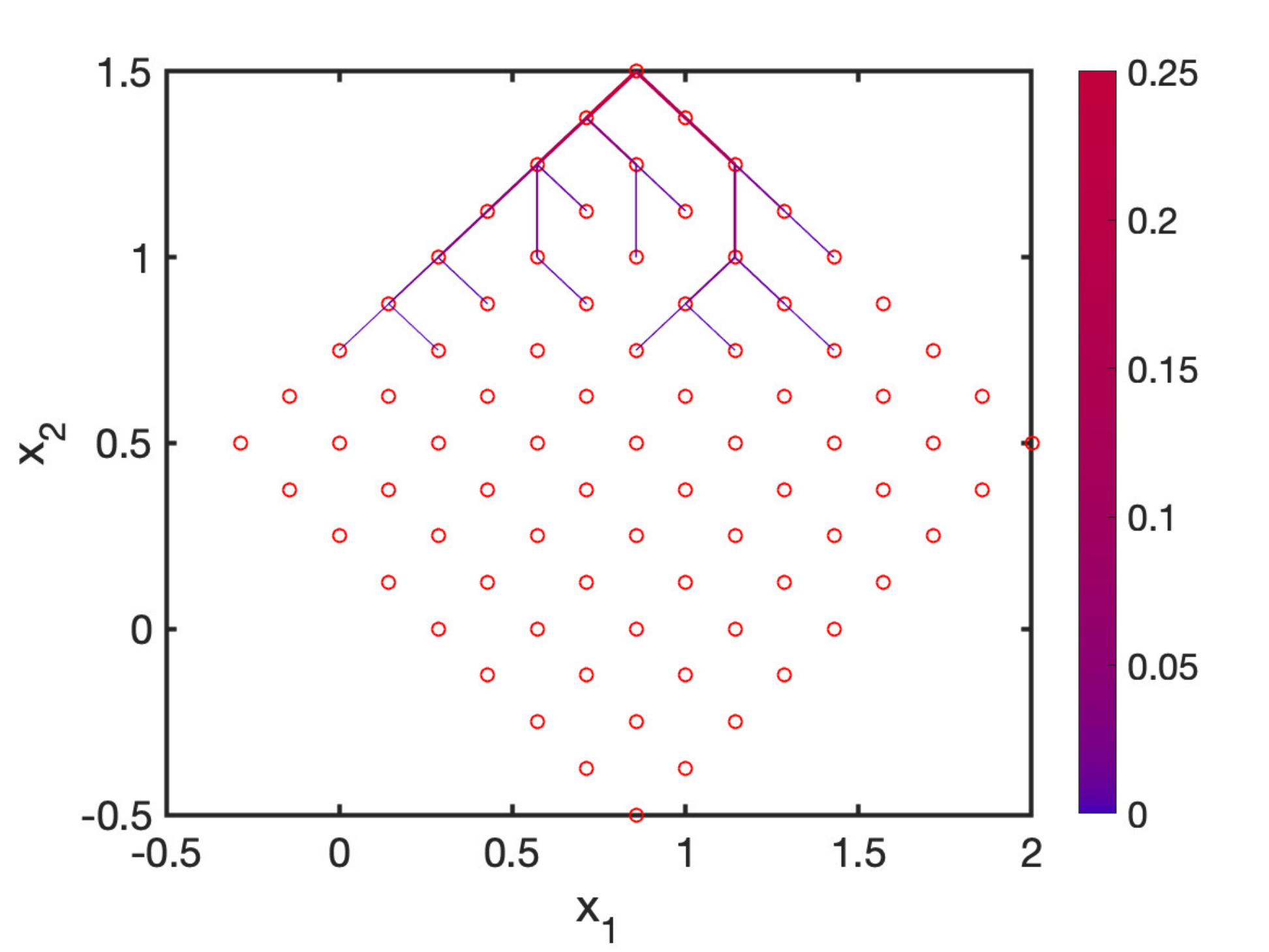}}
	\subfloat[$\epsilon=5$] {\includegraphics[width=0.24\textwidth]{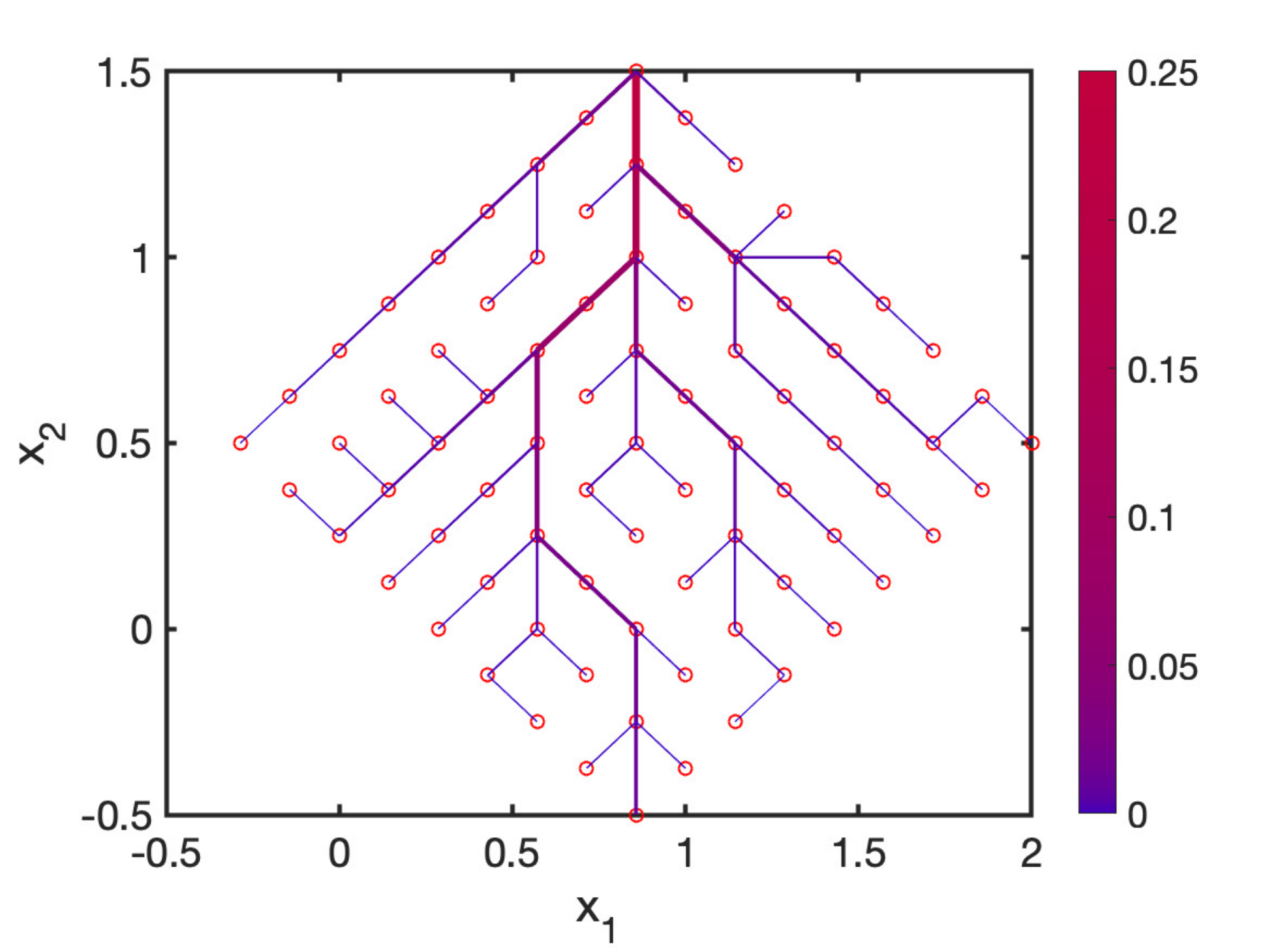}}
	\subfloat[$\epsilon=50$] {\includegraphics[width=0.24\textwidth]{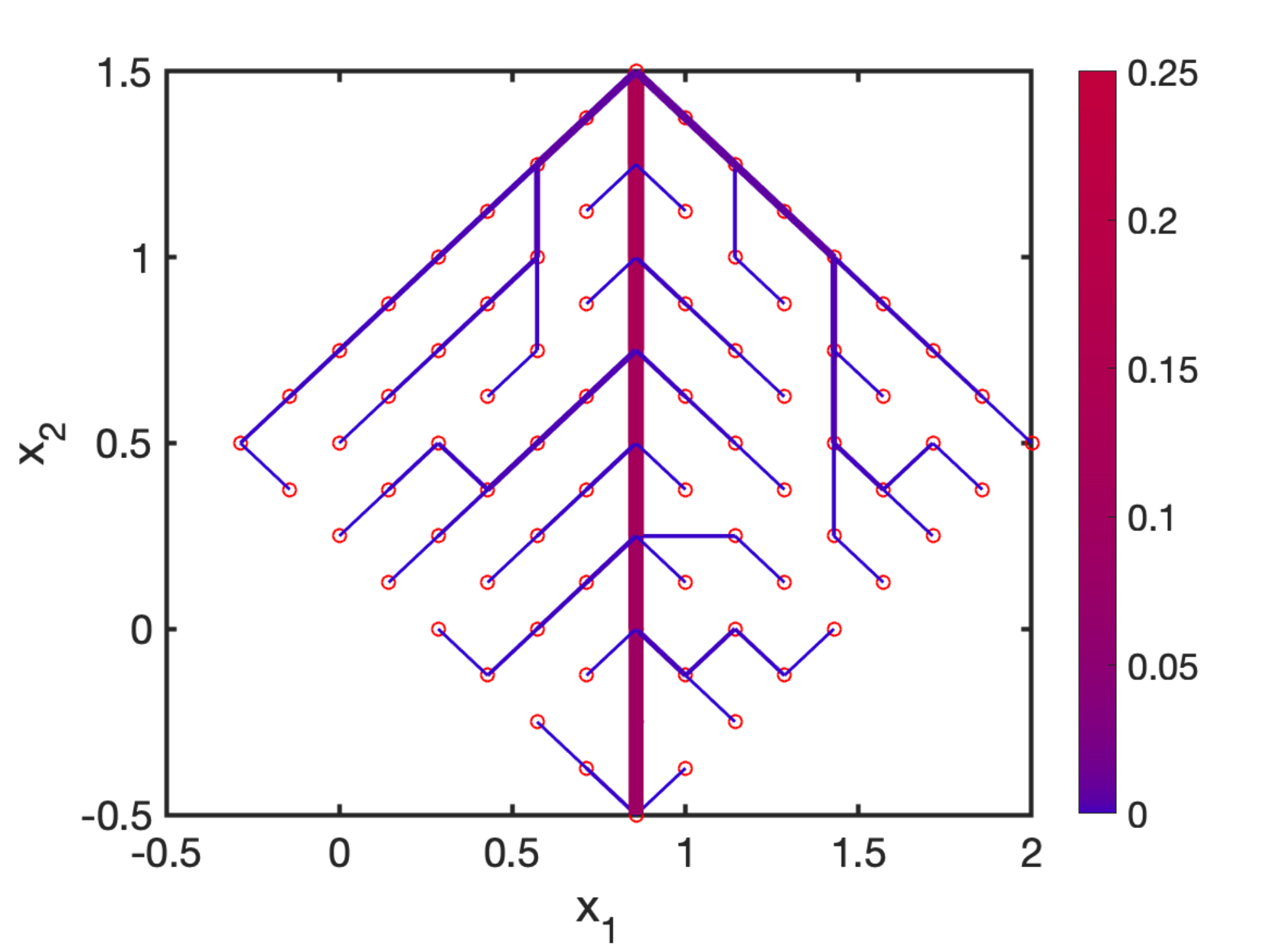}}
	\subfloat[$\epsilon=100$] {\includegraphics[width=0.24\textwidth]{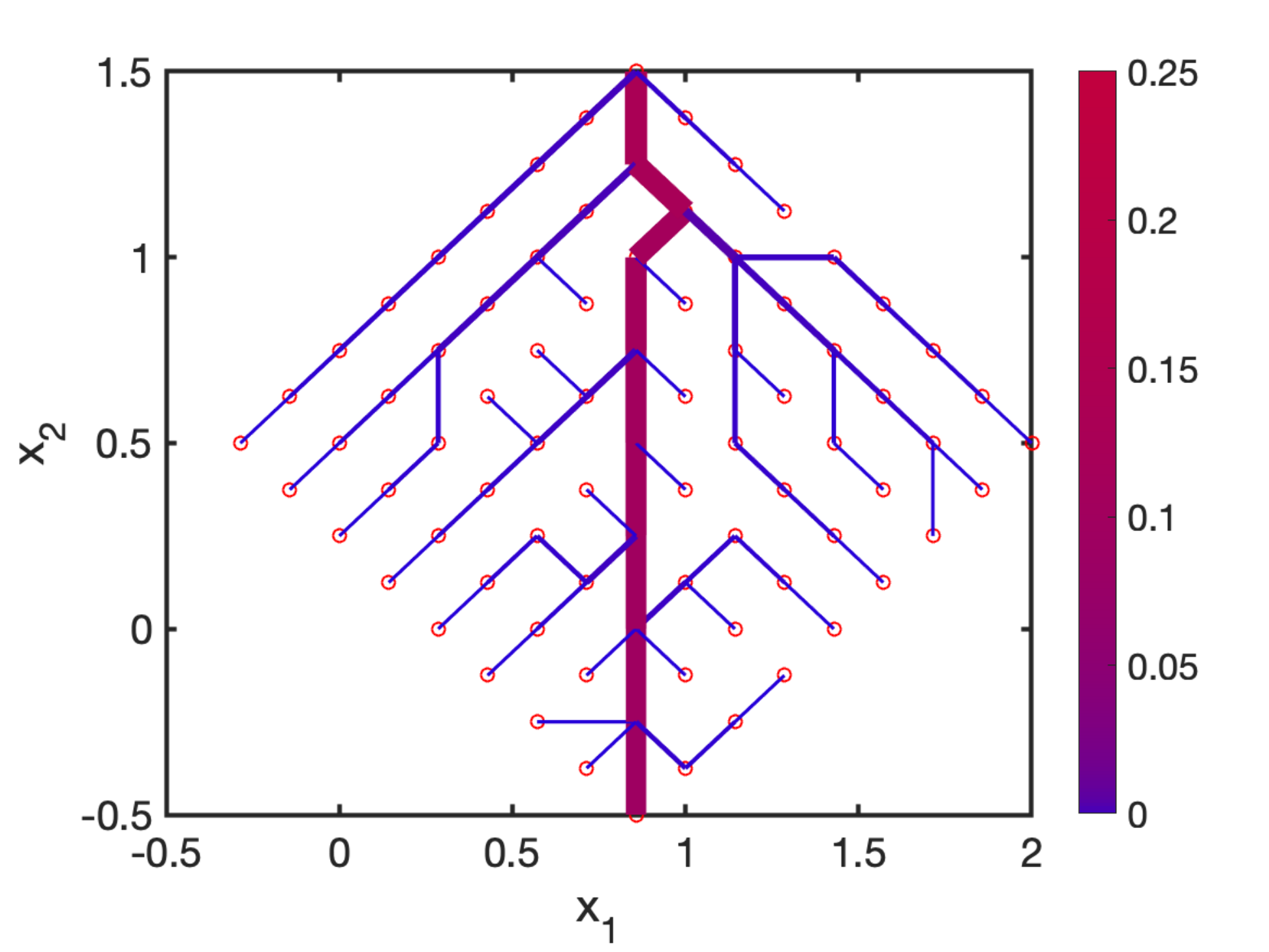}}
	\caption{Steady states for transport activity for initial transport activity $100\mathcal{U}(0,1)$ with source $10\epsilon$ and sink $\epsilon$.}\label{fig:pininitialdatavariationsources}
\end{figure}

In Figures \ref{fig:pininitialdatavariationrandom}--\ref{fig:pininitialdatavariationsources}, we vary the initial  auxin transport activity  and no longer consider the initial data $\bar{X}_{ij}$. In Figure \ref{fig:pininitialdatavariationrandom}, the steady states for the  transport activity  are shown where the initial  transport activity is chosen as $\theta+ 0.00001\epsilon$ for parameter $\epsilon \in \{0.5,5,50,100\}$ and a random variable $\theta$   with $\theta=1$ with probability $0.2$ and $\theta=0$ with probability $0.8$. In particular, the resulting patterns of the  transport activity  have no symmetries and the location of the mid-veins strongly depend on the choice of parameters, illustrating that model \eqref{eq:auxineqproduction}--\eqref{eq:pineq} can produce complex vein patterns. Note that the size of the stationary pattern increases as $\epsilon$ and, thus, as the absolute value of the initial transport activity increases.

In Figure \ref{fig:pininitialdatavariation}, we consider the initial transport activity $\epsilon \mathcal{U}(0,1)$ for $\epsilon \in\{0.5,1,5,100\}$. These numerical results demonstrate that model  \eqref{eq:auxineqproduction}--\eqref{eq:pineq} is capable to produce different complex stationary state, not only on subdomains as in Figure \ref{fig:pininitialdatavariationrandom}, but  on the entire underlying network. In particular, the stationary transport activity connects auxin sources and sinks.

In Figure \ref{fig:pininitialdatavariationsources}, we consider the same initial condition for the transport activity as in Figure \subref*{fig:pininitialdatavariation100}, i.e.\ $100 \mathcal{U}(0,1)$, but we vary the strengths $10\epsilon$ and $\epsilon$   of the auxin background source strengths $\xi_S$ and  sink strengths $\xi_I$, respectively, where  $\epsilon \in \{1,5,50,100\}$. One can clearly see in Figure~\ref{fig:pininitialdatavariationsources} that  the auxin sources and sinks are not strong enough for $\epsilon=1$ for transport activity to connect the top and bottom corners of the underlying network, while for larger values of $\epsilon$  mid-veins become  visible and get stronger as auxin sources and sinks increase. This shows that  complex stationary transport activity patterns with no symmetries and  major mid-veins can be obtained.

\begin{figure}[htbp]
	\centering
	\subfloat[$3$  sources] {\includegraphics[width=0.24\textwidth]{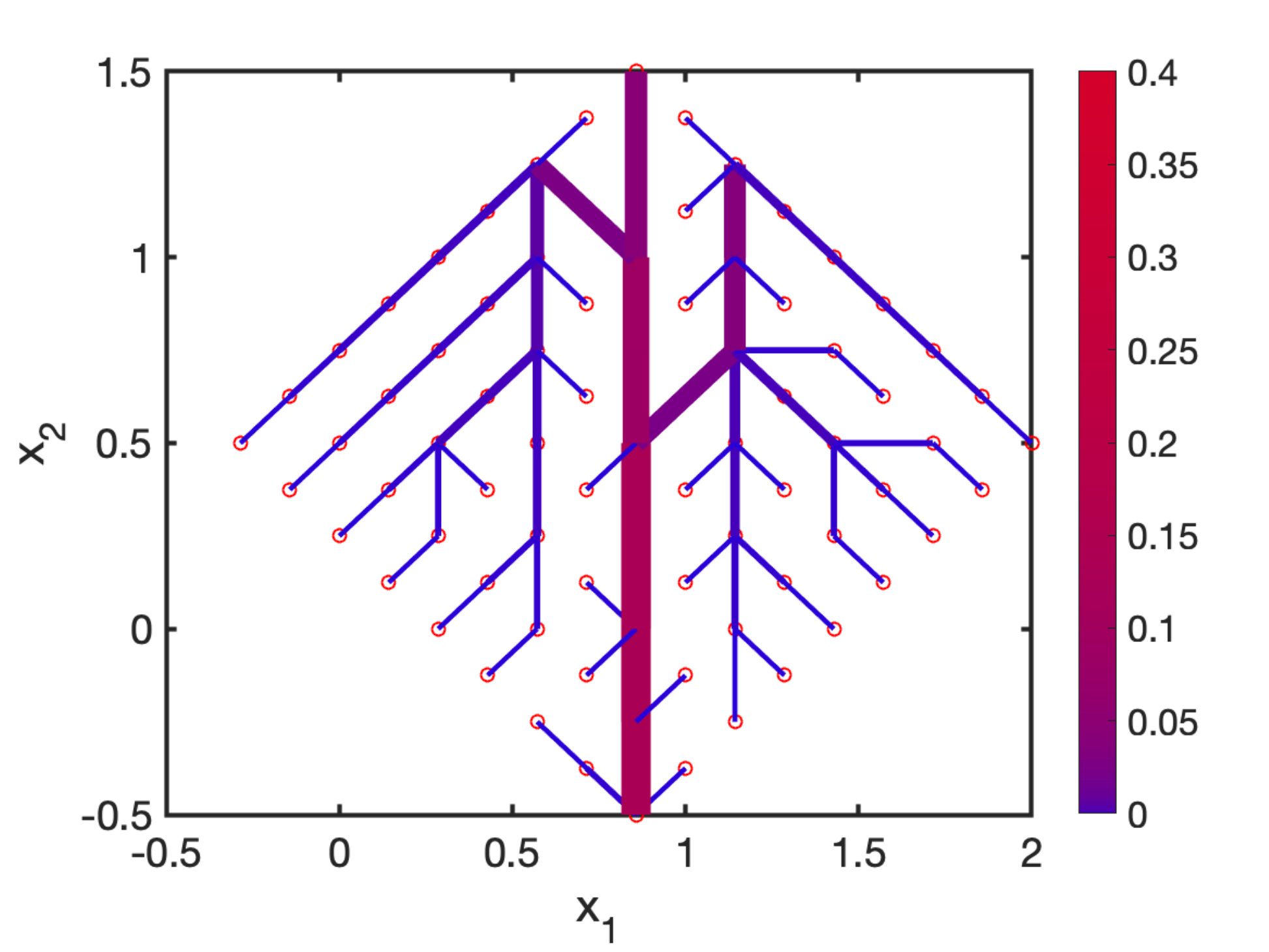}\label{fig:pininitialdatavariationnumbersources1}}
	\subfloat[$5$  soures] {\includegraphics[width=0.24\textwidth]{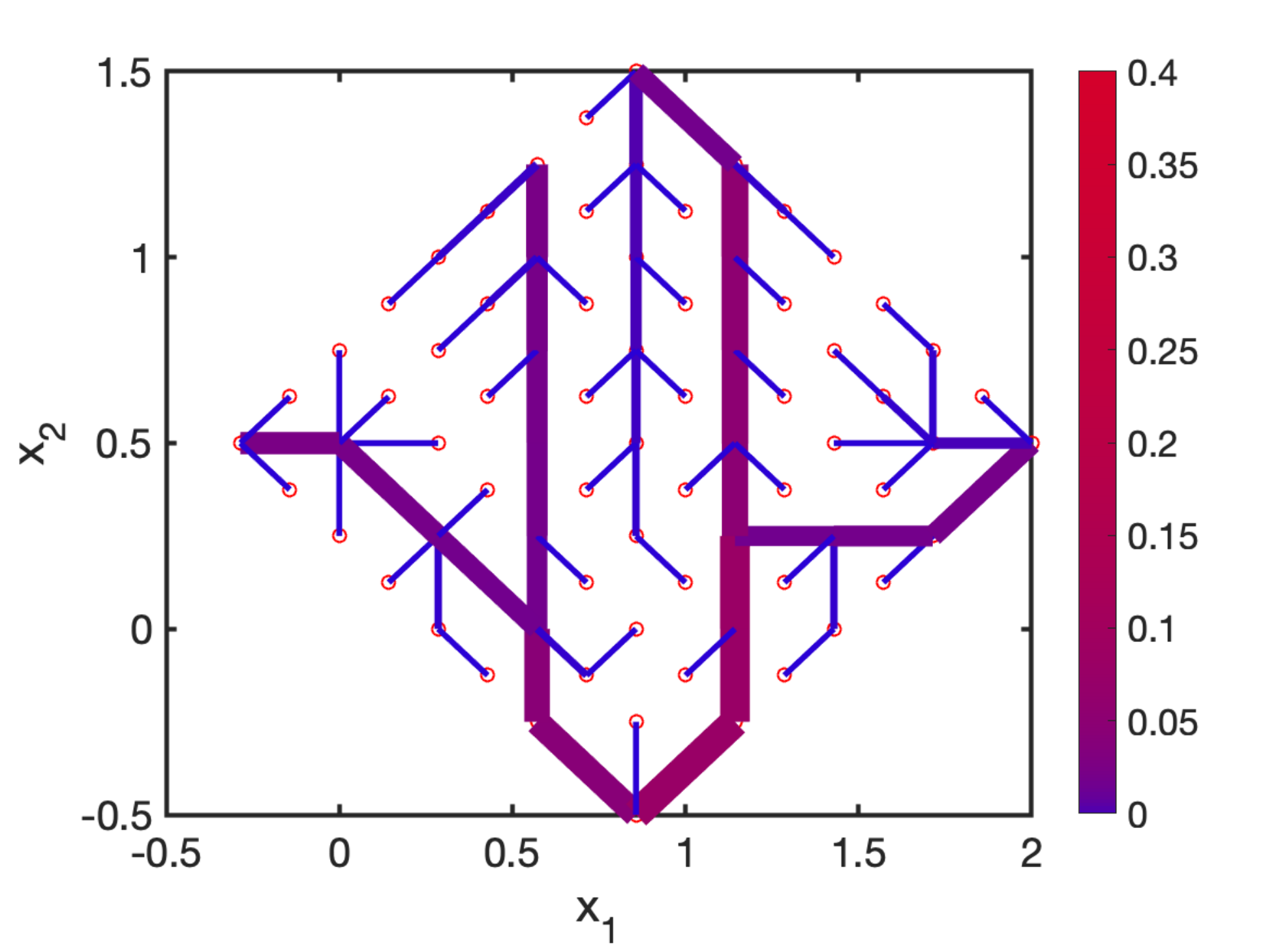}\label{fig:pininitialdatavariationnumbersources2}}
	\subfloat[$7$  sources] {\includegraphics[width=0.24\textwidth]{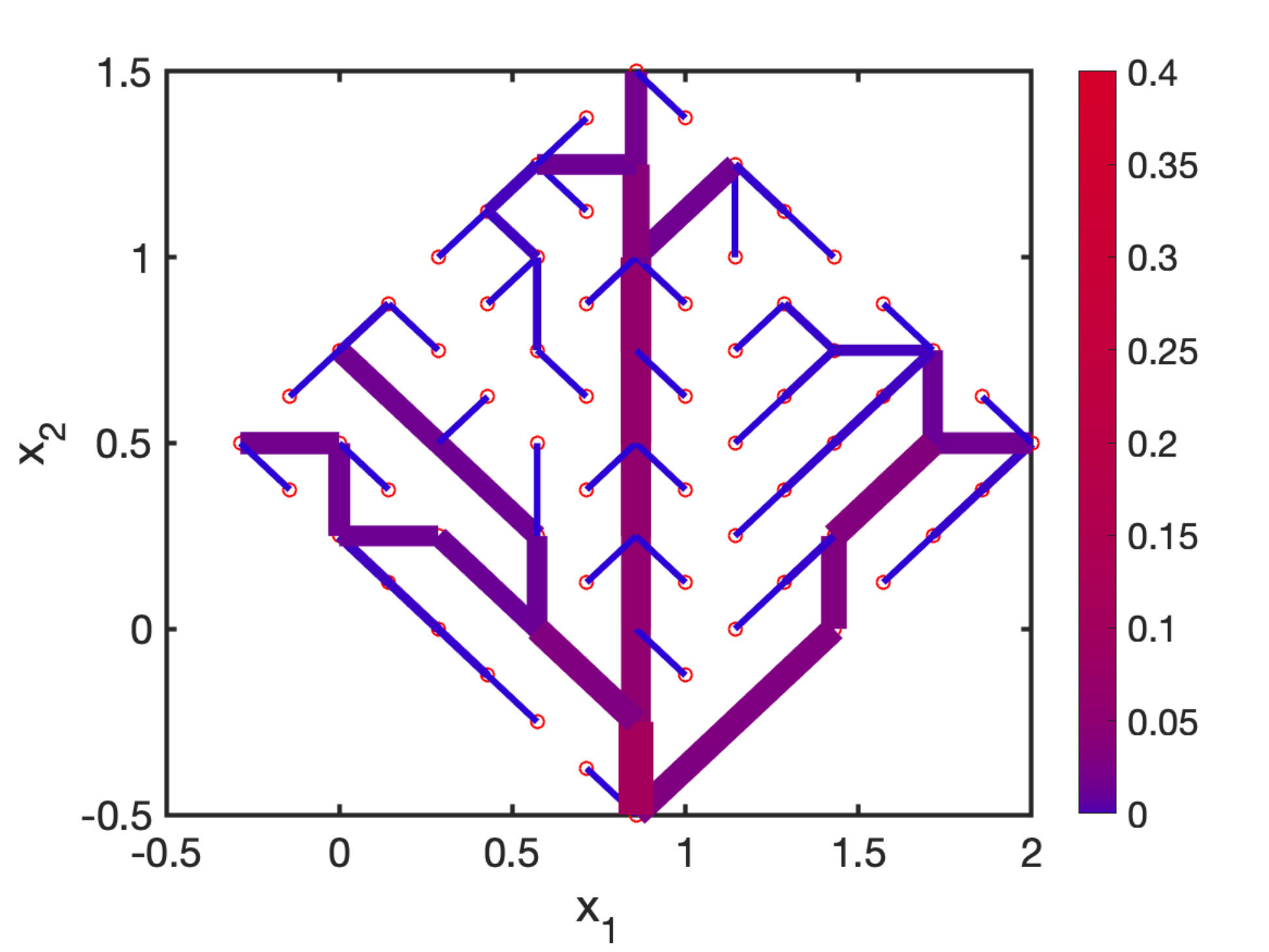}\label{fig:pininitialdatavariationnumbersources3}}
	\subfloat[$9$  sources] {\includegraphics[width=0.24\textwidth]{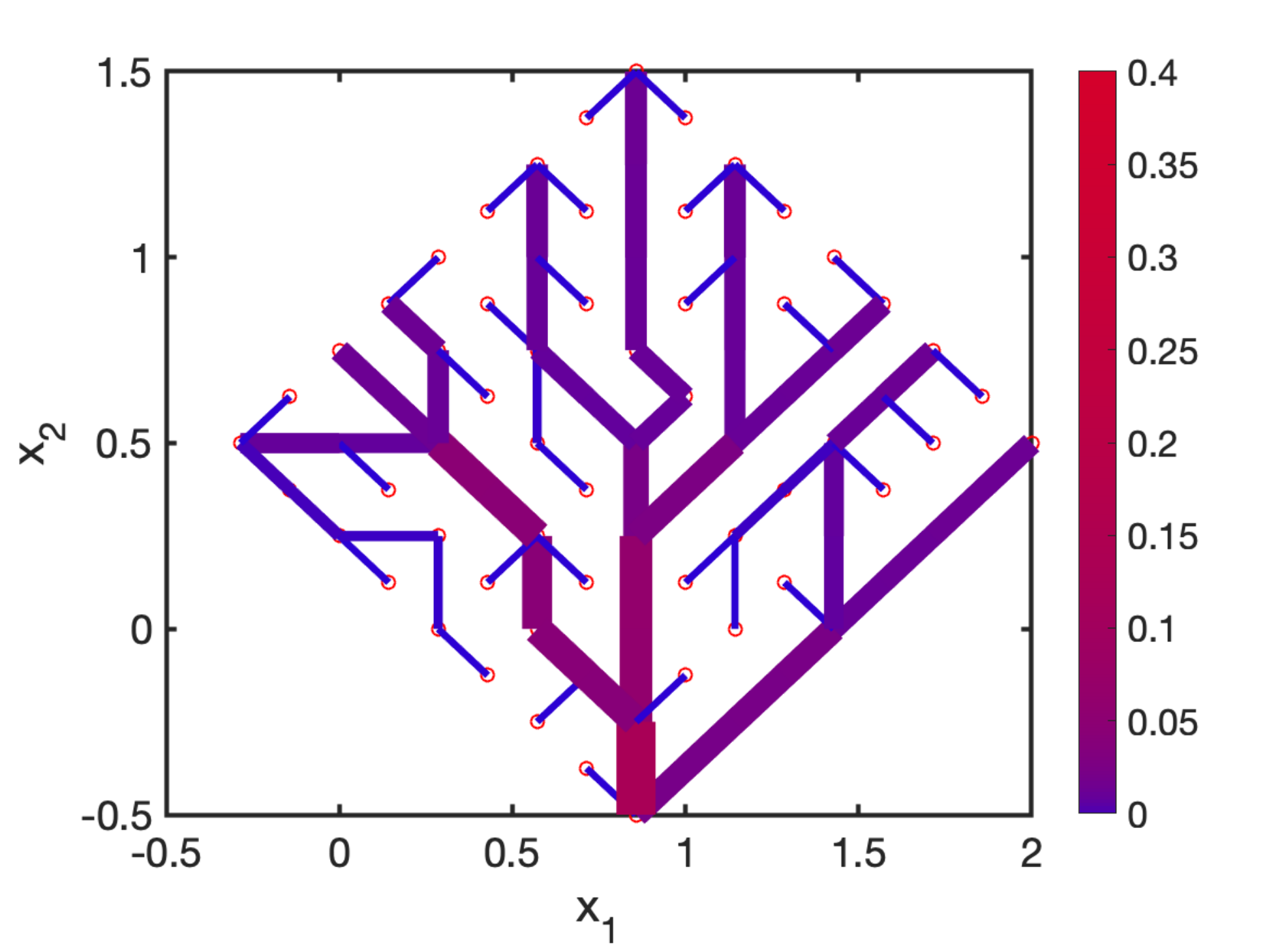}\label{fig:pininitialdatavariationnumbersources4}}
	\caption{Steady states for transport activity for initial transport activity $100\mathcal{U}(0,1)$ with different number of sources of strength $1000$ and sinks of strength $100$.}\label{fig:pininitialdatavariationnumbersources}
\end{figure}

\begin{figure}[htbp]
	\centering
	\subfloat[$3$ sources, $1$ sink] {\includegraphics[width=0.24\textwidth]{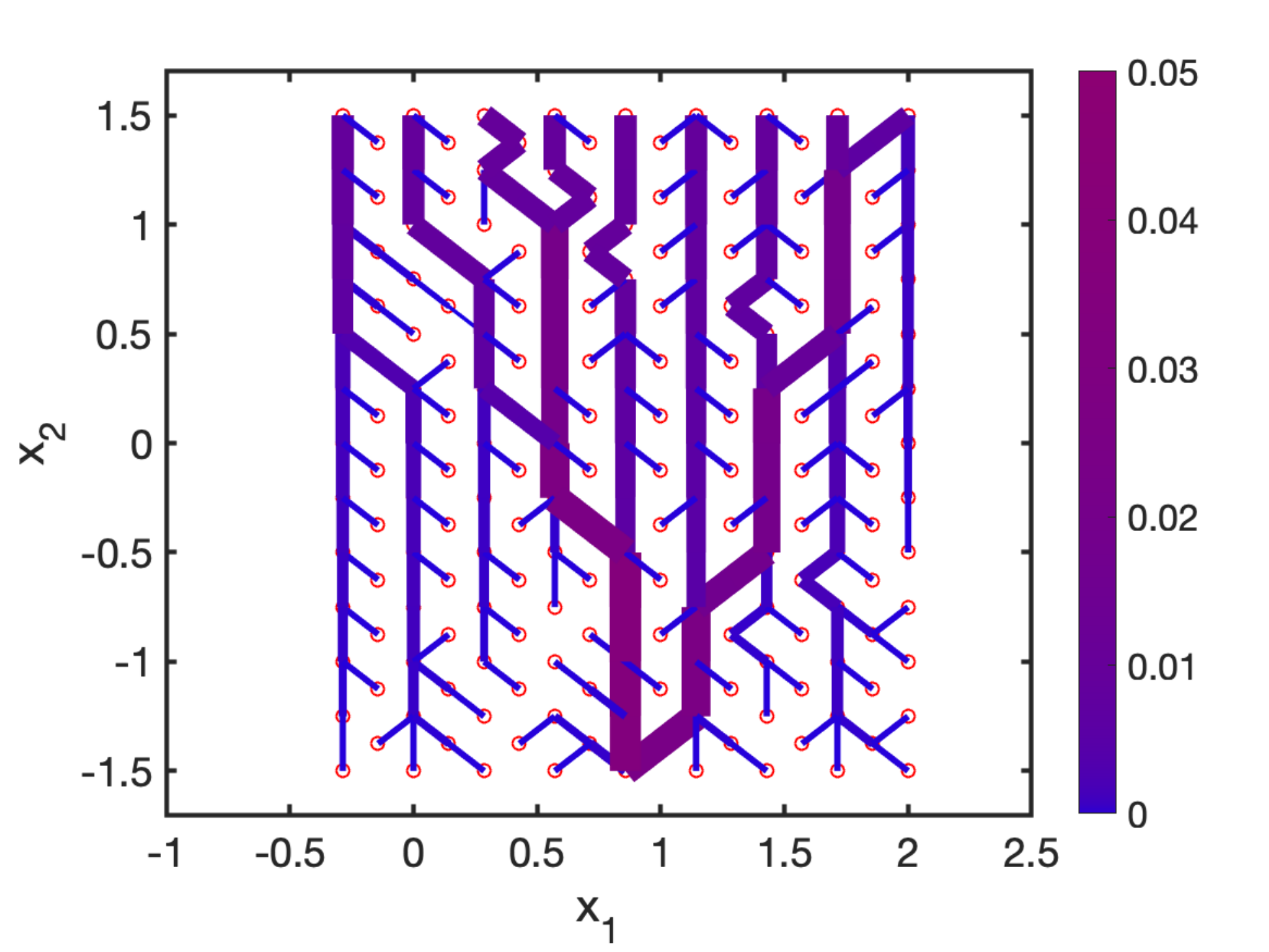}\label{fig:pininitialdatavariationnumbersourcesrectangle1}}
	\subfloat[$5$  soures, $1$ sink] {\includegraphics[width=0.24\textwidth]{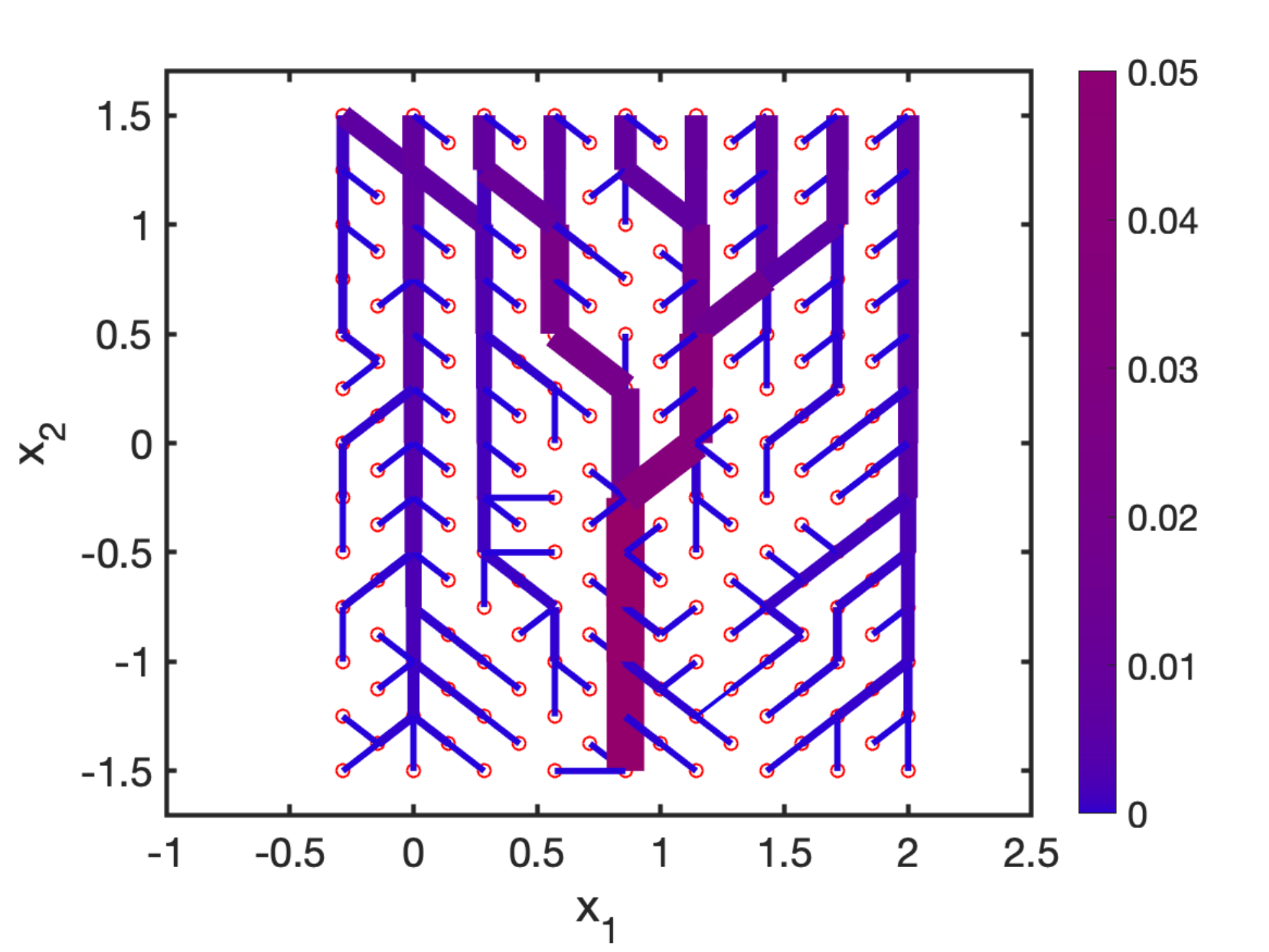}\label{fig:pininitialdatavariationnumbersourcesrectangle2}}
	\subfloat[$5$  sources, $3$ sinks] {\includegraphics[width=0.24\textwidth]{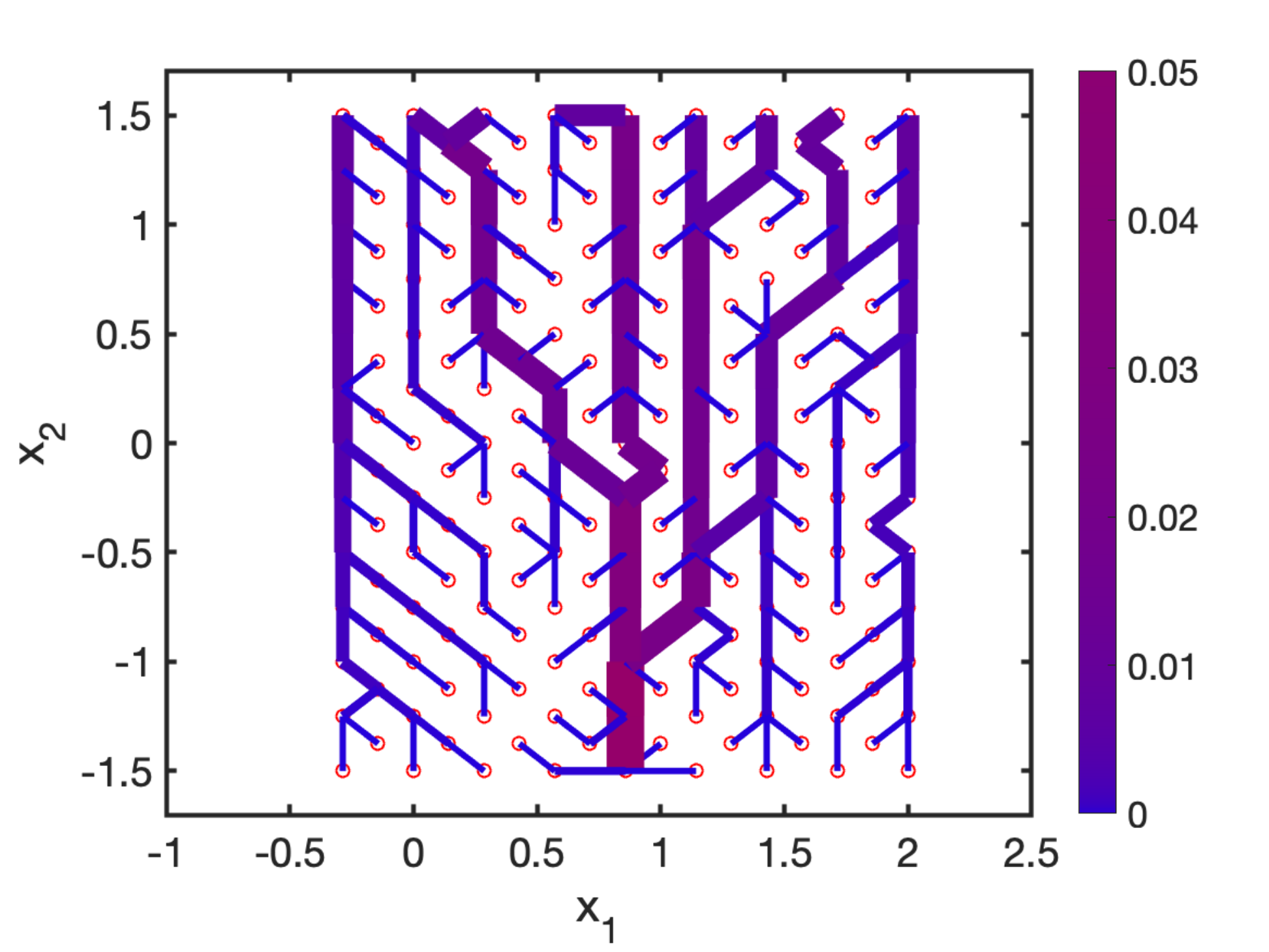}\label{fig:pininitialdatavariationnumbersourcesrectangle3}}
	\subfloat[$5$  sources, $5$ sinks] {\includegraphics[width=0.24\textwidth]{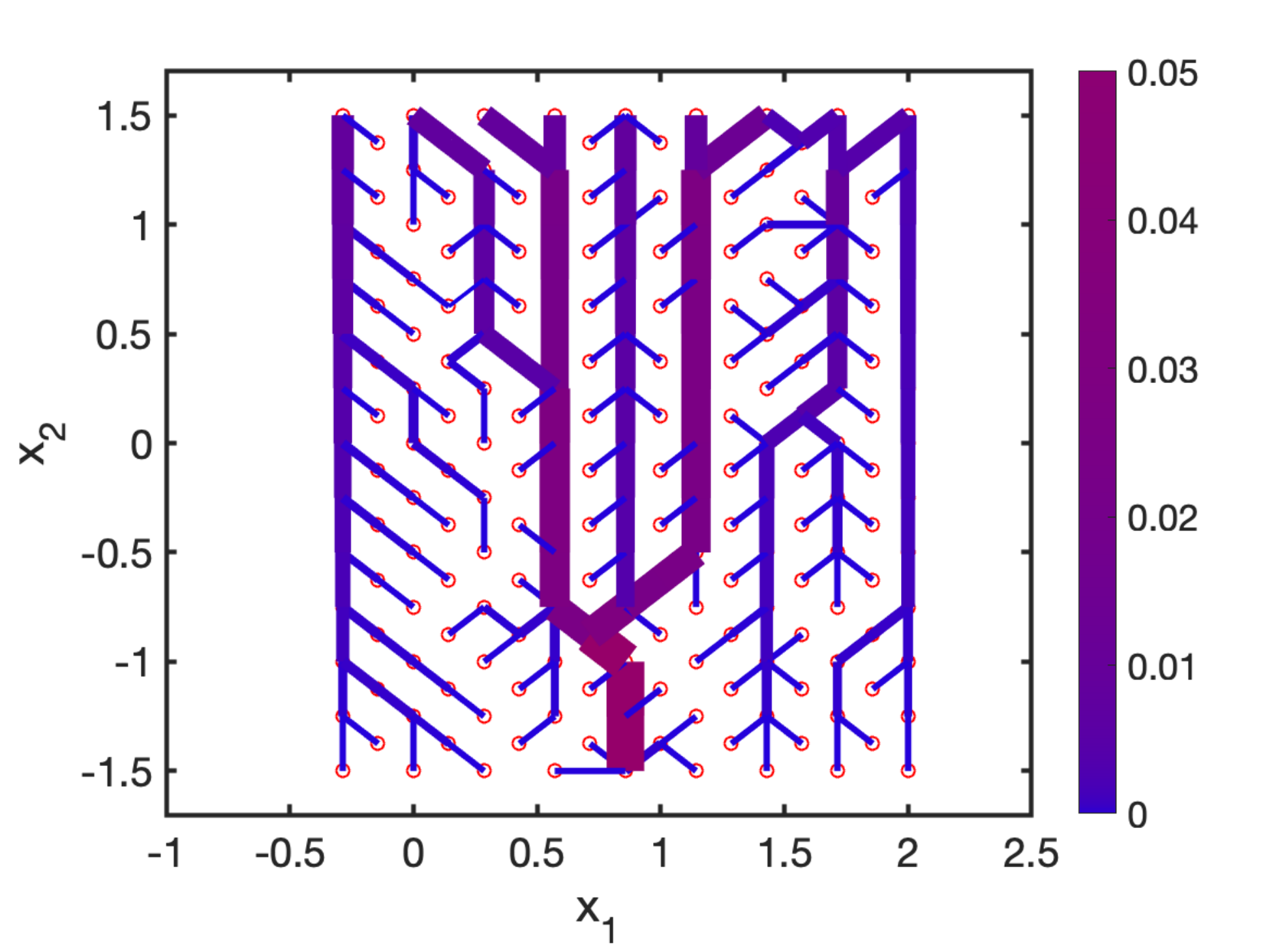}\label{fig:pininitialdatavariationnumbersourcesrectangle4}}
	\caption{Steady states for transport activity for initial transport activity $100\mathcal{U}(0,1)$ with different number of sources of strength $1000$ and different number of sinks of strength  $100$.}\label{fig:pininitialdatavariationnumbersourcesrectangle}
\end{figure}

In Figures \ref{fig:pininitialdatavariationnumbersources} and  \ref{fig:pininitialdatavariationnumbersourcesrectangle} we consider multiple sources and sinks for obtaining more realistic vein networks. Starting from a certain configuration of sources and sinks in Figures \subref*{fig:pininitialdatavariationnumbersources1} and \subref*{fig:pininitialdatavariationnumbersourcesrectangle1} we subsequently add sources and sinks in the subfigures further to the right. In Figure \ref{fig:pininitialdatavariationnumbersources} we consider a diamond grid as in most figures, but apart from a source at the top corner and a sink at the bottom corner of the grid, we add sources which are located symmetrically with respect to the longest vertical axis of the grid. Denoting the distance between the left and the top corner of grid by $l$, these sources are located on the boundary of the grid at a distance of $l/4$ from the top corner (Figures \subref*{fig:pininitialdatavariationnumbersources1}, \subref*{fig:pininitialdatavariationnumbersources2}, \subref*{fig:pininitialdatavariationnumbersources3}, \subref*{fig:pininitialdatavariationnumbersources4}),  the left corner (Figures \subref*{fig:pininitialdatavariationnumbersources2}, \subref*{fig:pininitialdatavariationnumbersources3}, \subref*{fig:pininitialdatavariationnumbersources4}) and at  distances of $3l/4$ and $5l/8$ from the top corner in Figures \subref*{fig:pininitialdatavariationnumbersources3}, \subref*{fig:pininitialdatavariationnumbersources4} and Figure \subref*{fig:pininitialdatavariationnumbersources4}, respectively. Similarly, the sources are located on the right side of the grid by symmetry of the source locations in each figure. One can clearly see that multiple sources result in a more complex transportation network between the sources and the sink in comparison to the simulation results in the previous figures with merely one point source. 

In Figure \ref{fig:pininitialdatavariationnumbersourcesrectangle} we consider a rectangular underlying grid with sources at the top and the bottom of the boundary of the grid. We denote the length between the left top and right top corner of the grid by $l$. We consider a sink in the middle of the bottom boundary and sources in the middle of the top boundary and at a distance of $l/4$ left and right of the middle on the top boundary in all subfigures of Figure  \ref{fig:pininitialdatavariationnumbersourcesrectangle}. Additional sources are located at the left top and the right top corner in Figures  \subref*{fig:pininitialdatavariationnumbersourcesrectangle2},  \subref*{fig:pininitialdatavariationnumbersourcesrectangle3}, \subref*{fig:pininitialdatavariationnumbersourcesrectangle4}. In Figures \subref*{fig:pininitialdatavariationnumbersourcesrectangle3}, \subref*{fig:pininitialdatavariationnumbersourcesrectangle4} additional sinks are added at the bottom boundary in a distance of $l/4$ left and right of the middle of the bottom boundary, while in Figure  \subref*{fig:pininitialdatavariationnumbersourcesrectangle4} additional sinks are considered in the left bottom and right bottom corner of the grid. In particular, the resulting patterns look very similar to those in leaves.

Model \eqref{eq:auxineqproduction}--\eqref{eq:pineq} describes the auxin transport with a positive feedback between auxin fluxes and auxin transporters where the auxin transporters are not necessarily polar.
The above numerical results illustrate that the model \eqref{eq:auxineqproduction}--\eqref{eq:pineq} is able to connect an auxin source-sink pair with a mid-vein and that branching vein patterns can also be produced. 	A nice feature of the model is that the veins end up with high auxin levels. This was not achieved with the original Mitchinson models  and this has been discussed in some detail. A solution to this has been to adapt the conservative approach $X_{tot} = \sum_{j \in \mathcal{N}(i)} X_{ij} = \text{const}$ for the auxin transporters which (together with feedback on the localisation of auxin transporters from auxin flux) can lead to high auxin in veins.

%
%The ability of the original flux-feedback model showed that if $\mathcal{P}_{ij} \propto |flux|^{\alpha}$ leads to vein formation if $\alpha \geq 2$ and 'homogeneous' flow between sources and sinks if $\alpha=1 (<2)$. I think it makes sense to discuss the $\kappa \gamma$ parameters in this context.

We want to stress here that our model \eqref{eq:auxineqproduction}--\eqref{eq:pineq} is able to generate a venation/transport network without a polar input, as seen in the case when auxin transporters are knocked out in the various numerical examples.

In reality, the venation patterns appear while the leaf is growing, and as such our simulations (and \review{the simulation results of} many previous PIN-based flux models  on static geometries) can only provide part of the answer. Changing the configuration of sources and sinks in the model is expected to lead to \review{different patterns} in the final leaf.

%%%%%%%%%%%%%%%%%%%%%%%%%%%%%%%%%%%%%%%%%%%%%%%%%%%%%%%%%%%
\def\P{X}
\section{The formal continuum limit}\label{sec:continuum}

The main reason for focusing on discrete models is that the patterns form when the leaves have very few cells, e.g.\ the (first) mid-vein forms when the leaf is about five cells wide. 
Cells split over time, resulting in a larger number of cells and network growth. Besides, there is an auxin peak at the tip before the high auxin/transport activity vein forms downwards from this. Still, this does not discard alternative mechanisms setting up an intitial pattern that connects the leaf tip with the vasculature in the stem (thought to be auxin sink). These phenomena can be modeled much better in a diffusion driven setting instead of the discrete setting and motivates us to consider the associated macroscopic model.

The goal of this section is to derive the formal macroscopic limit of the discrete model \eqref{eq:pineq}, \eqref{eq:auxineqproduction}
as the number of nodes and edges tends to infinity,
and to study the existence of weak solutions of the resulting PDE system. 
The derivation requires an appropriate rescaling of the auxin production equation \eqref{eq:auxineqproduction}.
Moreover, since the derivation of macroscopic limits of systems posed on general (unstructured) graphs is a highly nontrivial topic,
see, e.g., \cite{graphs}, we restrict ourselves to discrete graphs represented by regular equidistant grids, i.e.,
tessellations of a rectangular domain $\Omega\subset\R^d$, $d\in\N$, by congruent identical rectangles (in 2D) or cubes (in 3D)
with edges parallel to the axis. 
The results can be generalized to parallelotopes,
%A generalization of the result for parallelotopes will be given in Remark \ref{rem:parallelotopes}.
see \cite[Section 3]{HaKrMa} for details of the formal procedure applied to the Hu-Cai model \eqref{eq:GF-orig}--\eqref{eq:kirchhoff},
and \cite{HaKrMa2} for the rigorous procedure in the spatially one- and two-dimensional setting.

\subsection{The formal derivation of the continuum limit of the system \eqref{eq:pineq}, \eqref{eq:auxineqproduction}}\label{sec:derivation}
Given the graph $G=(\Vset,\Eset)$ as a rectangular tesselation of the rectangular domain $\Omega$,
let us denote the vertices left and right of vertex $i\in\Vset$ along the $k$-th spatial dimension by $(i-1)_k$ and, resp., $(i+1)_k$.
Moreover, let us denote $h_k>0$ the equidistant grid spacing in the $k$-th dimension.
The rescaled auxin production equation \eqref{eq:auxineqproduction} is then written as
\begin{align} \label{aux-rescaled}
\tot{a_i}{t} = S_i - I_ia_i + \delta\sum_{k=1}^d \frac{1}{h_k} \left(
X_{i,(i+1)_k}\frac{a_{(i+1)_k} - a_i}{h_k} - X_{i,(i-1)_k}\frac{a_i - a_{(i-1)_k}}{h_k}
\right)
\qquad \text{for~}i\in \Vset.
\end{align}
The rescaling of the sum on the right hand side by $h_k$ is reflecting the fact that the edges of the graph are inherently one-dimensional structures,
embedded into the $d$-dimensional space, cf. \cite[Section 3]{HaKrMa}.
A straightforward calculation reveals that \eqref{aux-rescaled} is a finite difference discretization of the parabolic equation
\begin{align}
\part{a}{t} = \delta \grad\cdot( \P\grad a) + S - I a, \label{cont1}
\end{align}
on the regular grid $G=(\Vset,\Eset)$, where $a=a(t,x)$ is a formal limit of the sequence of discrete auxin concentrations $(a_i)_{i\in\Vset}$
as $|\Vset|\to\infty$, and $I=I(x)$ is a formal limit of the sequence $(I_i)_{i\in\Vset}$.
Here, $\P=\P(t,x)$ is the diagonal tensor 
$\P =\text{diag} (\P_1 ,\ldots, \P_d)$
%\begin{align} \label{PTensor}
%\P = \begin{pmatrix} \P_1 & & \\ & \ddots & \\ & & \P_d \end{pmatrix},
%\end{align}
where $\P_k$ is the formal limit of the sequence $(\P_{ij})_{i,j\in\Vset}$
on edges $(i,j)\in\Eset$ oriented along the $k$-th spatial direction.
A formal continuum limit of \eqref{eq:pineq} yields the family of ODEs for $\P=\P(t,x)$,
\begin{align}
\part{\P_k}{t} = \left( \frac{|q_k|^\kappa}{\P_k^{\gamma+1}} - \tau\right) \P_k, \label{cont2}
\end{align}
with $q_k = \P_k\partial_{x_k} a$.
Note that the product $\P\grad a$ is the vector $\P\grad a = (\P_1\partial_{x_1} a, \dots, \P_d\partial_{x_d} a)$.

Observe that \eqref{cont2} is in fact a family of ODEs for $\P_k=\P_k(t,x)$, parametrized by $x\in\Omega$.
Consequently, in analogy to \cite{HaKrMa}, we introduce the diffusive terms $D^2 \laplace \P_k$ that model random fluctuations in the medium.
Thus, the updated version of \eqref{cont2} reads
\begin{align}
\part{\P_k}{t} = D^2\laplace \P_k + \left( \frac{|q_k|^\kappa}{\P_k^{\gamma+1}} - \tau\right) \P_k, \label{cont3}
\end{align}
with the diffusion coefficient $D^2>0$.

Biological observations suggest that the auxin dynamics takes place on a faster time scale
than the dynamics of the transporter proteins in the order of minutes for auxin movement \cite{delbarre:1996}, and in the order of hours for e.g. PIN1 reorientation \cite{Heisler2010}. 
Consequently, we consider a formal fast time scale limit of \eqref{cont1},
assuming large $\delta$, $S$ and $I$, which leads to the elliptic equation
\begin{align}
- \delta \grad\cdot( \P\grad a) = S - I a. \label{cont0}
\end{align}
The system \eqref{cont1}, \eqref{cont3} is equipped with the no-flux boundary condition
\begin{align}  \label{contBC}
\nu\cdot\P\grad a = 0, \qquad \nu\cdot\grad\P_k = 0 \qquad\mbox{on } \partial\Omega,\, k=1,\dots,d,
\end{align}
where $\nu=\nu(x)$ is the outer unit normal vector on $\partial\Omega$.
The no-flux boundary condition reflects the modeling assumption that there is no flow
of auxin or the auxin transporters through the boundary of the domain.
More general boundary conditions can be considered, leading to only slight modifications in the forthcoming analysis.
Moreover, we prescribe the initial datum for the auxin transporters
\begin{align} \label{contIC}
\P_k(0,x) = \P_k^0(x)\geq 0 \qquad\mbox{for } x\in\Omega,\,k=1,\dots,d.
\end{align}

\begin{remark}
	The choice to work with the elliptic-parabolic system \eqref{cont3}, \eqref{cont0}
	instead of the parabolic-parabolic system \eqref{cont1}, \eqref{cont3}
	simplifies the mathematical analysis, since one can apply the
	so-called weak-strong lemma for the elliptic equation \eqref{cont0},
	see Lemma \ref{lem:aux} below.
	The analysis of the full parabolic-parabolic PDE system \eqref{cont1}, \eqref{cont3}
	will be the subject of a further work.
\end{remark}

\subsection{Existence of weak solutions for the system \eqref{cont3}, \eqref{cont0}}
The weak formulation of \eqref{cont0}, subject to the no-flux boundary condition \eqref{contBC},
with a test function $\phi\in C^\infty(\Omega)$ reads
\begin{align}
\delta \int_\Omega (\P\grad a)\cdot\grad \phi \d x = \int_\Omega (S - Ia) \phi \d x, \label{weak1}
\end{align}
for almost all $t>0$,
and the weak formulation of \eqref{cont3}, \eqref{contBC} with a test function $\psi\in C^\infty(\Omega)$ is
\begin{align}
\tot{}{t} \int_\Omega {\P_k} \psi \d x = - D^2 \int_\Omega \grad \P_k \cdot \grad \psi \d x
+ \int_\Omega \left( {|\partial_{x_k} a|^\kappa}{\P_k^{\kappa-\gamma}} - \tau\P_k \right) \psi \d x,  \label{weak2}
\end{align}
for almost all $t>0$.
The system is subject to the initial datum \eqref{contIC} with
\begin{align}\label{weakIC}
\P_k^0 \in L^\infty(\Omega), \quad k=1,\dots,d.
\end{align}
We assume the uniform positivity $\P_k^0 \geq \bar \P^0 >0$ almost everywhere on $\Omega$,
which prevents degeneracy of the elliptic term $\grad\cdot( \P\grad a)$ in \eqref{cont1}.
Moreover, we assume that
\begin{align} \label{SI}
S\in L^2(\Omega),\qquad I\in L^\infty(\Omega) \text{ with } I(x) \geq \bar I > 0 \text{ almost everywhere on } \Omega.
\end{align}

To prove the existence of solutions of the system \eqref{weak1}, \eqref{weak2} subject to the initial condition \eqref{weakIC}
we shall use the Schauder fixed point iteration in an appropriate function space.
We start by proving suitable a-priori estimates. % for the two decoupled equations \eqref{weak1}, \eqref{weak2}.

\begin{lemma}\label{lem:exa}
	Let $S\in L^2(\Omega)$ and $I\in L^\infty(\Omega)$ verify \eqref{SI}.
	Let the diagonal tensor $\P\in L^2(\Omega)$ be uniformly positive on $\Omega$, i.e.,
	let there be $\bar \P >0$ such that $\P_k \geq \bar \P$ almost everywhere on $\Omega$,
	for $k=1,\dots,d$.
	Then there exists a unique solution $a\in H^1(\Omega)$ of \eqref{weak1}
	and a constant $C>0$ depending only $\delta$, $\bar\P$, $S$ and $\bar I$, such that
	\begin{align}  \label{a-apriori1}
	\Norm{a}_{H^1(\Omega)} \leq C.
	\end{align}
\end{lemma}

\begin{proof}
	Let us consider a sequence of uniformly positive diagonal tensors $\P^n \in L^\infty((0,T)\times\Omega)$,
	$\P^n_k \geq \bar \P$ almost everywhere on $\Omega$ for all $n\in\N$, such that $\P^n\to\P$ in the norm topology
	of $L^2((0,T)\times\Omega)$ as $n\to\infty$. For each $n\in\N$ a unique solution $a^n\in H^1(\Omega)$ of \eqref{weak1}
	is constructed using the Lax-Milgram Theorem, see, e.g., \cite{evans}.
	The continuity of the bilinear form $B: H^1(\Omega)\times H^1(\Omega) \to \R$ associated with \eqref{weak1},
	\[
	B(a,\phi) :=  \delta \int_\Omega (\P\grad a)\cdot\grad \phi \d x - \int_\Omega (S - Ia) \phi \d x,
	\]
	follows from a straightforward application of the Cauchy-Schwart inequality.
	The coercivity of $B$ follows from
	\[
	- \int_\Omega Sa \d x \geq - \frac{1}{4\bar I} \int_\Omega S^2 \d x - \bar I \int_\Omega a^2 \d x
	\]
	and the uniform boundedness $I(x) \geq \bar I$.
	Using $\phi:=a^n$ as a test function in \eqref{weak1} gives
	\[
	\delta \int_\Omega \grad a^n \cdot\P^n\grad a^n \d x = \int_\Omega S a^n \d x - \int_\Omega I (a^n)^2 \d x,
	\]
	%	and by \eqref{SI} and the Cauchy-Schwartz inequality we have
	%	\begin{align*} %\label{est1}
	%	\delta \int_\Omega \grad a^n \cdot\P^n\grad a^n \d x + \frac{\bar I}{2} \int_\Omega (a^n)^2 \d x
	%	\leq \frac1{2\bar I} \int_\Omega S^2 \d x.
	%	\end{align*}
	%	The uniform boundedness $\P^n_k \geq \bar \P>0$ gives then
	By \eqref{SI}, the Cauchy-Schwartz inequality and the uniform boundedness $\P^n_k \geq \bar \P>0$ we have
	\begin{align}  \label{est2}
	\delta \bar \P \int_\Omega |\grad a^n|^2 \d x + \frac{\bar I}{2} \int_\Omega (a^n)^2 \d x
	\leq \frac1{2\bar I} \int_\Omega S^2 \d x
	\end{align}
	and thus a uniform bound on $a^n$ in $H^1(\Omega)$.
	
	Consequently, we can extract a subsequence converging to some $a$ weakly in $H^1(\Omega)$
	and strongly in $L^2(\Omega)$.
	Then, it is trivial to pass to the limit in \eqref{weak1}, where the term $\P^n\grad a^n$
	converges to $\P\grad a$ due to the strong convergence of $\P^n$ in $L^2(\Omega)$.
	Consequently, the limiting object $a$ verifies the weak formulation \eqref{weak1}.
	Moreover, it satisfies the a-priori estimates \eqref{est2} due to the weak lower semicontinuity
	of the respective norms.
	Uniqueness of the solution follows from \eqref{est2} and the linearity of the equation.
\end{proof}

\begin{remark}
	With a straightforward modification of its proof, we shall apply Lemma \ref{lem:exa}
	for time-dependent permeability tensors $\P\in L^\infty(0,T; L^2(\Omega))$ in the sequel.
	We then obtain the unique solution $a\in L^2(0,T; H^1(\Omega))$ satisfying
	the uniform estimate
	\begin{align}  \label{est:a}
	\Norm{a}_{L^2(0,T;H^1(\Omega))} \leq C
	\end{align}
	with $C=C(\delta, \bar\P, S, \bar I)>0$.
\end{remark}

The following Lemma is an instance of the so-called weak-strong lemma
for elliptic problems, see, e.g. \cite[Lemma 1]{HaKrMa2}.
Here we formulate it in the time-dependent setting with $a=a(t,x)$.

\begin{lemma}\label{lem:aux}
	Fix $T>0$ and let $(\P^n)_{n\in\N} \subset L^\infty(0,T; L^2(\Omega))$ be a sequence of
	diagonal tensors in $\R^{d\times d}$ such that for some $\bar\P>0$,
	$\P^n_k \geq \bar \P >0$ almost everywhere on $(0,T)\times \Omega$, $k=1,\dots,d$, $n\in\N$.
	Moreover, assume that $\P^n\to \P$ in the norm topology of $L^2((0,T)\times\Omega)$.
	Let $(a^n)_{n\in\N}$ be a sequence of weak solutions of \eqref{weak1}
	with the permeability tensors $\P^n$.
	Then $\grad a^n$ converges to $\grad a$ strongly in $L^q((0,T)\times\Omega)$ for any $q<2$,
	where $a$ is the solution of \eqref{weak1} with permeability tensor $\P$.
\end{lemma}

\begin{proof}
	Due to the uniform estimate on $a^n$ in $L^2(0,T; H^1(\Omega))$ of Lemma \ref{lem:exa},
	$a^n$ that converges weakly in $L^2(0,T; H^1(\Omega))$ to some $a$.
	%Clearly, the weak solutions of \eqref{weak} also satisfy the distributional formulation of \eqref{cont1}.
	Since $a^n\to a$ strongly in $L^2((0,T)\times \Omega)$, we can pass to the limit $n\to\infty$ in \eqref{weak1}.
	%\(  \label{a-distr}
	%   - \int_0^T \int_\Omega a \partial_t\phi \d x\d t - \int_\Omega a^0 \phi(t=0) \d x =
	%   - \delta \int_0^T \int_\Omega \P \grad a \cdot\grad\phi \d x\d t + \int_0^T \int_\Omega (S-Ia) \phi \d x\d t
	%   \qquad\mbox{for all } \phi\in C_0^\infty([0,T)\times\Omega).
	%\)
	With the uniform estimate on $\sqrt{\P^n}\grad a^n$ in $L^2((0,T)\times\Omega)$ provided by \eqref{est:a},
	the weak lower semicontinuity of the $L^2$-norm implies
	\begin{align}  \label{aux:wlsc}
	%   \int_\Omega |\grad a|^2 &\leq& \bar\P^{-1} \int_\Omega \P \grad a \cdot\grad a \d x \\
	%   &=& \bar\P^{-1} \int_\Omega |\sqrt{\P}\grad a|^2 \d x
	%   \leq \bar\P^{-1} \liminf_{n\to\infty}  \int_\Omega |\sqrt{\P^n}\grad a^n|^2 \d x < +\infty,
	\int_0^T \int_\Omega \P \grad a \cdot\grad a \d x\d t = \int_0^T \int_\Omega |\sqrt{\P}\grad a|^2 \d x\d t
	\leq \liminf_{n\to\infty}  \int_0^T \int_\Omega |\sqrt{\P^n}\grad a^n|^2 \d x\d t < +\infty
	\end{align}
	for almost all $t>0$.
	Consequently, we can use $a$ as a test function in the time-integrated version of \eqref{weak1} to obtain
	\[
	\delta \int_0^T \int_\Omega \P \grad a \cdot\grad a \d x\d t = \int_0^T \int_\Omega (S-Ia) a \d x \d t.
	% \qquad\mbox{for all } \phi\in C_0^\infty(\Omega).
	\]
	Then, using $a^N$ as a test function in \eqref{weak1} with $\P^n$, we have
	\begin{align*}
	\lim_{N\to\infty} \delta \int_0^T \int_\Omega \P^n \grad a^n \cdot\grad a^n \d x %&=	\lim_{N\to\infty}  \int_0^T \int_\Omega (S-Ia^n) a^n \d x\d t \\
	=\int_0^T \int_\Omega (S-Ia) a \d x\d t
	= \delta \int_0^T \int_\Omega \P \grad a \cdot\grad a \d x\d t.
	\end{align*}
	Consequently,
	\[
	\int_0^T \int_\Omega |\sqrt{\P}\grad a|^2 \d x\d t = \lim_{n\to\infty} \int_0^T \int_\Omega |\sqrt{\P^n}\grad a^n|^2 \d x\d t,
	\]
	so that we have the strong convergence of $\sqrt{\P^n}\grad a^n$ to $\sqrt{\P}\grad a$ in $L^2((0,T)\times\Omega)$.
	Now we write, 
	\begin{align*}
	\int_0^T \int_\Omega |\partial_{x_k} a^n - \partial_{x_k} a| \d x\d t
	&\leq \bar\P^{-1/2} \int_0^T \int_\Omega |\sqrt{\P_k}\partial_{x_k} a^n - \sqrt{\P_k}\partial_{x_k} a| \d x\d t \\
	&\leq \bar\P^{-1/2} \Norm{\grad a^n}_{L^2((0,T)\times\Omega)} \Norm{\sqrt{\P_k^n}-\sqrt{\P_k}}_{L^2((0,T)\times\Omega)} \\
	&\quad + \bar\P^{-1/2} \int_0^T \int_\Omega |\sqrt{\P^n_k}\partial_{x_k} a^n - \sqrt{\P_k}\partial_{x_k} a| \d x\d t,
	\end{align*}
	for $k=1,\dots,d$, and the first term of the right-hand side converges to zero due to the assumed strong convergence of $\P^n$ in $L^2((0,T)\times\Omega)$,
	while the second term does so due to the strong convergence of $\sqrt{\P^n}\grad a^n$.
	Thus, we have the strong convergence of $\grad a^n$ to $\grad a$ in $L^1((0,T)\times\Omega)$.
	Since $\grad a^n$ is also uniformly bounded in $L^2((0,T)\times\Omega)$, a simple consequence of the interpolation
	inequality \cite[Chapter 1]{roubicek} implies strong convergence in $L^q((0,T)\times\Omega)$ for $q<2$.
\end{proof}

\begin{lemma}\label{lem:exP}
	Fix $T>0$ and let $\grad a\in L^2((0,T)\times\Omega)$.
	Let $\kappa > \gamma$ and, 
	\begin{align} \label{kappa-d}
	%\begin{aligned}
	\kappa < 2 \quad\mbox{for } d\in \{1,2\},\qquad
	\kappa \leq \frac{\gamma+5}{4} \quad\mbox{for } d=3,
	%\end{aligned}
	\end{align}
	depending on the space dimension $d$.	Then there exists a unique solution
	\[
	\P_k\in L^2(0,T; H^1(\Omega)) \cap L^\infty(0,T; L^2(\Omega)) \cap C([0,T); H^{-1}(\Omega)), \qquad k=1,\dots,d,
	\]
	of  \eqref{weak2} subject to the initial datum \eqref{weakIC} with $\P_k^0 \geq \bar \P^0 >0$ almost everywhere on $\Omega$.
	Moreover, the solution stays uniformly bounded away from zero on $(0,T)\times\Omega$, i.e.,
	there exists $\bar\P>0$ depending on $\bar \P^0$, $T$, $D^2$ and $\tau$, but independent of $a$, such that
	\begin{align}  \label{a-priori-positive}
	\P_k \geq \bar\P >0 \qquad\mbox{almost everywhere on } (0,T)\times\Omega.
	\end{align}
	Moreover, there exists a constant $K_0>0$ independent of $\P$ and $a$ such that
	\begin{align} \label{P-apriori1}
	\Norm{\P_k}^2_{L^\infty(0,T; L^2(\Omega))} \leq \Norm{\P_k^0}_{L^2(\Omega)}^2 + K_0 \Norm{\partial_{x_k} a}_{L^2((0,T)\times \Omega)}^2
	\end{align}
	and, 	for $k=1,\dots,d$,
	\begin{align} \label{P-apriori2}
	\Norm{\grad \P_k}^2_{L^2(0,T; L^2(\Omega))} \leq \Norm{\P_k^0}_{L^2(\Omega)}^2 + K_0 \Norm{\partial_{x_k} a}_{L^2((0,T)\times \Omega)}^2.
	\end{align}
\end{lemma}

\begin{remark}
	Observe that the necessary condition
	for the mutual validity of the assumptions $\kappa > \gamma$ and \eqref{kappa-d}
	is $\gamma, \kappa<2$ for $d\in\{1,2\}$ and $\gamma, \kappa\leq 5/3$ for $d=3$.
\end{remark}

\begin{proof}
	%$\frac{|q_k|^\kappa}{\P_k^{\gamma+1}} \P_k = |\partial_{x_k} a|^\kappa \P_k^{\kappa-\gamma} \sgn(\P_k)$
	%so that for $\kappa\geq\gamma>0$ and $\grad a\in L^2((0,T)\times\Omega)$ the classical parabolic theory provides the local in time existence of  a weak solution of \eqref{weak2}.
	Let us fix $k \in \{1,\dots,d\}$ and use $\psi:=\P_k$ as a test function in \eqref{weak2},
	\begin{align} \label{est21}
	\frac12 \tot{}{t} \int_\Omega \P_k^2 \d x = -D^2 \int_\Omega |\grad\P_k|^2 \d x
	+ \int_\Omega |\partial_{x_k} a|^\kappa \P_k^{\kappa - \gamma +1} \d x - \tau \int_\Omega \P_k^2 \d x,
	\end{align}
	where we used the identity $q_k = \P_k \partial_{x_k} a$.
	Using the H\"older inequality with exponents $p$ and $p'$, $\frac1p + \frac1{p'} = 1$, we have
	\begin{align}  \label{Holder}
	\int_\Omega |\partial_{x_k} a|^\kappa \P_k^{\kappa - \gamma +1} \d x \leq
	C_\eps \int_\Omega |\partial_{x_k} a|^{\kappa p} \d x + \eps \int_\Omega |\P_k|^{(\kappa - \gamma +1)p'} \d x
	\end{align}
	for $\eps>0$ and a suitable constant $C_\eps$.
	Due to the assumed $L^2$-integrability of $\partial_{x_k} a$, we choose $\kappa p = 2$, so that $p' = \frac{2}{2-\kappa}$.
	%with the convention $p'=\infty$ if $\kappa=2$, in which case we replace \eqref{Holder} by
	%\(  \label{noHolder}
	%  \int_\Omega |\partial_{x_k} a|^\kappa \P_k^{\kappa - \gamma +1} \d x \leq
	%     \Norm{\P_k}_{L^\infty(\Omega)} \int_\Omega |\partial_{x_k} a|^2 \d x.
	%\)
	Denote $\alpha := (\kappa - \gamma +1)p'$ and observe that $\alpha>0$ due to the assumption $\kappa > \gamma$.
	Let us  distinguish the following two cases:
	If $\alpha\leq 2$, then by the H\"older inequality we have
	\[
	\int_\Omega |\P_k|^\alpha \d x \leq C_\Omega \int_\Omega |\P_k|^2 \d x,
	\]
	so that  \eqref{est21} and \eqref{Holder} imply
	\[
	\frac12 \tot{}{t} \int_\Omega \P_k^2 \d x \leq - D^2 \int_\Omega |\grad\P_k|^2 \d x
	+ C_\eps \int_\Omega |\partial_{x_k} a|^2 \d x - (\tau - \eps C_\Omega) \int_\Omega \P_k^2 \d x,
	\]
	and choosing $\eps>0$ such that $\tau - \eps C_\Omega >0$ directly implies the a-priori estimates \eqref{P-apriori1} and \eqref{P-apriori2}.
	On the other hand, if $\alpha >2$, we apply the Sobolev inequality \cite{evans}
	\[
	\int_\Omega |\P_k|^{\alpha} \d x \leq C_S \left( \int_\Omega |\grad \P_k|^2 \d x + \int_\Omega |\P_k|^2 \d x \right)
	\]
	with $C_S=C_S(\Omega)$ the Sobolev constant. Depending on the space dimension, we have:
	\begin{itemize}
		%\item
		%For $d=1$,
		%\[
		%  \Norm{\P_k}_{L^\infty(\Omega)} \leq C_S \left( \int_\Omega |\grad \P_k|^2 \d x + \int_\Omega |\P_k|^2 \d x \right)
		%\]
		%and, consequently, we may choose $\alpha=\infty$, i.e., $p=1$ and we admit $\kappa\leq 2$ in \eqref{Holder}.
		\item
		For $d\in\{1,2\}$,
		\begin{align}  \label{SobIneq}
		\Norm{\P_k}_{L^\alpha(\Omega)} \leq C_S \left( \int_\Omega |\grad \P_k|^2 \d x + \int_\Omega |\P_k|^2 \d x \right)
		\end{align}
		for any $\alpha<\infty$, i.e., we admit any $p>1$ and, consequently, $\kappa < 2$.
		\item
		For $d=3$ we have \eqref{SobIneq} for $\alpha\leq 6$, i.e., we need
		$(\kappa-\gamma+1)p' = \frac{2(\kappa-\gamma+1)}{2-\kappa} \leq 6$,
		which gives the condition $\kappa \leq \frac{\gamma+5}{4}$.
		%Note that the assumption $\kappa>\gamma$ implies that $\kappa < \frac{5}{3}$.
	\end{itemize}
	Consequently, we have
	\[
	\frac12 \tot{}{t} \int_\Omega \P_k^2 \d x \leq - (D^2 - \eps C_S) \int_\Omega |\grad\P_k|^2 \d x
	+ C_\eps \int_\Omega |\partial_{x_k} a|^2 \d x - (\tau - \eps C_S) \int_\Omega \P_k^2 \d x,
	\]
	and choosing $\eps>0$ such that $\eps C_S < \min\{D^2, \tau\}$ directly implies the a-priori estimates \eqref{P-apriori1} and \eqref{P-apriori2}.
	The uniform positivity \eqref{a-priori-positive} follows from the fact that solutions $u=u(t,x)$ of the linear parabolic equation
	$\part{u}{t} = D^2\laplace u - \tau u$
	are subsolutions to \eqref{cont2}, and they remain uniformly positive on bounded time intervals for uniformly positive initial data,
	see, e.g., \cite{evans}.
	
	Finally, note that we have the identity (in distributional sense)
	\[
	\part{\P_k}{t} = D^2\laplace \P_k + |\partial_{x_k} a|^\kappa \P_k^{\kappa-\gamma} - \tau\P_k.
	\]
	An easy calculation reveals that, for the aforementioned range of $\kappa$ and $\gamma$,
	\[
	|\partial_{x_k} a|^\kappa \P_k^{\kappa-\gamma} \in L^1(0,T; L^{6/5}(\Omega)) \subset L^1(0,T; H^{-1}(\Omega)),
	\]
	implying $\part{\P_k}{t} \in L^1(0,T; H^{-1}(\Omega))$,
	so that $\P_k \in C([0,T); H^{-1}(\Omega))$, see, e.g., \cite[Chapter 7]{roubicek}.
\end{proof}

\def\calB{\mathcal{B}}
\begin{theorem}\label{thm:ex_cont}
	Fix $T>0$ and let $\kappa > \gamma$, and, in dependence of the space dimension $d$,
	\begin{align} \label{kappa-d2}
	\begin{aligned}
	\kappa &< \frac{\gamma+4}{3} \quad\mbox{for } d\in \{1,2\},\qquad
	\kappa &< \frac{\gamma+5}{4} \quad\mbox{for } d=3.
	\end{aligned}
	\end{align}
	Then the system \eqref{weak1}--\eqref{weak2} subject to the initial datum \eqref{weakIC}
	with $\P_k^0 \geq \bar \P^0 >0$ almost everywhere on $\Omega$
	admits a weak solution $(\P,a)$ on $(0,T)$ such that
	\begin{align}  \label{weak_sol}
	\begin{aligned}
	\P_k &\in L^\infty(0,T; L^2(\Omega)) \cap L^2(0,T; H^1(\Omega)) \cap C([0,T); H^{-1}(\Omega)), \\ %,\qquad k=1,\dots,d,\\
	a &\in L^\infty(0,T; L^2(\Omega)) \cap L^2(0,T; H^1(\Omega)) \cap C([0,T); W^{-1,4/3}(\Omega)).
	\end{aligned}
	\end{align}
	%\comment{Do we have continuity in time in a "better" topology?}
\end{theorem}

\begin{proof}
	We construct a solution using the Schauder fix-point theorem on the set
	\begin{align*}
	\calB_T := \Bigl\{ \P \in (L^\infty(0,T; L^2(\Omega)))^{d\times d}_\mathrm{diag};\;
	\Norm{\P_k}^2_{L^\infty(0,T; L^2(\Omega))} \leq \Norm{\P_k^0}_{L^2(\Omega)}^2 + K_0 B_T^2,\;  \\
	%     \Norm{\grad u}^2_{L^2((0,T)\times\Omega)} \leq \Norm{\P_k^0}_{L^2(\Omega)}^2 + C_0 B_T^2
	\P_k \geq \bar\P\mbox{ almost everywhere on } (0,T)\times\Omega,\; k=1,\dots,d
	\Bigr\}.
	\end{align*}
	Here $(L^\infty(0,T; L^2(\Omega)))^{d\times d}_\mathrm{diag}$ denotes the space
	of diagonal $d\times d$-tensors with entries in $L^\infty(0,T; L^2(\Omega))$, and
	$K_0$ and $\bar\P$ are the constant defined in Lemma \ref{lem:exP};
	note that they depend only on $\bar \P^0$, $T$, and the parameters $\kappa$, $\gamma$, $D^2$ and $\tau$.
	Moreover, we denoted
	\[
	B^2_T := \frac1{2\delta\bar\P} \left( (Te^T+1) \Norm{a^0}_{L^2(\Omega)}^2 + Te^T \Norm{S}_{L^2(\Omega)}^2 \right).
	\]
	The set $\calB_T$ shall be equipped with the norm topology of $L^2((0,T)\times\Omega)$.
	Obviously, $\calB_T$ is nonempty, convex and closed.
	We define the mapping $\Phi: \calB_T \to L^\infty(0,T; L^2(\Omega))$,
	\( %\label{Phi}
	\Phi: \P\in\calB_T \mapsto \tilde \P,
	\)
	where given $\P\in\calB_T$ we construct $a$ the unique weak solution of \eqref{weak1} by Lemma \ref{lem:exa},
	and, subsequently, construct $\tilde\P$ as the unique weak solution of \eqref{weak2} by Lemma \ref{lem:exP}.
	Clearly, due to the a-priori estimates \eqref{a-apriori1} and \eqref{P-apriori1}, $\tilde \P\in\calB_T$.
	
	To prove the continuity of the mapping $\Phi$, let us consider a sequence $(\P_n)_{n\in\N} \subset\calB_T$,
	converging to $\P\in\calB_T$ in the norm topology of $L^2((0,T)\times\Omega)$.
	Denote $(a_n)_{n\in\N}$ and, resp., $a$, the solutions of \eqref{weak1}
	corresponding to $\P_n$ and, resp., $\P$. 
	Then, due to Lemma \ref{lem:aux}, $\grad a_n$ converges to $\grad a$
	in the norm topology of $L^q((0,T)\times\Omega)$ for any $q<2$.
	Let $\tilde \P_n:=\Phi(\P_n)$ and $\tilde \P:=\Phi(\P)$.
	Due to Lemma \ref{lem:exP} and the Aubin-Lions theorem, a subsequence of $\tilde \P_n$
	converges strongly to some $\tilde\P^*$ in $L^2(0,T; L^q(\Omega))$ with $q<\infty$ if $d\in\{1,2\}$
	and $q=6$ if $d=3$. The limit passage $n\to\infty$ in \eqref{weak2} is trivial for the linear terms.
	For the term $|\partial_{x_k} a_n|^\kappa \tilde\P_n^{\kappa-\gamma}$ we observe that,
	due to Lemma \ref{lem:aux}, the term $|\partial_{x_k} a_n|^\kappa$ converges to $|\partial_{x_k} a|^\kappa$
	in the norm topology of $L^q((0,T)\times\Omega)$ for $q<{2}/{\kappa}$. Moreover:
	\begin{itemize}
		\item
		For $d\in\{1,2\}$, the interpolation inequality between $L^\infty(0,T; L^2(\Omega))$ and $L^2(0,T; L^q(\Omega))$ with $q<\infty$
		implies that $\tilde\P^n$ is uniformly bounded, and thus converges, in the norm topology of $L^q((0,T)\times\Omega)$ for $q<4$.
		Consequently, since $\kappa<2$, the product $|\partial_{x_k} a_n|^\kappa \tilde\P_n^{\kappa-\gamma}$ converges strongly in (at least) $L^1((0,T)\times\Omega)$
		to $|\partial_{x_k} a|^\kappa (\tilde\P^*_n)^{\kappa-\gamma}$ if
		$	\frac{\kappa}{2} + \frac{\kappa-\gamma}{4} < 1$,
		which is equivalent to $\kappa<\frac{\gamma+4}{3}$.
		\item
		For $d=3$ the interpolation inequality between $L^\infty(0,T; L^2(\Omega)$ and $L^2(0,T; L^6(\Omega))$
		implies that $\tilde\P^n$ is uniformly bounded in the norm topology of $L^{10/3}((0,T)\times\Omega)$.
		Then the sufficient condition for $L^1$-convergence of the product $|\partial_{x_k} a_n|^\kappa \tilde\P_n^{\kappa-\gamma}$ reads
		$\frac{\kappa}{2} + \frac{3(\kappa-\gamma)}{10} < 1$,
		which is equivalent to $\kappa<\frac{10+3\gamma}{8}$. This condition is weaker than \eqref{kappa-d2}.
	\end{itemize}
	By the uniqueness of solutions of \eqref{weak1}, we conclude that $\tilde\P^* = \tilde\P$,
	i.e., the mapping $\Phi$ is continuous on $\calB_T$ with respect to the norm topology of $L^2((0,T)\times\Omega)$.
	
	To prove the compactness of the mapping $\Phi$, we employ the Aubin-Lions lemma \cite{Aubin}.
	Let us again consider a sequence $(\P_n)_{n\in\N}\subset\calB_T$ and denote $\tilde \P_n := \Phi(\P_n)$.
	Due to the a-priori estimates \eqref{a-apriori1} and \eqref{P-apriori1}, \eqref{P-apriori2},
	the sequence $\tilde \P_n$ is bounded in $L^\infty(0,T; L^2(\Omega))$
	and in $L^2(0,T; H^1(\Omega))$. Moreover, $\partial_t \tilde \P_n$ is bounded
	in $L^1(0,T; H^{-1}(\Omega))$.
	Then, since $H^1(\Omega)$ is compactly embedded into $L^2(\Omega)$ and $L^2(\Omega) \subset H^{-1}(\Omega)$,
	the Aubin-Lions theorem provides the relative compactness of the sequence $\tilde \P_n$ with respect to the norm topology of $L^2((0,T)\times\Omega))$.
	Consequently, the Schauder fix-point theorem provides a solution $(\P,a)$ of the system
	\eqref{weak1}--\eqref{weakIC}, satisfying \eqref{weak_sol}.
\end{proof}

\begin{remark}
	For the case $\kappa=\gamma=2$ the system \eqref{cont3} simplifies to
	\(
	\part{\P_k}{t} = D^2\laplace \P_k + {(\partial_{x_k} a)^2} - \tau\P_k. \label{cont4}
	\)
	Then, \eqref{cont0}, \eqref{cont4} is similar to the system studied in \cite{HaKrMa} and \cite{HaKrMa2},
	the main difference being that the permeability tensor in the elliptic equation
	is of the form $rI + X$ in  \cite{HaKrMa},  \cite{HaKrMa2}, where $r>0$ is a constant.
	The significant property of \eqref{cont0}, \eqref{cont4} is its energy-dissipation structure.
	Indeed, defining
	\[
	\mathcal{E}[\P] := \frac{D^2}{2} \sum_{k=1}^d \int_\Omega |\grad\P_k|^2  \d x + \int_\Omega \grad a\cdot\P\grad a \d x
	+ \tau \sum_{k=1}^d \int_\Omega \P_k^2 \d x,
	\]
	where $a=a[\P]$ is the unique weak solution of \eqref{cont0},
	a simple calculation (see \cite[Lemma 3]{HaKrMa}) reveals that, %formally,
	\[
	\tot{}{t} \mathcal{E}[\P] = - \sum_{k=1}^d \int_\Omega \left( \part{\P_k}{t} \right)^2 \d x
	\]
	along the solutions of \eqref{cont0}, \eqref{cont4}.
	The energy dissipation naturally provides uniform a-priori estimates on $\P$ and $a$ in the energy space.
	However, these still do not allow us to extend the validity of Theorem \ref{thm:ex_cont}
	to $\kappa=\gamma=2$. The problem is that in the proof of continuity of the fix-point mapping $\Phi$,
	it is not clear how to pass to the (weak) limit in the sequence $(\partial_{x_k} a)^2$.
	Note that Lemma \ref{lem:aux} only provides (strong) convergence of $\partial_{x_k} a$ in $L^q((0,T)\times\Omega)$
	with $q<2$.
\end{remark}

\def\calA{\mathcal{A}}
\begin{remark}[Steady states of the system \eqref{cont3}, \eqref{cont0} with $D^2=0$]
	The steady states of the system \eqref{cont3}, \eqref{cont0} with $D^2=0$ satisfy, in the weak sense,
	\begin{align}
	\delta \grad\cdot( \P\grad a) + S - Ia = 0, \label{steady1}\\
	|\partial_{x_k} a|^\kappa \P_k^{\kappa-\gamma} - \tau \P_k = 0, \label{steady2}
	\end{align}
	for $k=1,\dots,d$, with $q_k = \P_k\partial_{x_k} a$.
	For $\kappa > \gamma > 0$, \eqref{steady2} implies that there exist measurable sets
	$\calA_k\subset\Omega$, $k=1,\dots,d$, such that
	\[
	\P_k =  \left( \frac{|\partial_{x_k} a|^\kappa}{\tau} \right)^\frac{1}{\gamma-\kappa+1} \chi_{k},
	\]
	where $\chi_{k}=\chi_{k}(x)$ is the characteristic function of $\calA_k$.
	Inserting this into \eqref{steady1}, we obtain
	\begin{align} \label{steadyA}
	- \delta \tau^\frac{1}{\kappa-\gamma-1} \sum_{k=1}^d \partial_{x_k} \left( \chi_{k} |\partial_{x_k} a|^\frac{\kappa}{\gamma-\kappa+1} \partial_{x_k} a \right)
	= S - Ia.
	\end{align}
	Due to the presence of the characteristic functions $\chi_{k}$, this is a strongly degenerate elliptic equation,
	rendering its analysis a very challenging task, which we leave for a future work.
	Let us only note that the degeneracy in \eqref{steadyA} induces strong nonuniqueness of its solutions.
	Consequently, it is necessary to equip \eqref{steadyA} with suitable selection criteria in order
	to obtain unique solutions. This is to be done through further modeling inputs.
	For $\kappa = \gamma > 0$, contrarily, \eqref{steady2} gives $\P_k = \tau^{-1} |\partial_{x_k} a|^\kappa$,
	and \eqref{steady1} reads
	\begin{align} \label{steadyB}
	- \delta \tau^{-1} \sum_{k=1}^d \partial_{x_k} \left( |\partial_{x_k} a|^\kappa \partial_{x_k} a \right) = S - Ia.
	\end{align}
	Equipped with the no-flux boundary condition \eqref{contBC}, its weak formulation reads
	\begin{align}  \label{steadyBw}
	\delta \tau^{-1} \sum_{k=1}^d \int_\Omega |\partial_{x_k} a|^\kappa (\partial_{x_k} a) (\partial_{x_k} \psi) \d x
	+ \int_\Omega (a-S)\psi \d x = 0
	\end{align}
	for all test functions $\psi\in C^\infty(\Omega)$.
	Weak solutions $a\in W^{1,\kappa+2}(\Omega)$ of \eqref{steadyBw} are constructed as the global minima of the functional
	$\mathcal{F}:W^{1,\kappa+2} \to \mathbb{R}$,
	\[
	\mathcal{F}[a] := \frac{\delta \tau^{-1}}{\kappa+2} \sum_{k=1}^d
	\int_\Omega |\partial_{x_k} a|^{\kappa+2} \d x
	+ \frac{1}{2}\int_\Omega a^2 \d x - \int_\Omega S a \d x.
	\]
	Obviously, for $\kappa>0$ the functional is uniformly convex.
	Moreover, a straightforward application of the Cauchy-Schwartz inequality implies boundedness below
	and coercivity of $\mathcal{F}$ with respect to the norm of $W^{1,\kappa+2}(\Omega)$.
	Then the classical theory (see, e.g., \cite{evans}) provides the existence of a unique minimizer $a\in W^{1,\kappa+2}(\Omega)$ of $\mathcal{F}$,
	which is the unique solution of the corresponding Euler-Lagrange equation \eqref{steadyBw}.
\end{remark}

%\section*{Acknowledgments}
%LMK was supported by the UK Engineering and Physical Sciences Research Council (EPSRC) grant
%EP/L016516/1 and the German National Academic Foundation (Studienstiftung des Deutschen Volkes). 

\section{Conclusion}
In this paper, we proposed a new dynamic modelling framework for leaf venation, which is not dependent on polar localisation of auxin transporters, i.e.\ the transport capacity across a cell wall does not have to be asymmetric. Given that it is still an open question how you get leaf veins, also in the absence of PIN-based transport activity, we argue that the current work is of interest since it is the first model, to our knowledge, trying to address this question. Due to its new description of possible mechanisms in leaf venation, our model is of interest to the modelling community. Our work can be regarded as a general modelling framework for auxin transport, which can be equipped or extended
with various biologically relevant features that would then produce experimentally verifiable hypotheses.
The main advantage is the rather simple form of the model, allowing a rigorous mathematical analysis,
which is one of the main aims of our paper. Moreover, it facilitates the derivation of a continuum limit,
which can capture network growth and is expected to exhibit a much richer patterning capacity,
bearing again potential for delivering testable hypotheses. The analytical and numerical study of the continuum model is currently a work in progress.

\section*{Data Accessibility}
The data set containing the MATLAB code necessary to reproduce the computational results is available at the DOI link \url{https://doi.org/10.17863/CAM.40619}.

\section*{Acknowledgments}
HJ is supported by the Gatsby Charitable Foundation (grant GAT3395-PR4). LMK is supported by the EPSRC grant EP/L016516/1 and the German National Academic Foundation.
	
\bibliographystyle{plain}
\bibliography{references}

\end{document}